\documentclass[final]{siamart171218}

\usepackage{lipsum}
\usepackage{amsfonts}
\usepackage{graphicx}
\usepackage{epstopdf}
\usepackage{algorithm}
\usepackage[noend]{algpseudocode}
\usepackage{nccmath}
\ifpdf
  \DeclareGraphicsExtensions{.eps,.pdf,.png,.jpg}
\else
  \DeclareGraphicsExtensions{.eps}
\fi

\newcommand{\abs}[1]{\left\vert#1\right\vert}
\newcommand{\norm}[1]{\left\lVert#1\right\rVert}
\newcommand{\Prob}{\mathbb{P}}
\newcommand{\X}{\mathbf{X}}

\newcommand{\x}{\mathbf{x}}

\newcommand{\y}{\mathbf{y}}

\renewcommand{\v}{\mathbf{v}_{s,k-1}}

\newcommand{\z}{\mathbf{z}}

\renewcommand{\P}{\mathbf{P}_{s,k}}

\newcommand{\I}{\mathbf{I}}

\newcommand{\Expectation}{\mathbb{E}}

\newcommand{\convas}{\overset{a.s.}{\longrightarrow}}
\newcommand{\convp}{\overset{\Prob}{\longrightarrow}}
\newcommand{\T}{\mathsf{T}}

\newcommand{\PPT}{\mathbf{P}_k\mathbf{P}_k^{\T}}
\newcommand{\CondExp}[1]{\Expectation\left[#1 \mid \mathcal{F}_{s,k} \right]}

\newcommand{\reals}{\mathbb{R}}
\newtheorem{remark}{Remark}
\algnewcommand{\Initialize}[1]{%
	\State \textbf{Initialize:}
	\Statex \hspace*{\algorithmicindent}\parbox[t]{.8\linewidth}{\raggedright #1}
}

\algnewcommand{\Inputs}[1]{%
	\State \textbf{Inputs:}
	\Statex \hspace*{\algorithmicindent}\parbox[t]{.8\linewidth}{\raggedright #1}
}
\algnewcommand{\Outputs}[1]{%
	\State \textbf{Outputs:}
	\Statex \hspace*{\algorithmicindent}\parbox[t]{.8\linewidth}{\raggedright #1}
}

\theorembodyfont{\upshape}
\newsiamthm{assumptions}{Assumptions}

\newsiamremark{hypothesis}{Hypothesis}
\crefname{hypothesis}{Hypothesis}{Hypotheses}
\newsiamthm{claim}{Claim}

\headers{Stochastic Subspace Descent}{D. Kozak, S. Becker, A. Doostan, L. Tenorio}

\title{Stochastic Subspace Descent\thanks{Submitted to the editors DATE.
\funding{AD acknowledges funding by the US Department of Energy’s Office of Science Advanced Scientific Computing Research, Award DE-SC0006402 and NSF Grant CMMI-145460. LT acknowledges funding by NSF grant DMS-1723005. SB acknowledges funding by NSF grant DMS-1819251.}}}

\author{ David Kozak\thanks{Department of Applied Mathematics and Statistics, Colorado School of Mines, Golden, CO}
\and
Stephen Becker\thanks{Department of Applied Mathematics, University of Colorado, Boulder, CO}
\and Alireza Doostan\thanks{Aerospace Engineering Sciences Department, University of Colorado, Boulder, CO}
\and Luis Tenorio\footnotemark[4]}

\usepackage{amsopn}

\newcommand{\defeq}{\stackrel{\text{\tiny def}}{=}}

\newcommand{\Ib}{\textbf{I}}
\newcommand{\fb}{\textbf{f}}
\newcommand{\bP}{\textbf{P}}
\newcommand{\Sigmab}{\mathbf{\Sigma}}
\newcommand{\zero}{\mathbf{0}}
\newcommand{\bA}{\textbf{A}}
\newcommand{\R}{\mathbb{R}}
\DeclareMathOperator*{\argmin}{argmin}

\ifpdf
\hypersetup{
  pdftitle={Stochastic Subspace Descent},
  pdfauthor={D. Kozak, S. Becker, A. Doostan, L. Tenorio}
}
\fi

\usepackage{enumitem}

\usepackage{bm,psfrag}

\begin{document}

\maketitle

\begin{abstract}
We present two stochastic descent algorithms that apply to unconstrained optimization and are particularly efficient when the objective function is slow to evaluate and gradients are not easily obtained, as in some PDE-constrained optimization and machine learning problems.
The basic algorithm projects the gradient onto a random subspace at each iteration, similar to coordinate descent but without restricting directional derivatives to be along the axes. This algorithm is previously known but we provide new analysis. We also extend the popular SVRG method to this framework without requiring finite-sum objective functions.
We provide proofs of convergence under various convexity assumptions and show favorable results when compared to gradient descent and BFGS on non-convex problems from the machine learning and shape optimization literature.
\end{abstract}



\section{Introduction}

We consider optimization problems of the form 
\begin{equation}\label{eq: minimize f}
 \min_{\x \in \reals^d} f(\x) \quad \text{or} \qquad   \arg\min_{\x \in \reals^d} f(\x),
\end{equation}
where $f : \reals^d \to \reals$ has $\lambda$-Lipschitz derivative. We also consider additional restrictions such as convexity or $\gamma$-strong convexity of $f$, restrictions that will be made clear as they are required. In the most basic form discussed in Section \ref{sect: SSD}, we project the gradient onto a random $\ell$-dimensional subspace and descend along that subspace as in gradient descent. We use the following iteration scheme
\begin{equation}\label{eq: iterations}
    \x_{k+1} = \x_k - \alpha \PPT \nabla f(\x_k),
\end{equation}
where $\mathbf{P} \in \reals^{d \times \ell}$, for some $\ell < d$, is a random matrix with the properties $\Expectation~\mathbf{P}\mathbf{P}^\top  = \I_d$ and $\mathbf{P}^\top \mathbf{P} = (d/\ell) \I_{\ell}$. This algorithm has convergence properties consistent with gradient descent (see Theorem \ref{thm:convergence}) and allows us to leverage techniques from the substantial literature surrounding stochastic gradient descent (see Theorem \ref{thm:VarianceReducedRandomGradient}). Note that \eqref{eq: iterations} with a diagonal $\mathbf{P}\mathbf{P}^\top $ reduces to randomized block-coordinate descent \cite{bertsekas1999nonlinear}. 

A particular case of \eqref{eq: minimize f} is Empirical Risk Minimization (ERM) commonly used in machine learning, where $ f(\x)=(1/n)\sum_{i=1}^n f_i(\x)$ and $n$ is typically very large. Hence this problem is amenable to iterative stochastic methods that approximate $f(\x)$ using $S$ randomly sampled observations, $(i_s)_{s=1}^{S} \in \{1, \ldots, n\}$, at each iteration with $(1/S)\sum_{s=1}^S f_{i_s}(\x)$ where $S\ll n$. While the methods we discuss do not require this type of finite-sum structure, they can be used for such problems. In fact, the canonical parameter estimation problem in machine learning is one application we envision for our work. 

There are important classes of functions that do not fit into the ERM framework and therefore do not benefit from stochastic algorithms tailored to ERM. Partial Differential Equation (PDE) constrained optimization is one such example, and except in special circumstances (such as \cite{haber2012effective}), a stochastic approach leveraging the ERM structure (such as stochastic gradient descent and its variants) does not provide any benefits. This is because in PDE-constrained optimization the cost of evaluating each $\nabla f_i(\x)$ is often identical to the cost of evaluating $\nabla f(\x)$. A primary goal of this paper is to use variations on the iterative scheme \eqref{eq: iterations} to adapt recent developments from the machine learning literature so that stochastic methods can be used effectively in the PDE-constrained optimization setting. Along the way we present results that may be useful to the machine learning community. An example that comes directly from the machine learning literature is Gaussian processes, specifically the sparse Gaussian process framework of \cite{snelson2006sparse, titsias2009variational}.

\paragraph{PDE-constrained optimization}
Partial differential equations are frequently used to model physical phenomena. Successful application of PDEs to modeling is contingent upon appropriate discretization and parameter estimation. Parameter estimation in this setting arises in optimal control, or whenever the parameters of the PDE are unknown, as in inverse problems. Algorithmic and hardware advances for PDE-constrained optimization have allowed for previously impossible modeling capabilities. Examples include fluid dynamics models with millions of parameters for tracking atmospheric contaminants \cite{Flath2011}, modeling the flow of the Antarctic ice sheet \cite{isaac2015scalable, petra2014computational}, parameter estimation in seismic inversion \cite{Abacioglu2001, bui2013computational}, groundwater hydrology \cite{WRCR:WRCR23173}, experimental design \cite{horesh2010optimal, haber2012numerical}, and atmospheric remote sensing \cite{LikelihoodInformedDimensionReduction_Alessio}. 

\paragraph{Gaussian processes} 
Gaussian processes are an important class of stochastic processes useful for modeling when there are correlations in input space. Here we use them to model an unknown function in the context of regression. Results from functional analysis allow for the flexibility of working in function space using only machinery from finite dimensional linear algebra. However, the applications of Gaussian processes are somewhat hamstrung in many modern settings because their time complexity scales as $\mathcal{O}(n^3)$ and their storage as $\mathcal{O}(n^2)$. One recourse is to approximate the Gaussian process, allowing time complexity to be reduced to $\mathcal{O}(nm^2)$ with storage requirements of $\mathcal{O}(nm)$ where $m \ll n$ is the number of points used in lieu of the full data set. Methods have been developed \cite{snelson2006sparse, titsias2009variational} which place these $m$ inducing points, also called landmark points, along the domain at points different from the original inputs; optimal placement of these points is a continuous optimization problem with dimension at least the number of inducing points to be placed. Such a framework places the burden on the optimization procedure as improperly placed points may result in poor approximations.
\subsection{Related work}
Despite being among the easiest to understand and oldest variants of gradient descent, subspace methods (by far the most common of which is coordinate descent) have been conspicuously absent from the optimization literature. Such under-representation may be partially attributed to the fact that until recently the performance of these types of algorithms did not justify the difficult theoretical work.

\paragraph{Coordinate descent schemes}
The simplest variant of subspace descent is a deterministic method that cycles over the coordinates. This method is popular because many problems have structure that makes a coordinate update very cheap. Convergence results for coordinate descent requires challenging analysis and the class of functions for which it converges is restricted; indeed, \cite{warga1963minimizing, powell1973search} provide simple examples for which the method fails to converge while simpler-to-analyze methods such as gradient descent will succeed. 
Choosing the coordinates in an appropriate, random manner can lead to results on par with gradient descent \cite{nesterov2012efficiency,richtarik2014iteration}. More recently, much emphasis has been placed on accelerating coordinate descent methods \cite{allen2016even,hanzely2018accelerated}, but all  these improvements require knowledge of the Lipschitz constants of the partial derivatives of the functions and/or special structure in the function to make updates inexpensive and to choose the sampling scheme. See \cite{wright2015coordinate} for a survey of recent results. 



\paragraph{Optimization without derivatives}\label{sec:DFO}
Our methods use directions $\bP^\T\nabla f(\x)$, where $\bP$ is $d\times \ell$ with $\ell \ll d$, which is equivalent to taking $\ell$ directional derivatives of $f$, and can therefore be estimated by forward- or centered-differences in $\ell+1$ or $2\ell$ function evaluations, respectively. We assume the error due to the finite difference approximation is negligible, which is realistic for low-precision optimization. Recent analysis has shown that many first-order methods (e.g., gradient descent, SVRG) maintain asymptotic convergence properties \cite{liu2018zeroth} when the gradients are estimated via finite differences. Because our methods only evaluate $f(\x)$ and not $\nabla f(\x)$ directly, they are in the class of derivative-free optimization (DFO) and ``zeroth-order'' optimization.  

To be clear, when $\nabla f(\x)$ is readily available, zeroth-order optimization methods are not competitive with first-order and second-order methods (e.g., Newton's method). For example, if $f(\x) = \|\bA\x-\textbf{b}\|^2$ 
then evaluating $f(\x)$ and evaluating $\nabla f(\x) = 2\bA^T(\bA\x-\textbf{b})$ have nearly the same computational cost, namely $\mathcal{O}(dm)$.  In fact, such a statement is true regardless of the structure of $f$: by using reverse-mode automatic differentiation (AD), one can theoretically evaluate $\nabla f(\x)$ in about four-times the cost of evaluating $f(\x)$, regardless of the dimension $d$ \cite{griewank2008evaluating}. In the context of PDE-constrained optimization, the popular adjoint-state method, which is a form of AD applied to either the continuous or discretized PDE, also evaluates $\nabla f(\x)$ in time independent of the dimension.
However, there are many situations when AD and the adjoint-state method are inefficient or not applicable. 
Finding the adjoint equation requires a careful derivation (which depends on the PDE as well as on initial and boundary conditions), and then a numerical method must be implemented to solve it, which takes considerable development time. For this reason, complicated codes that are often updated with new features, such as weather models, rarely have the capability to compute a full gradient.
There are software packages that solve for the adjoint automatically, or run AD, but these require a programming environment that restricts the user, and may not be efficient in parallel high-performance computing environments. Even when the gradient is calculated automatically via software, the downside of these reverse-mode methods is a potential explosion of memory when creating temporary intermediate variables.
For example, in unsteady fluid flow, the naive adjoint method requires storing the entire time-dependent PDE-solution~\cite{NASA_NielsenDiskin} and hybrid check-pointing schemes~\cite{wang2009minimal} are the subject of active research. Furthermore, in black-box simulations that only return function evaluations given inputs (see \cite[\S1.2]{connDFO} for examples) the derivative may not be available by any means.

 There is a plethora of DFO algorithms, including 
 grid search, Nelder-Mead, (quasi-) Monte-Carlo sampling,  simulated annealing and MCMC methods~\cite{SimulatedAnnealing}. 
 Modern algorithms include randomized methods, Evolution Strategies (ES) such as CMA-ES~\cite{CMA-ES2001}, 
 Hit-and-Run~\cite{Hit-And-Run} and random cutting planes~\cite{RandomCuttingPlane}.
 Textbook DFO methods (\cite[Algo.\ 10.3]{connDFO}, 
\cite[Algo.\ 9.1]{NocedalWright}) are based on interpolation and trust-regions. A limitation of all these methods is that they do not scale well to high-dimensions (beyond 10 or 100).
For this reason, and motivated by the success of stochastic gradient methods in machine learning, we turn to cheaper methods that scale well to high dimensions.

\paragraph{Stochastic DFO}\label{sec:SDFO}
Our basic stochastic subspace descent (SSD) method \eqref{eq: iterations} has been previously explored under the name ``random gradient,'' ``random pursuit,'' ``directional search,'' and ``random search''.
The algorithm dates back to the 1970s, 
  with some analysis (cf.\ 
\cite[Ch.\ 6]{ErmolievWets1988} and
\cite{GavianoRandomSearch1975,SolisWets81}), but it never achieved prominence since  zeroth-order methods are not competitive with first-order methods when the gradient is available. Most analysis has focused on the specific case $\ell=1$ \cite{nesterov2011random, nesterov2017random, stich2013optimization, LewisDirectionalSearch}.
More recently, the random gradient method has seen renewed interest. \cite{LewisDirectionalSearch} analyzes the case when $f$ is quadratic, and \cite{stich2013optimization} provides an analysis (assuming a line search oracle). A variant of the method, using $\ell=1$ and $\bP$ comprised of IID Gaussian entries with  mean zero and variance such that $\Expectation~\bP\bP^\top=\I$, was popularized by Nesterov \cite{nesterov2011random,nesterov2017random},
%
 and various proximal, acceleration and noise-tolerant extensions and analysis appeared in \cite{dvurechensky2017randomized,dvurechensky2018accelerated,GhadimiLan13,ChenWildSTARS2015}.
The novelty of our work, beyond the convergence results and the new variance reduced SSD algorithm, is that we focus specifically on $\ell>1$, and provide practical demonstrations.
%
Another variant of random gradient has recently been proposed in the reinforcement learning community. The Google Brain Robotics team used orthogonal sampling to train artificial intelligence systems~\cite{choromanski2018structured}, but treated it as a heuristic to approximate Gaussian smoothing. 
 Previous use of orthogonal sampling by the OpenAI team~\cite{salimans2017evolution} used variance reduction, but only in the sense of antithetic pairs.
 
We note that the methods we propose may appear similar to stochastic gradient descent (SGD) methods popular in machine learning, but the analysis is completely different. Indeed, even under strong convexity and with a finite-sum structure, SGD methods require diminishing stepsizes and hence converge slowly or else converge only to within a ball around the solution~\cite{NocedalBottou}.

Recent work relevant to the SSD scheme is the SEGA algorithm \cite{SEGA} which takes similar random measurements but then uses them in a  different way, following their ``sketch-and-project'' methodology. Given a sketch $\bP^\T \nabla f(\x_k)$ and an estimate of the gradient at the previous step $\textbf{g}_{k-1} \approx \nabla f(\x_{k-1})$, the simplest version of SEGA finds a new estimate
\begin{equation*}
    \textbf{g}_k = \argmin_{\textbf{g}}\, \|\textbf{g} - \textbf{g}_{k-1}\|_2 \quad \text{subject to } \bP^\T\textbf{g} = \bP^\T\nabla f(\x_k),
\end{equation*} 
followed by a debiasing step.
Another recent related work is \cite{duchi2015optimal} which considers the case $\ell=1$ and focuses on technical issues related to the small bias introduced by estimation of directional derivatives by finite differences, and presents a different type of result from ours.

\paragraph{Alternatives}
As a baseline, one could use $\mathcal{O}(d)$ function evaluations to estimate $\nabla f(\x)$ using finite differences, which is clearly costly when $d$ is large and evaluating $f(\x)$ is expensive. Once $\nabla f(\x)$ is computed, one can run gradient descent, accelerated variants~\cite{Nesterov1983}, non-linear conjugate gradient methods \cite{cg851}, or quasi-Newton methods like BFGS and its limited-memory variant~\cite{NocedalWright}. We compare with these baselines in the numerical results section, particularly with plain gradient descent as it is less sensitive to step-size and errors compared to accelerated variants, and to BFGS since it is widely used.  We do not consider Newton, nor inexact Newton (Newton-Krylov) solvers such as \cite{biegler2003large, biros2005parallel, biros2005parallel2}, because these require estimating the Hessian as well.
Some authors have recently investigated estimating the Hessian stochastically (or adapting BFGS to work naturally in the stochastic sampling case), such as 
\cite{ByrdPaper, Mokhtari, Schraudolph2007ASQ,  StochasticLBFGS_Moritz,berahas2018derivative,TakacNIPS16,GowerBFGS,MaGoldfarbStochasticQuasiNewton}, but none apply directly to our setup as they require an ERM structure.

 \subsection{Structure of this document and contributions}
In Section \ref{sect: SSD} we investigate convergence of the stochastic subspace descent method under various convexity assumptions, and provide rates of convergence. In Section \ref{sect: SVRG} we recall the Stochastic Variance Reduced Gradient (SVRG) method of \cite{SVRG}, and show how control variates reduce the variance in stochastic subspace descent and improve the rate of convergence. Along the way, we provide a simple lemma that relates convergence in expectation at a linear rate with almost sure convergence of our method. It can also be used to show that the $L_1$ convergence of SVRG and other popular first order methods can be trivially extended to almost sure convergence.  In Section \ref{subsect:synthetic-data} we provide empirical results on a simulated function that Nesterov has dubbed ``the worst function in the world'' \cite{nesterov2013introductory}. In Section \ref{subsect: experiments-GP} we find the optimal placement of inducing points for sparse Gaussian processes in the framework of \cite{titsias2009variational}. As a final empirical demonstration, in Section \ref{subsect: experiments-plate}  we test our algorithms in the PDE-constrained optimization setting on a shape optimization problem. 
For the sake of readability, proofs are relegated to the Appendix. More detailed and didactic proofs of the same results can be found in the supplementary material along with additional information that may be useful for intuition.

In this document, uppercase boldfaced letters represent matrices, lowercase boldfaced letters are vectors. The vector norm is assumed to be the Euclidean 2-norm, and the matrix norm is the operator norm. 

\section{Main results}\label{sect: MainResults} 
\subsection{SSD} \label{sect: SSD}
We now provide conditions under which function evaluations $f(\x_k)$ of stochastic subspace descent converge to a function evaluation at the optimum $f(\x_*)$. In the case of a unique optimum we also provide conditions for the iterates $\x_k$ to converge to the optimum $\x_*$. Stochastic subspace descent, so-called because at each iteration the method descends in a random low-dimensional subspace of a vector space containing the gradient, is a gradient-free method as it only requires computation of directional derivatives at each iteration without requiring direct access to the gradient.
In practice we use scaled Haar-distributed random matrices $(\mathbf{P}_k)$ to define randomized directions (in the subspace) along which to descend at each iteration. However, neither Theorem \ref{thm:convergence}, nor the subsequent theorems require Haar-distributed matrices in particular as long as the random matrices satisfy Assumption (A0) below. See the supplementary material for an algorithm to generate Haar-distributed matrices, and a brief discussion of the advantages of using Haar over random coordinate descent type schemes. 
\begin{assumptions} For the remainder of this section we make use of the following assumptions on the matrices $\mathbf{P}_k$ and the function to be optimized. 
\begin{enumerate}
    \item[(A0)\,\,]  $\mathbf{P}_k \in \reals^{d \times \ell}$, $k=1,2, \ldots$, are independent random matrices such that $\Expectation~ \mathbf{P}_k\mathbf{P}_k^\top  = \I_d$  and $\mathbf{P}^\top _k\mathbf{P}_k=(d/\ell)~\I_{\ell}$ with $d > \ell$.
\end{enumerate}
\begin{enumerate}[label=A\theenumi]
	\item[\refstepcounter{enumi}(A\number\value{enumi})\,\,]    $f: \reals^d \to \reals$ is continuously-differentiable with a $\lambda$-Lipschitz first derivative.
	\item[\refstepcounter{enumi}(A\number\value{enumi})\,\,]The function $f$ attains its minimum $f_*$.
	\item[\refstepcounter{enumi}(A\number\value{enumi})\,\,] For some $0<\gamma \leq \lambda$ (where $\lambda$ is the Lipschitz constant in (A1)) and all $\x \in \reals^d$, the function $f$ satisfies the Polyak-Lojasiewicz (PL) inequality: 
			\begin{equation}\label{eq:PL}
			f(\x) - f_* \leq  \norm{\nabla f(\x)}^2/(2\gamma).
			\end{equation} 
    \item[(\theenumi')]   $f$ is $\gamma$-strongly-convex for some $\gamma > 0$ and all $\x \in \reals^d$. Note, $\lambda \geq \gamma$ where $\lambda$ is the Lipschitz constant in (A1).
    \item[(\theenumi'')]   $f$ is convex and attains its minimum $f_*$ on a domain $\mathcal{D}$, and there is an $R>0$ such that for the parameter initialization $\x_0$,	$\max_{\x,~\x_* \in \mathcal{D}} \{\norm{\x-\x_*} : f(\x) \leq f(\x_0) \} \leq R$.
\end{enumerate}
\end{assumptions}
Coercivity of $f$ implies the existence of the constant $R$ in (A3''), hence Assumption (A3'') is equivalent to coercivity of $f$. For the results below, particularly the rate in Theorem \ref{thm: convergence-convex}, we require knowledge of the value of $R$. Also note that (A3') implies (A3).

Algorithm \ref{Alg: SSD} provides basic pseudocode to implement stochastic subspace descent using the recursion \eqref{eq: iterations}.
\begin{algorithm}[ht]
	\caption{Stochastic subspace descent (SSD)}\label{Alg: SSD}
	\begin{algorithmic}[1]
		\State{\textbf{Inputs:} $\alpha, \ell$ \Comment{step size, subspace rank}}
		\State{\textbf{Initialize:} $\x_0$ \Comment{arbitrary initialization}}
		\For {k = 1, 2, \ldots}
		\State Generate $\textbf{P}_k$
		\State  Apply recursion \eqref{eq: iterations}
		\EndFor
		\end{algorithmic}
\end{algorithm}
\begin{theorem}[Convergence of SSD]\label{thm:convergence}
 Assume (A0), (A1), (A2), (A3) and let $\x_0$ be an arbitrary initialization. Then recursion \eqref{eq: iterations} with $0<\alpha <2\ell/d\lambda$ results in  $f(\x_k) \convas f_*$ and $f(\x_k) \overset{L^1}{\to} f_*$. 
\end{theorem}

Formal definitions of the different modes of convergence are provided in the supplementary material. Since non-convex functions may satisfy the PL-inequality \eqref{eq:PL}, Theorem \ref{thm:convergence} provides a convergence result for well-behaved non-convex functions, or (see Corollary \ref{corr:strong-convexity} below) linear convergence results for non-strongly-convex functions such as linear least squares with a matrix that is not full column rank \cite[\S2.3]{Karimi2016LinearCO}.

\begin{corollary}[Convergence under strong-convexity and rate of convergence]\label{corr:strong-convexity}
\vspace*{-.9em}
\begin{enumerate}[label=(\roman*)]
    \item Assume (A0), (A1), (A2), (A3') and let $\x_0$ be an arbitrary initialization. Then recursion \eqref{eq: iterations} with $0<\alpha <2\ell/d\lambda$ results in $\x_k \convas \x_*$ where $\x_*$ is the unique minimizer of $f$.
    
    \item Assume (A0), (A1), (A2), and (A3). Then with  $\alpha=l/d\lambda$, the recursion \eqref{eq: iterations} attains the following expected rate of convergence
	\begin{equation}
	0 \leq \Expectation f(\x_k)-f_*  \leq \beta^k (f(\x_0)-f_*) ,\quad \beta = 1 - \ell \gamma/d \lambda.
	\end{equation}

\end{enumerate}
\end{corollary}

We note that if $\ell = d$ we recover a textbook rate of convergence for gradient descent \cite[\S 9.3]{BoydVandenbergheBook}. Under the more restrictive assumption of strong-convexity the result of Corollary \ref{corr:strong-convexity} is much stronger than Theorem \ref{thm:convergence}; we get almost sure convergence of the function evaluations and of the iterates to the optimal solution at a linear rate. In inverse problems it is typically the convergence of $\x_k$ not of $f(\x_k)$ that is sought. Furthermore, if either assumption (A3) or (A3') is satisfied, SSD has a linear rate of convergence. In addition to the convergence rate, we provide guidance for selecting an optimal step-size, which in theory alleviates the need for the user to tune any hyperparameters. The step-size is easily interpreted as the quotient of the sample rate $\ell/d$ and the global Lipschitz constant. This is consistent with our intuition: more directional derivatives make for a better approximation of the gradient. In practice the Lipschitz and strong-convexity constants are rarely known and a line search can be used.

Thus far, we have provided results for cases when the function satisfies the PL-inequality or is strongly-convex. The proof in the convex case is different, but substantively similar to a proof of coordinate descent on convex functions found in \cite{wright2015coordinate}.

\begin{theorem}[Convergence under convexity]\label{thm: convergence-convex}
Assume (A0), (A1), (A2), (A3''). Then recursion \eqref{eq: iterations}  with $0<\alpha <2\ell/d\lambda$ gives $\Expectation f(\x_k) - f_*  \leq 2d \lambda R^2/k \ell$.
\end{theorem}

\subsection{Variance Reduced Stochastic Subspace Descent (VRSSD)}\label{sect: SVRG}
In its most basic formulation SGD suffers from a sub-linear rate of convergence. This is due to a requirement for the step-size to gradually decay to assure convergence of the iterates. Several methods have been designed to address this short-coming, one of which is SVRG \cite{SVRG}, which controls the variance of the gradient estimate using a control variate, and admits convergence with a constant step-size, resulting in a linear rate of convergence. In this section we recall some results regarding SVRG from \cite{SVRG} before proving almost sure convergence of the iterates and function evaluations, and ultimately adapting it to the framework of SSD. SVRG minimizes the functional $f(\x)=(1/n)\sum_{i=1}^n f_i(\x)$ via the iterative scheme
\begin{equation}\label{eq: SVRG1}
\x_{s,k} = \x_{s,k-1} - \alpha\left(\nabla f_{i_k}(\x_{s,k-1})- \nabla f_{i_k}(\tilde{\x}_s) + \tilde{\boldsymbol{\mu}}_s\right),
\end{equation}
where  $i_k \in \{1, \ldots, n \}$ is chosen uniformly at random, $\tilde{\x}_s$ is a past iterate updated every $m$ iterations, and the full gradient is also calculated every $m$ iterations: 
\begin{equation}\label{eq: SVRG2}
\tilde{\boldsymbol{\mu}}_s=\nabla f(\tilde{\x}_s) = (1/n)\medmath{\sum_{i=1}^n} \nabla f_{i}(\tilde{\x}_s).
\end{equation}
 In \cite{SVRG} a proof is presented for the case where $\tilde{\x}_s$ is updated every $m$ iterations by setting $\tilde{\x}_s = \x_{s,J_s}$, where $J_s$ is chosen uniformly at random from $\{1, \ldots, m\}$. The authors call this updating scheme ``Option 2'' and recommend an alternative ``Option 1'', without proof, that involves setting $J_s=m$ instead. A proof of convergence for Option 1 is provided in \cite{Tan2016BarzilaiBorweinSS}. The SVRG method is presented in Algorithm \ref{alg:SVRG}, and the convergence result of \cite{SVRG}, is presented in Theorem \ref{SVRG}.
 
 \begin{algorithm}
	\caption{SVRG(\cite{SVRG})}\label{alg:SVRG}
	\begin{algorithmic}
		\State{\textbf{Inputs:} $m, \alpha$ \Comment{memory parameter, step-size}}
		\State{\textbf{Initialize:} $\tilde{\x}_0$}
		
		\For {$s = 1, 2, \ldots$}
		\State $\x_{s,0} = \tilde{\x}_{s-1}$
		\State $\tilde{\boldsymbol{\mu}}_s= \frac{1}{n}\sum_{i=1}^n \nabla f_{i}(\tilde{\x}_s).$
		\For {k = 1, \ldots, m}
		\State Choose $i_k \in \{1, \ldots, n \}$ uniformly at random
		\State $\x_{s,k} = \x_{s,k-1} - \alpha\left(\nabla f_{i_k}(\x_{s,k-1})- \nabla f_{i_k}(\tilde{\x}_s) + \tilde{\boldsymbol{\mu}}_s\right) $
		\EndFor
		\State Option 1: $J_s = m$
		\State Option 2: $J_s = \text{unif}\{1, \ldots, m\}$
		\State {$\tilde{\x}_s = \x_{s,J_s}$}
		\EndFor
	\end{algorithmic}
\end{algorithm}
\begin{theorem}[Convergence of SVRG \cite{SVRG}]\label{SVRG}
	Assume that each $f_i$ is convex with $\lambda$-Lipschitz derivative and that $f(\x)$ is $\gamma$-strongly-convex. Assume that the memory parameter $m$ is sufficiently large so that
	\begin{equation*}
	\beta = \frac{1}{\gamma \alpha (1-2\lambda\alpha)m} + \frac{2\lambda\alpha}{1-2\lambda\alpha} <1.
	\end{equation*}
	Then, for the sequences \eqref{eq: SVRG1} and \eqref{eq: SVRG2}, using Option 2 to update $\tilde{\x}_s$, we obtain the rate of convergence: $0 \leq \Expectation f(\tilde{\x}_s) - f(\x_*) \leq \beta^s\left[f(\tilde{\x}_0) - f(\x_*)\right]$.
	
\end{theorem}

The proof, omitted here, can be found in \cite{SVRG}. Since $f(\tilde{\x}_0) - f(\x_*)$ is a constant and $\lim_{s\to\infty}\beta^{s}=0$, Theorem \ref{SVRG} implies that $f(\tilde{\x}_s) \overset{L^1}{\longrightarrow} f(\x_*)$ . We can also show almost sure convergence. This is a particular application of the following simple lemma:

\begin{lemma}[a.s.\ convergence from a linear rate of convergence]\label{thm:SVRG-AS}
	Consider minimizing a function $f:\reals^d \to \reals$, and set $f_* \defeq \min_{\x\in \reals^d} f(\x)$. If a minimization algorithm converges in expectation at a linear rate; that is, if there exists $\beta \in (0,1)$ and a constant $C>0$ such that
	\begin{equation*}
	0 \leq \mathbb{E}f(\x_k)-f_* \leq C\beta^k,
	\end{equation*} 
then $f(\x_k) \convas f_*$.		
If in addition, $f$ is strongly-convex, then $ \x_k \convas \x_*$.

\end{lemma}
The distinction between modes of convergence can be important, see for instance \cite[p.4]{ferguson2017course} for a simple example of a sequence that converges in $L^r$ for any $r >0$ but does not converge almost surely. As discussed in \cite{lin2015universal}, the class of algorithms that satisfy the conditions required to invoke Lemma \ref{thm:SVRG-AS} is substantial. This class includes many of the most popular gradient-based optimization methods such as block-coordinate descent, SAG \cite{SAG}, SAGA \cite{SAGA}, SDCA \cite{SDCA}, SVRG \cite{SVRG}, Finito/MISO \cite{FINITO, MISO}, and their proximal variants. Though the conditions of Lemma \ref{thm:SVRG-AS} are sufficient for almost sure convergence, they are not necessary: there are optimization algorithms, such as stochastic gradient descent, that converge almost surely for which the rate of convergence is sub-linear. Of the aforementioned algorithms only SVRG is relevant to our present discussion:

\begin{remark}[a.s.\ convergence of SVRG]\label{remark:svrg-as}
	For the SVRG algorithm discussed in Theorem \ref{SVRG},  $f(\tilde{\x}_s) \convas f(\x_*)$ and $\tilde{\x}_s \convas \x_*$.
\end{remark}

\begin{proof}
	The problems for which the theory of SVRG can be applied are strongly-convex. Also, by construction $\beta \in (0,1)$, and $f(\tilde{\x}_0)-f(\x_*) \geq 0$. Thus the result follows from Lemma \ref{thm:SVRG-AS}.
\end{proof}

As stated, SVRG is stochastic in that it randomly samples observations from the dataset. Our goal is different in that we are interested in random sampling along the dimensions of input space, motivating a new algorithm we call Variance Reduced Stochastic Subspace Descent (VRSSD). The main algorithm is described in Algorithm \ref{alg:varprod}. Notice that the similarity of SVRG and VRSSD is in their use of a control variate; in VRSSD we also incorporate a parameter $\eta$ to weight the influence of the control variate.

\begin{algorithm}
	\caption{Variance reduced stochastic subspace descent}\label{alg:varprod}
	\begin{algorithmic}
		\State{\textbf{Inputs}: $m, \alpha$ \Comment{memory parameter, step-size}}
		\State{\textbf{Initialize}: $\tilde{\x}_0$}
		
		\For {$s = 1, 2, \ldots$}
		\State $\x_{s,0} = \tilde{\x}_{s-1}$
		\For {k = 1, \ldots, m}
		\State \textbf{Generate} $\P$
		\State $\x_{s,k} = \x_{s,k-1} - \alpha [\P\P^\top \nabla f(\x_{s,k-1}) -\eta_{s,k-1} (\P\P^\top -\I)\nabla f(\tilde{\x}_{s-1}) ]$
		\EndFor
		\State  Option 1: $J_s = m$
		\State  Option 2: $J_s = \text{unif}\{1, \ldots, m\}$
		\State {$\tilde{\x}_s = \x_{s,J_s}$}
		\EndFor
	\end{algorithmic}
\end{algorithm}

Theorem \ref{thm:VarianceReducedRandomGradient} allows for, but (unlike SVRG) does not require, $f$ to have a finite-sum representation and as a consequence it does not require Lipschitz-continuity of the individual samples $f_i$ as in Theorem \ref{SVRG}. 

\begin{theorem}[a.s.\ convergence of VRSSD]\label{thm:VarianceReducedRandomGradient}
Let $f: \reals^d \to \reals$ be a  $\gamma$-strongly convex continuously-differentiable function  with $\lambda$-Lipschitz gradient for some $\lambda \geq \gamma > 0$. Define $\rho = d/\ell$ and let $\x_*$ be the minimizer of $f$. 

\begin{enumerate}[label=(\roman*)]
    \item Fix a memory parameter $m$, a sampling rate $\rho>2$, and a step-size $\alpha$ such that
	\begin{equation*}
	\beta =  \frac{1}{\alpha \gamma m (1 - \alpha \lambda \rho)} + \frac{\alpha \lambda (\rho-1)}{1- \alpha \lambda \rho}<1.	\end{equation*}
	Then, Algorithm \ref{alg:varprod} with Option 2 and $\eta\equiv 1$ results in:
	\begin{equation*}
	0 \leq \Expectation f(\tilde{\x}_s) - f(\x_*) \leq \beta^s \left[f(\tilde{\x}_0) - f(\x_*)\right],
	\end{equation*}
	and $\tilde{\x}_s \convas \x_*$.
	\item
    Fix $m$, $\rho>2$, and $\alpha$ such that
	\begin{equation*}
	\beta =\frac{1}{\alpha \gamma m (1 - \alpha \lambda \rho)} < 1.
	\end{equation*}
Then, with the control variate parameter: $\eta_{s,k-1} =  \nabla f(\tilde{\x}_{s-1})^\top  \nabla f(\x_{s,k-1})/\norm{\nabla f(\tilde{\x}_{s-1})}^2$,
Algorithm \ref{alg:varprod} using Option 2 leads to:
	\begin{equation*}
	0\leq \Expectation f(\tilde{\x}_s) - f(\x_*)) \leq \beta^s \left[f(\tilde{\x}_0) - f(\x_*))\right],
	\end{equation*}
	and $\tilde{\x}_s \convas \x_*$.
	\end{enumerate}
\end{theorem}

For appropriately chosen parameters $m$ and $\rho$, the variance reduction introduced in Theorem \ref{thm:VarianceReducedRandomGradient} leads to an improved convergence rate compared to the simpler method presented in Theorem \ref{thm:convergence}. The benefits can be seen in practice as shown in Section \ref{subsect: experiments-plate}. Part ($ii$) of Theorem \ref{thm:VarianceReducedRandomGradient} shows that it is possible to further improve the rate by selecting a control variate parameter $\eta$ to optimally weight the impact of the full gradient held in memory. In theory we see the following effects on the rate for VRSSD with the optimal control variate, all of which are consistent with our intuition: (i) a large sampling rate $\rho$ improves per-iteration convergence; (ii) a large $m$ improves per-iteration convergence, but each iteration is only counted after $m$ inner loops have been completed, so the effect of $m$ cancels (in theory); (iii) poor conditioning ($\gamma \approx 0$) has a deleterious effect on the convergence rate commensurate with its effect on the convergence rate of gradient descent; (iv) a large Lipschitz constant must be offset by a small step-size; and (v), it is possible to find the optimal step for a known local Lipschitz constant and strong-convexity constant, and a fixed $m$.

In many of the problems for which SSD is useful it is undesirable to compute the term $\nabla f(\x_{s,k-1})$ to find the optimal control variate as this would amount to $d$ additional function evaluations per iteration. In practice we therefore recommend approximating it with the current estimate of the gradient instead, noting that
\begin{equation}\label{eqn: VRSSD-approximation}
\nabla f(\x_{s,k-1})\approx \mathbf{P}_{s,k}\mathbf{P}_{s,k}^\top \nabla f(\x_{s,k-1}).
\end{equation}
Though this approximation is not supported by the theory, empirical performance in Section \ref{sect: EmpiricalResults} suggests that such an approximation is reasonable. We chose the approximation in \eqref{eqn: VRSSD-approximation} because it is an unbiased estimate of the optimal control variate parameter at each iteration and has the advantage that no additional function evaluations are required for each iteration. Of course, some of the intuitions outlined above may not hold; this in itself is not necessarily a problem as the Lipschitz and strong convexity constants are rarely known and step-size is determined by line search; in our experiments we use an Armijo backtracking line search as discussed in \cite{NocedalWright}. 

\section{Experimental results}\label{sect: EmpiricalResults} In this section we provide results for a synthetic problem, a problem from the machine learning literature, and a shape-optimization problem. For the sake of fair comparisons, we compare to zeroth-order variants of gradient descent and BFGS where the gradient is calculated via forward-differences. In the synthetic problem we also compare to a randomized block-coordinate descent. For VRSSD we use the approximation  \eqref{eqn: VRSSD-approximation} of the optimal $\eta$, and we use Option 1 for updating $\tilde{\x}_s$.

\subsection{Synthetic data}\label{subsect:synthetic-data}
\begin{figure}[ht]
      \centering
      \includegraphics[width=.42\linewidth]{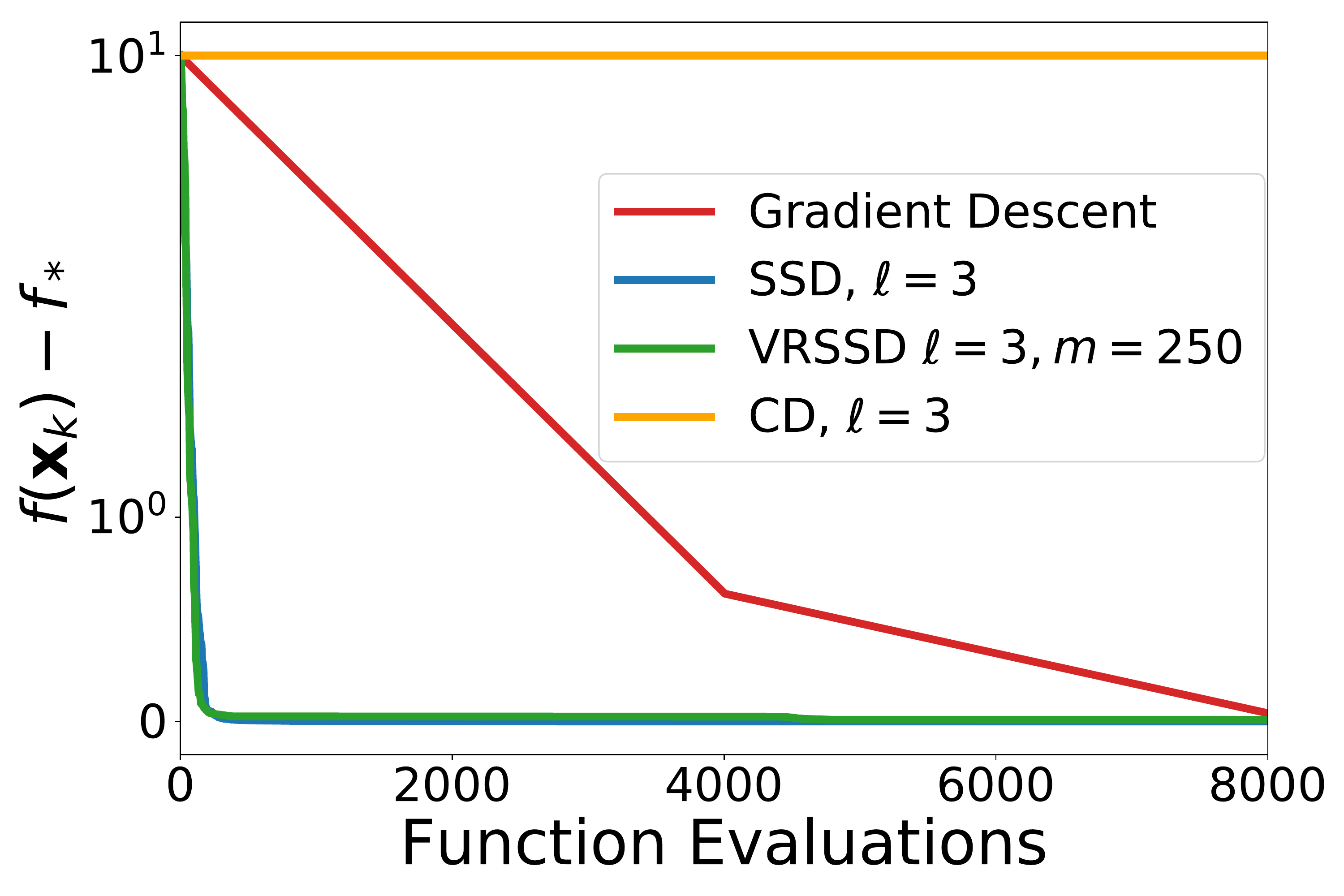}
      \includegraphics[width=.42\linewidth]{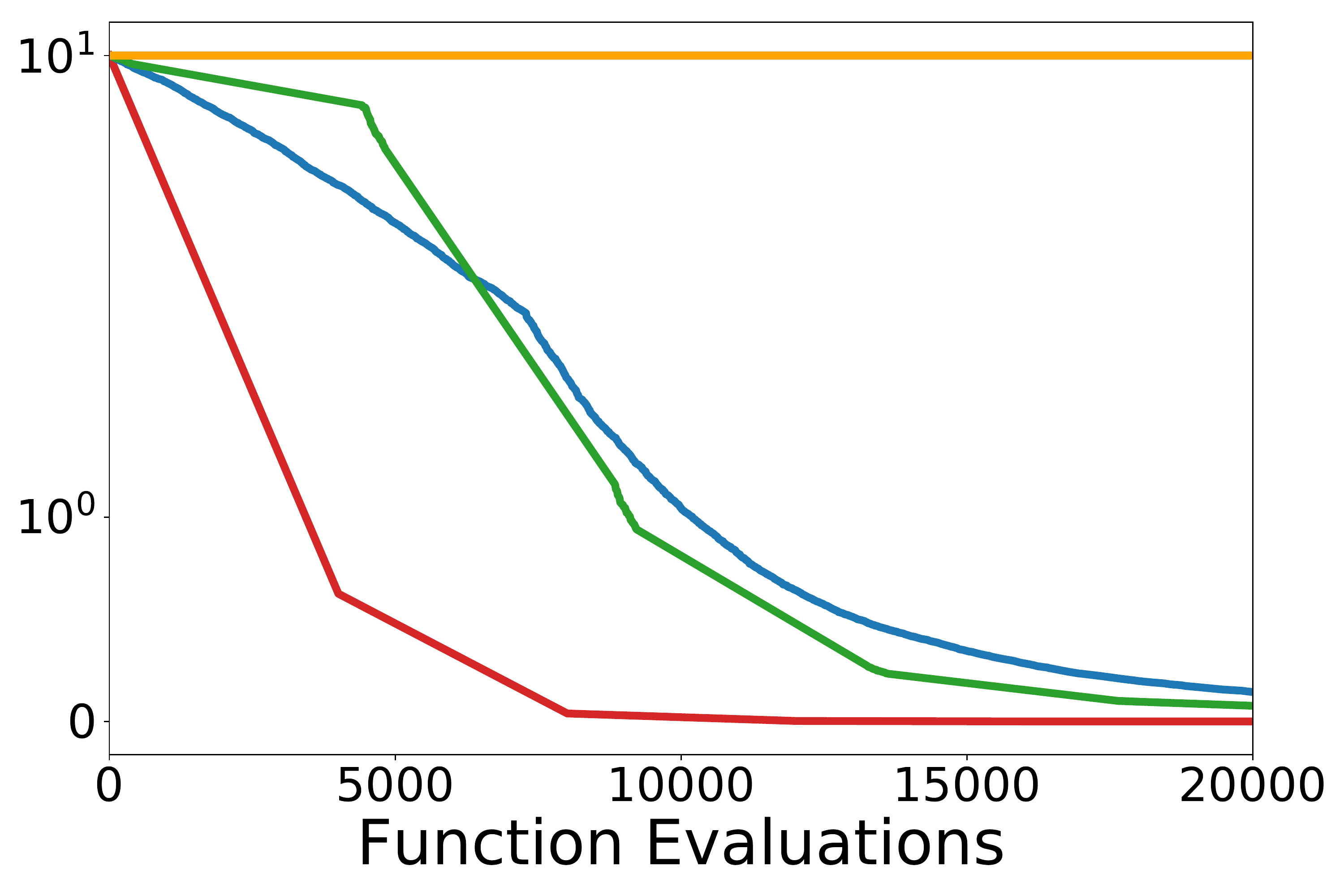}
      \vspace*{-1em}
      \caption{Left: A function from the class \eqref{eq: NesterovWorst} with low intrinsic dimension and high global Lipschitz constant: $r=10,~d=4001,~\lambda=80$. SSD is not visible as it lies directly under the VRSSD curve. Right: A function from the same class with high intrinsic dimension and low global Lipschitz constant: $r=2000,~d=4001,~\lambda=0.8$. $CD$ in this case is randomized block-coordinate descent. }
      \label{fig: NesterovExample}
  \end{figure}

We explore the performance of the algorithms on Nesterov's ``worst function in the world'' \cite{nesterov2013introductory}.  Fix a Lipschitz constant $\lambda>0$ and let
\begin{equation}\label{eq: NesterovWorst}
    f_{\lambda,r}(\x) = \lambda (( x_1^2 + \medmath{\sum_{i=1}^{r-1}} \left( x_i - x_{i+1}\right)^2 + x_r^2 )/2 - x_1 )/4,
\end{equation}
where $x_i$ represents the $i^{\text{th}}$ coordinate of $\x$ and $r < d$ is a constant integer that defines the intrinsic dimension of the problem. This function is convex with global minimum $f_*= -\lambda r/8(r+1)$.
Nesterov has shown that a broad class of iterative first-order solvers perform poorly when minimizing $f_{\lambda,r}$ with $x_0 = 0$. We compare our results to randomized block-coordinate descent and gradient descent for two illustrative examples from the family of functions defined by \eqref{eq: NesterovWorst}. We do not compare to BFGS as we do not run enough iterations for it to differ substantially from gradient descent. For the randomized methods we use $\ell=3$; in Section \ref{subsect: experiments-GP} we explore the impact that different choices of $\ell$ have on convergence. We expect coordinate descent to perform poorly as there is only a single coordinate of this function that provides a descent direction from the initialization; for high-dimensional problems it becomes exceedingly unlikely that the ``correct'' coordinate is chosen. Also, since it is a \emph{block} coordinate descent with $\ell=3$, the impact of the other two coordinates chosen outweighs the benefit, and the step-size gets pushed to zero; hence, in both examples shown in Figure \ref{fig: NesterovExample} we see no improvement from coordinate descent. Had we used instead $\ell=1$, with coordinate descent we could expect improvement in the objective once every $d$ iterations, rather than at every iteration as in the other methods.

If the function is intrinsically low-dimensional but is embedded in a high-dimensional space, then SSD and VRSSD perform exceptionally well (see left plot of Figure \ref{fig: NesterovExample}). This holds true for the Gaussian process example of Section \ref{subsect: experiments-GP} as well. On the other hand, if the space is truly high-dimensional then the stochastic methods do not appear to provide any benefit. 

\begin{figure}[ht]
      \centering
      \includegraphics[width=.40\linewidth]{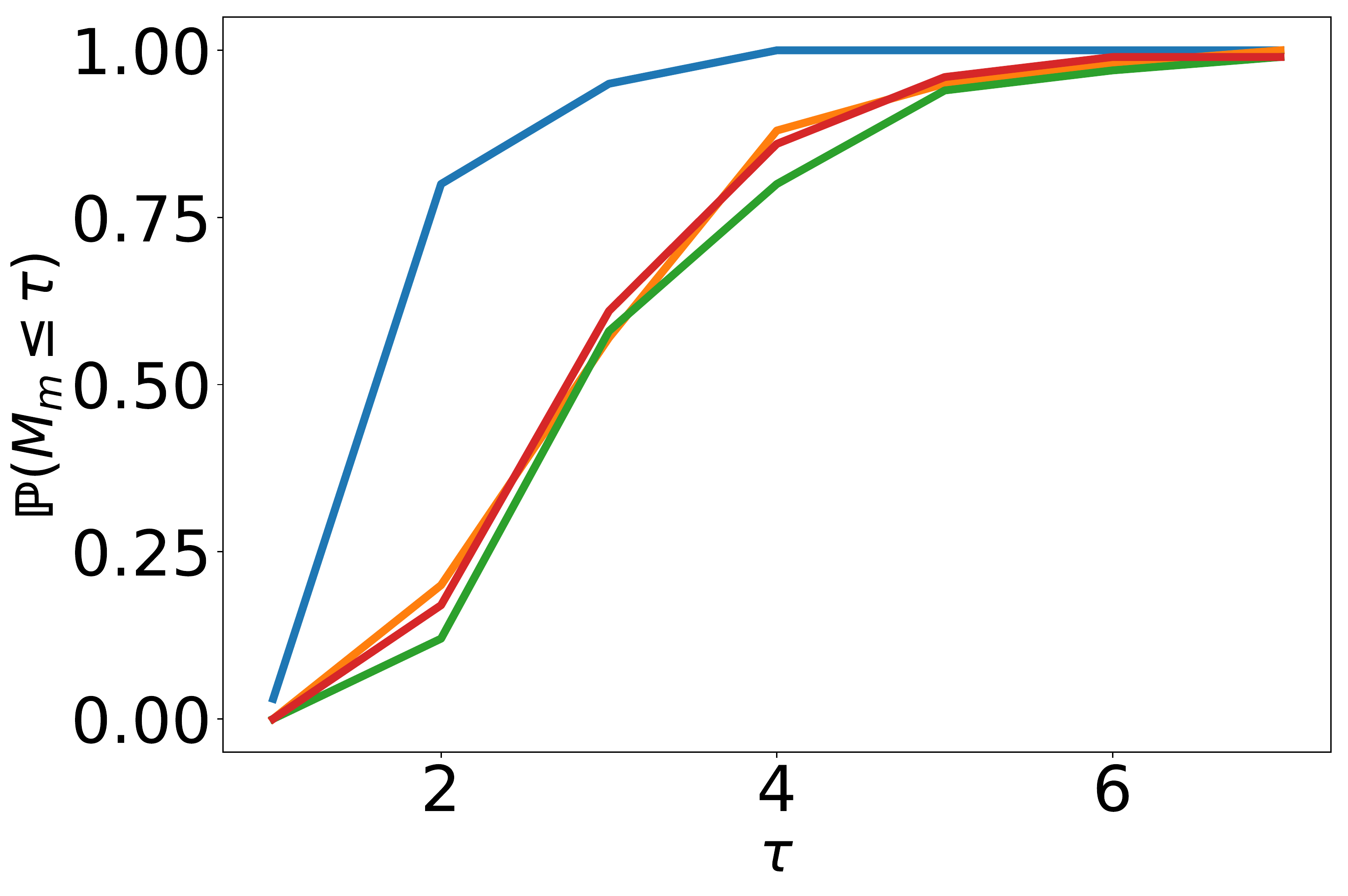}
        \includegraphics[width=.40\linewidth]{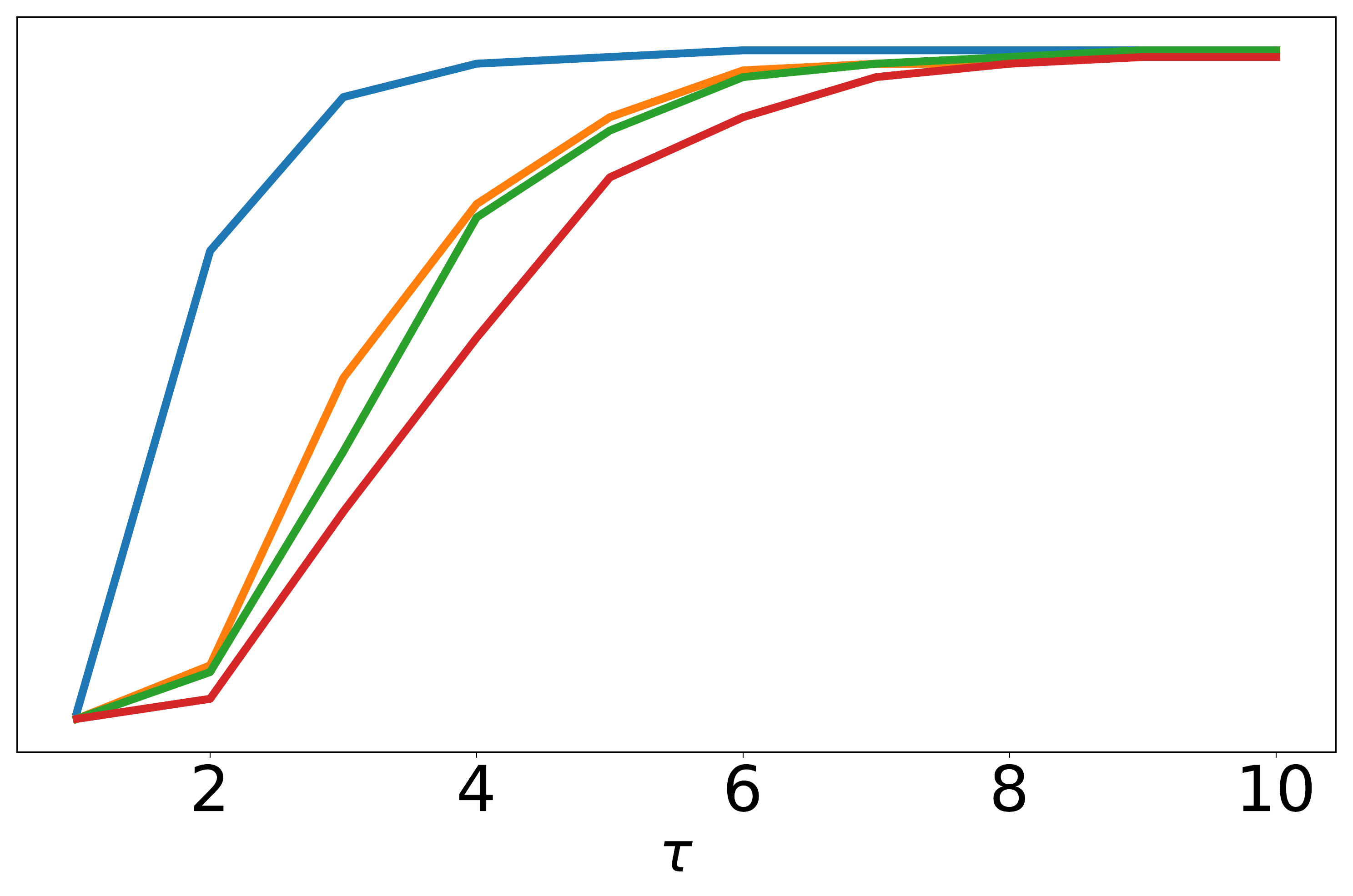}
        \includegraphics[width=.40\linewidth]{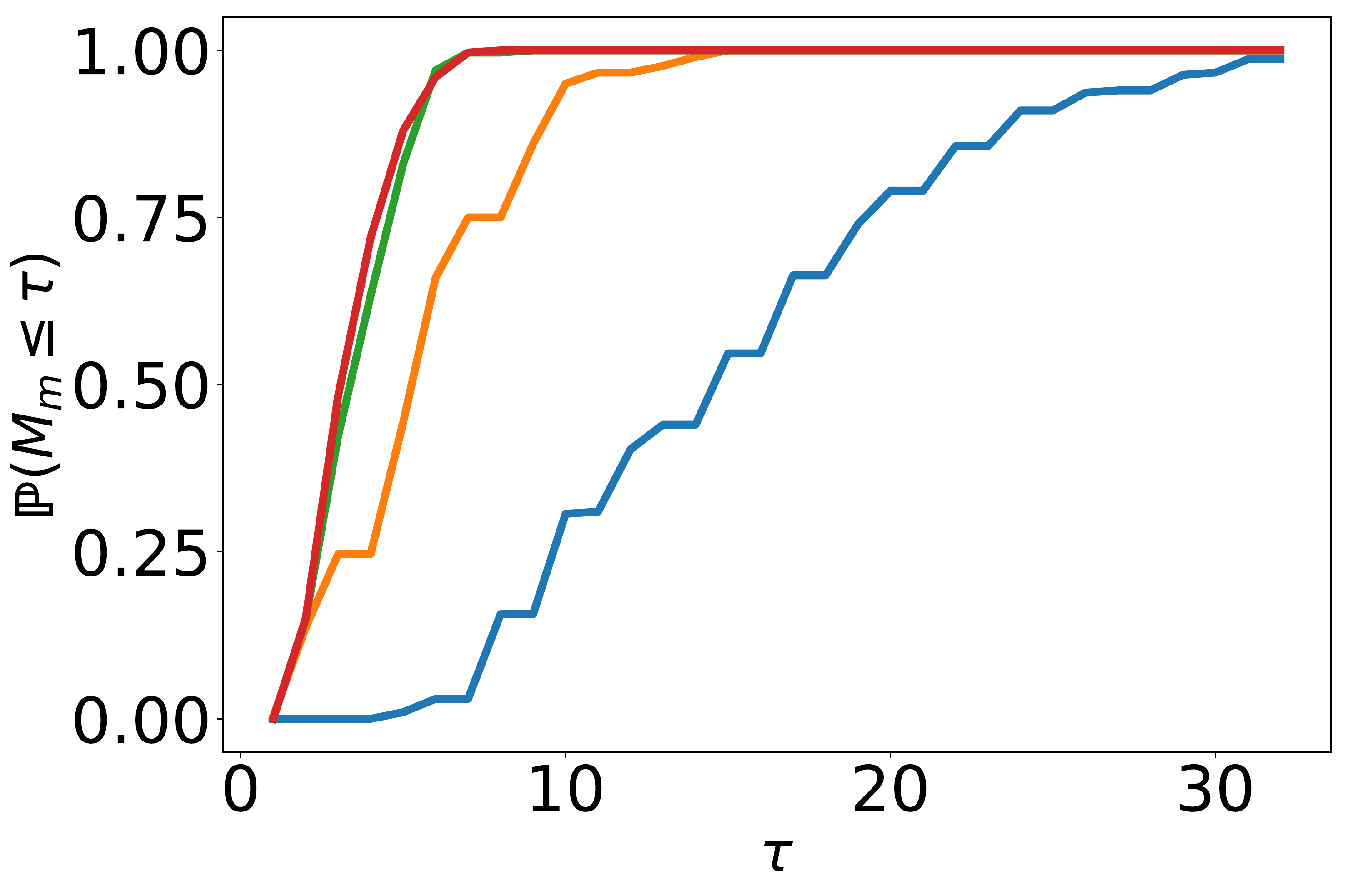}
        \includegraphics[width=.40\linewidth]{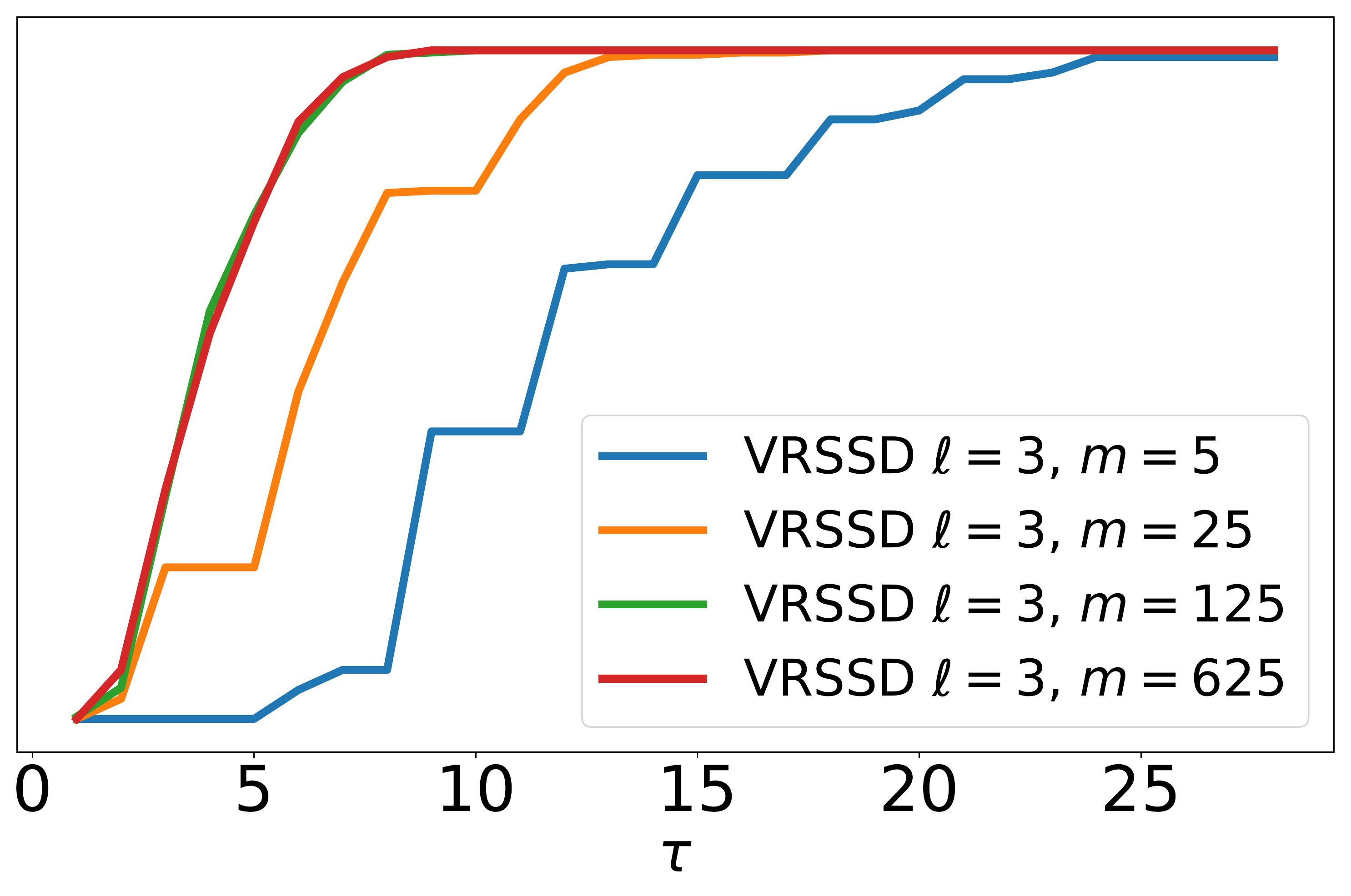}
    \vspace*{-1em}
      \caption{Top-left: L=800, r=50. Top-right: L=0.8, r=50. Bottom-left: L=800, r=5. Bottom-right: L=0.8, r=5. $M_{m}$ is the number of function evaluations required to attain a pre-specified objective value for various values of $m$. $\tau$ is the ratio of the number of function evaluations required for the fastest realization compared with the number required for any realization. That is, $\tau=10$ corresponds to a solution that requires 10 times as many function evaluation as the fastest solution.}
      \label{fig: VRSSD-performanceprofile}
  \end{figure}

There are hyperparameters to set in our methods. We reserve discussion of $\ell$ for Section \ref{subsect: experiments-GP}. In Figure \ref{fig: VRSSD-performanceprofile} we use performance profiles \cite{dolan2002benchmarking} to examine the impact of varying $m$ in VRSSD. A performance profile is conducted by running each parameterization with 300 randomized restarts with termination after some pre-specified tolerance for accuracy has been reached. We count the proportion of realizations from each parameterization that achieves the specified tolerance within $\tau$ function evaluations, where $\tau=1$ is the fewest function evaluations required in any of the trials, $\tau=2$ is twice as many function evaluations, etc. We use a function from the family \eqref{eq: NesterovWorst} and fix $d=101$. The threshold we use for the objective value is within $10^{-3}$ of the optimal objective.

Figure \ref{fig: VRSSD-performanceprofile} shows that the Lipschitz constant has little bearing on the relative performance of the different parameterizations, but the intrinsic dimension has a large role to play. The top plots in Figure \ref{fig: VRSSD-performanceprofile} paint a bleak picture for these stochastic methods: the more similar to gradient descent ($m \approx 1$) the better their performance. Conversely, the bottom plots show that randomized methods provide a distinct advantage when the data is intrinsically low-dimensional. Recent work \cite{fefferman2016testing} takes a rigorous look at the manifold hypothesis which suggests that most real, high-dimensional data tend to lie near a low-dimensional manifold. This is corroborated in \cite{udell2019big} where the authors provide theoretical rationale for the observed low-rank structure of large matrices. These two theoretical papers, along with ample empirical evidence \cite{NIPS2002_2203, Bell2007LessonsFT, levina2005maximum, tropp2017practical}, suggest that the lower plots in Figure \ref{fig: VRSSD-performanceprofile}, and the left plot in Figure \ref{fig: NesterovExample} are more representative of problems that arise in practice. We have presented two extremes in an effort to provide a heuristic for choosing $m$ if the intrinsic dimension can be approximated. \cite{levina2005maximum} provides a method for estimating the intrinsic dimension via maximum likelihood. More recently, \cite{ceruti2014danco} provides a method for estimating the intrinsic dimension of an arbitrary data set by exploiting concentration inequalities of norms and angles for random samples in high dimensional space. In \cite{ceruti2014danco} the authors also conduct a thorough review of the literature and compare to several other methods for estimating intrinsic dimension.

\subsection{Sparse Gaussian process parameter estimation}\label{subsect: experiments-GP}
We apply the methods developed in Section \ref{sect: MainResults} to hyperparameter estimation for Gaussian processes used in regression,
where the goal is inference on a function $f:\reals^d\to\reals$  based on noisy observations at $n$ points
$\z_1,\ldots,\z_n$.  The $n$ observations are modeled as $y_i = f(\z_i) + \epsilon_i$, where the function $f$ is a zero-mean Gaussian process with covariance function 
$\mathbb{C}\mathrm{ov}(f(\z_i),f(\z_j)) = K(\z_i, \z_j; \boldsymbol{\theta})$,
and $K(\cdot\,,\cdot\,;\boldsymbol{\theta})$ is a symmetric positive-definite kernel that depends on a vector of parameters $\boldsymbol{\theta}$.
The process $f$ is assumed to be independent of the noise vector $(\epsilon_1,\ldots,\epsilon_n)^\top \sim N(\zero,\sigma^2\Ib)$ with unknown variance $\sigma^2$.
The covariance matrix of the vector $\fb=(f(\z_1),\ldots,f(\z_n))^\top $ will be written as $\boldsymbol{\Sigma}_{\fb} =  \mathbb{V}\mathrm{ar}(\fb)$. So, $(\boldsymbol{\Sigma}_{\fb})_{ij} = K(\z_i, \z_j;\boldsymbol{\theta})$.
Maximum likelihood estimates of the parameters $\boldsymbol{\Theta} = [\boldsymbol{\theta}, \sigma^2]$ are obtained by maximizing the 
log-marginal likelihood of observations 
$\y=(y_1,\ldots,y_n)^\top $ with density $p_{\y}$ \cite{Rasmussen2009GaussianPF}: $\ell(\boldsymbol{\Theta};\y) = \log p_{\y}(\y;\boldsymbol{\Theta})$. %
Difficulties arise in finding the optimal parameterization of Gaussian processes when the number of observations is large: due to the inversion and determinant calculation in  $\ell(\boldsymbol{\Theta};\y)$, the cost of this maximization is $\mathcal{O}(n^3)$. Various approximations exist, perhaps most popular among them in the machine learning community is that of \cite{titsias2009variational}. The basic idea is as follows: Let $\widetilde{\z}_1,\ldots,\widetilde{\z}_m$, $m<n$, be points
in $\reals^d$ different from the original $\z_1,\ldots,\z_n$ (these are called ``inducing points''), and let
$\widetilde{\fb} = (f(\widetilde{\z}_1),\ldots f(\widetilde{\z}_m))^\top $. Let $g$ be a multivariate density function
on $\reals^d$, and define the density $h(\fb,\widetilde{\fb}) = p_{\fb\mid\widetilde{\fb}}(\fb\mid\widetilde{\fb}) g(\widetilde{\fb})$
on $\reals^{n+m}$, where $p_{\fb\mid\widetilde{\fb}}$ is the conditional density of $\fb$ given $\widetilde{\fb}$.
The Gaussian density that minimizes the Kullback-Leibler divergence between the density $h$ and the posterior density
of $(\fb,\widetilde{\fb})$ given $\y$ leads to the following lower bound for the loglikelihood \cite{titsias2009variational}:
\begin{equation}\label{eq:llbound}
   \ell(\boldsymbol{\Theta};\y) \geq F(\widetilde{\z}_1,\ldots\widetilde{\z}_m,\boldsymbol{\Theta})=\widetilde{\ell}(\boldsymbol{\Theta};\y) - \mathrm{tr}
   (\mathbb{V}\mathrm{ar}(\fb\mid \widetilde{\fb}))/2\sigma^2.
\end{equation}
Here
$\widetilde{\ell}$ is the loglikelihood of the multivariate Gaussian $N(\zero, \widehat{\boldsymbol{\Sigma}}_{\fb})$, where $
\widehat{\Sigmab}_{\fb} = \Sigmab_{\fb} - \mathbb{V}\mathrm{ar}(\fb\mid \widetilde{\fb})=
\mathbb{C}\mathrm{ov}(\fb,\widetilde{\fb}) \,\Sigmab_{\widetilde{\fb}}^{-1}\,\mathbb{C}\mathrm{ov}(\widetilde{\fb},\fb)
$ 
is the  the Nystr\"om approximation of $\Sigmab_{\fb}$ introduced in \cite{williams2001using}. 
Note that
$
\mathrm{tr}(\mathbb{V}\mathrm{ar}(\fb\mid \widetilde{\fb}))=
\Expectation\,\|\,\Expectation\,(\fb\mid \widetilde{\fb}) - \fb\,\|^2.
$
That is, the Nystr\"om approximation is controlled by the mean-squared error of the prediction of $\fb$ based
on $\widetilde{\fb}$. The trace term in \eqref{eq:llbound} controls the size of this mean-square error. We note that the Nystr\"om approximation is a low-rank non-negative update of $\boldsymbol{\Sigma}_{\fb}$ and that optimal approximations of this form were studied in \cite{spantini2015optimal, spantini2017goal} using Rao's geodesic distance between SPD matrices as a measure of optimality. We do not pursue this approach here but recognize it as a potential alternative to the method of \cite{titsias2009variational} that we use. Derivations of the conditional means and variances of a multivariate Gaussian can be found in, for example,
\cite{LuisBook}. 

Gradient-based methods are used to find an optimal placement of the $m$ inducing points by maximizing the lower
bound in \eqref{eq:llbound}, which we re-state as a function of $\x = [\widetilde{\z}_1, \ldots, \widetilde{\z}_m; \boldsymbol{\Theta}]$ to be consistent with notation in previous sections:
\begin{equation}\label{eqn: VariationalUpperBound}
F(\x) = \widetilde{\ell}(\boldsymbol{\Theta};\y) -\,\mathrm{tr}
   (\mathbb{V}\mathrm{ar}(\fb\mid \widetilde{\fb}))/2\sigma^2.
\end{equation}

\begin{figure}[ht]
      \centering
      \includegraphics[width=.4\linewidth]{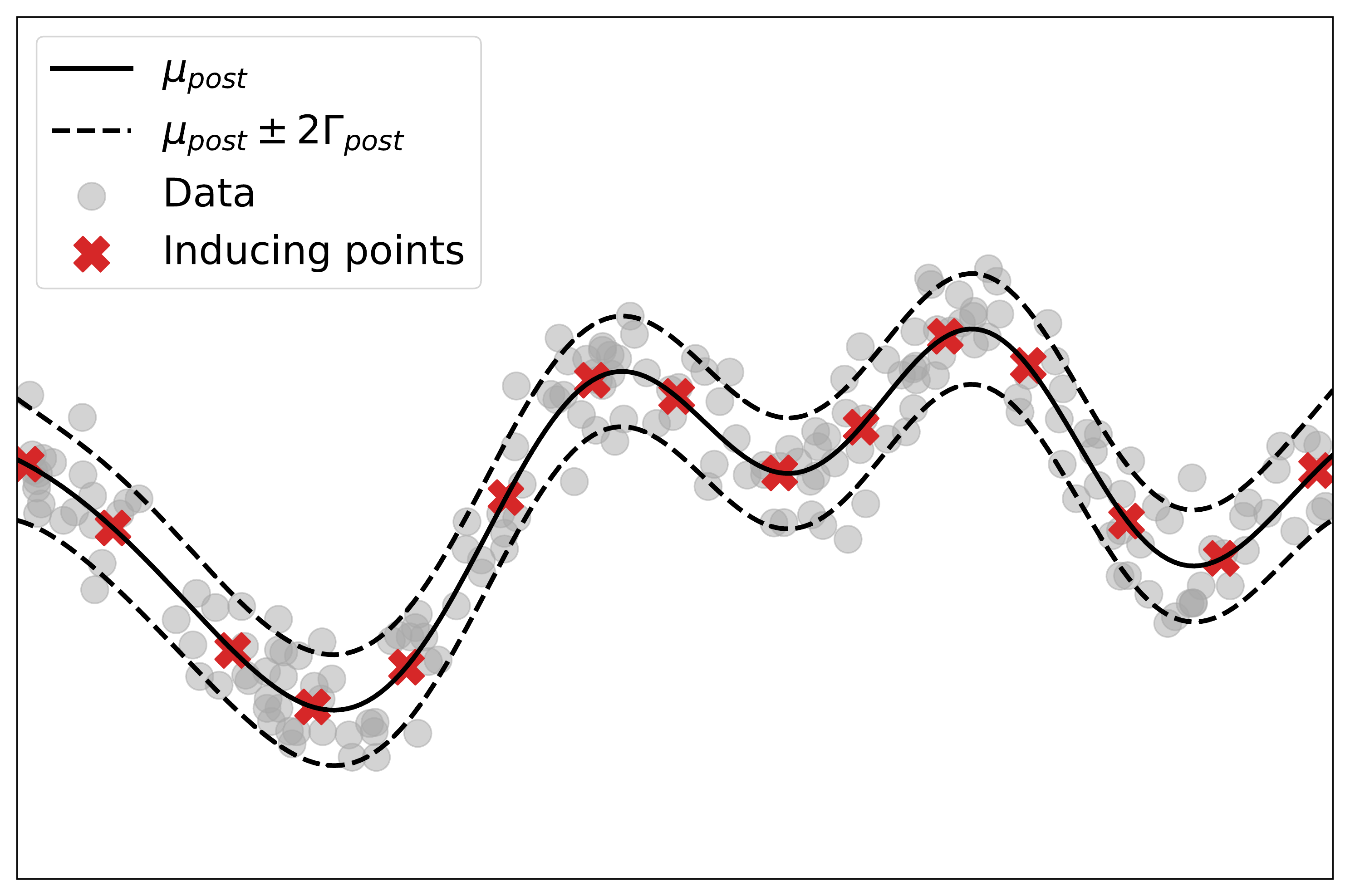}
      \includegraphics[width=.4\linewidth]{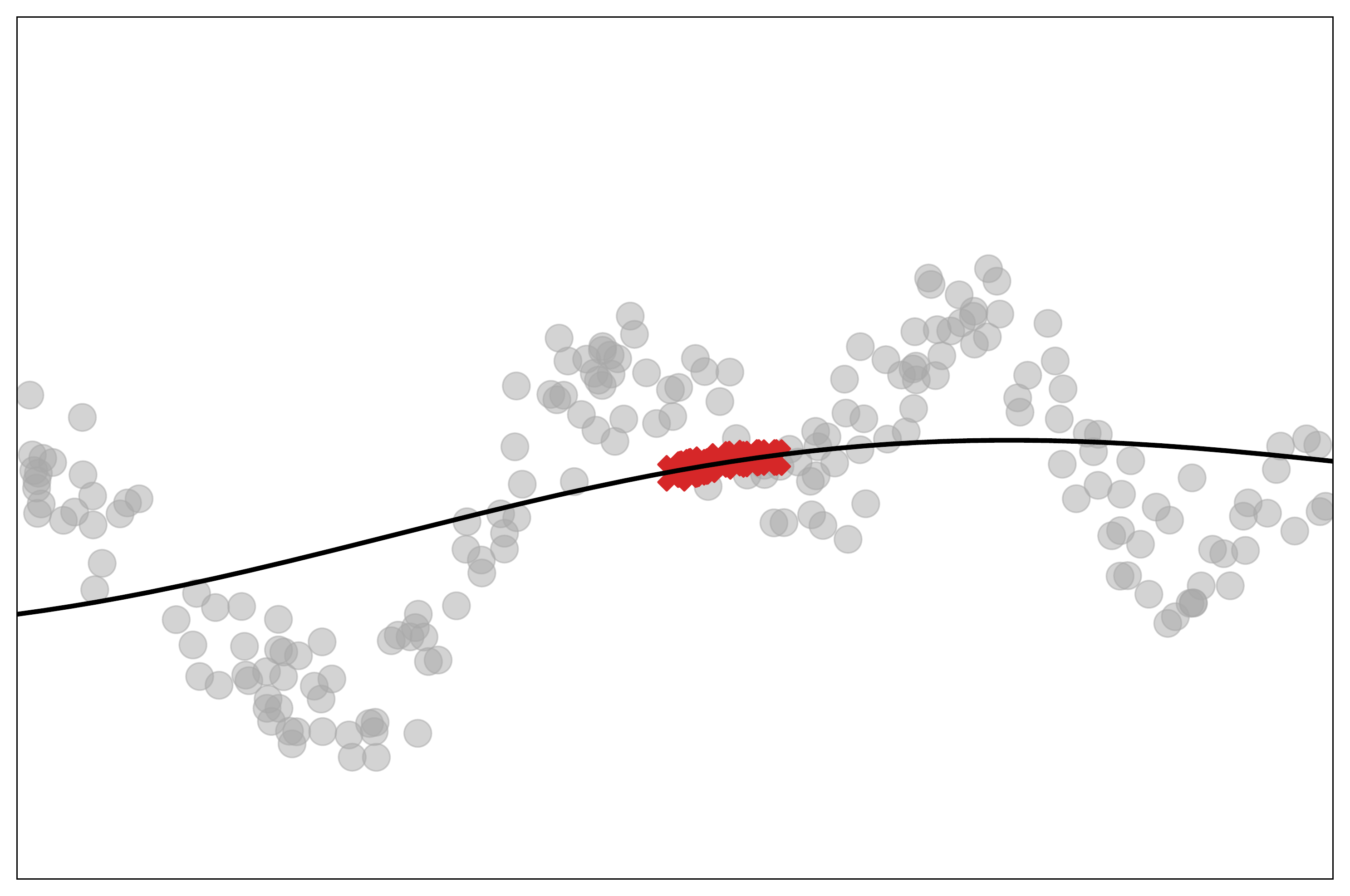}
      \vspace*{-.5em}
      \caption{Left: Gaussian process regression with optimal placement of inducing points and hyperparameter settings, determined with multiple random restarts of a BFGS solver. Right: Gaussian process regression with  initial inducing points and hyperparameter settings, chosen adversarially with the inducing points clustered together.}
      \label{fig: GPExample}
  \end{figure}

 Practically speaking, the optimization problem is $(md+\abs{\boldsymbol{\theta}}+1)$-dimensional: $md$ for $m$ inducing points in $\reals^d$, $\abs{\boldsymbol{\theta}}$ for the kernel hyperparmeters, and 1 for the noise variance. By moving to this high-dimensional optimization problem the time complexity is reduced to $\mathcal{O}(nm^2)$ and the storage costs to $\mathcal{O}(nm)$. The first term on the right-hand-side of \eqref{eqn: VariationalUpperBound} can be regarded as the loss functional, and the second term as a regularization functional penalizing a large mean-squared error. Figure \ref{fig: GPExample} is a test problem using the same data as in \cite{titsias2009variational} with fifteen inducing points. No recommendation is made in \cite{titsias2009variational} regarding which optimization routine to use. In the left panel of Figure \ref{fig: GPExample} we show the result for the placement of the inducing points and the setting of the hyperparameters, found using a BFGS solver on equation \eqref{eqn: VariationalUpperBound}.  The Gaussian process obtained with arbitrarily initialized hyperparameters and clustered inducing point locations is shown in the right pane of Figure \ref{fig: GPExample}.
  We use the variations of Algorithm \ref{alg:varprod} introduced in this paper to determine the placement of the inducing points and the settings of the hyperparameters for the same data used in \cite{titsias2009variational} and by its intellectual predecessor \cite{snelson2006sparse}. The results are compared with gradient descent and BFGS solvers. The problem is non-convex and does not satisfy the PL-inequality so the theory introduced in this paper may not provide insight into what to expect. As is common with non-convex optimization problems we perform multiple random restarts to improve the chance of converging to the global minimum. 
  
\begin{figure}[ht]
      \centering
      \includegraphics[width=.36\linewidth]{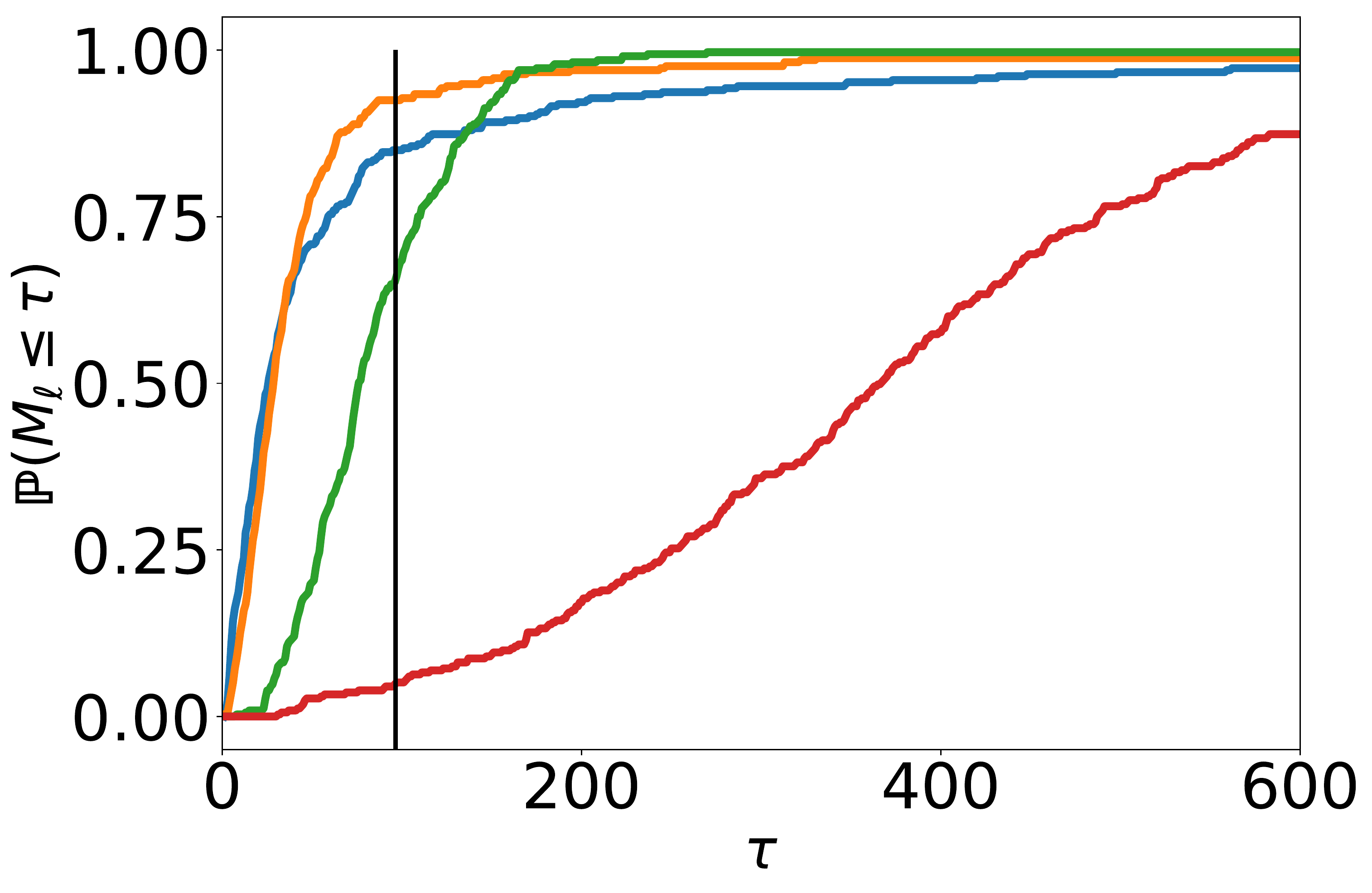}
        \includegraphics[width=.38\linewidth]{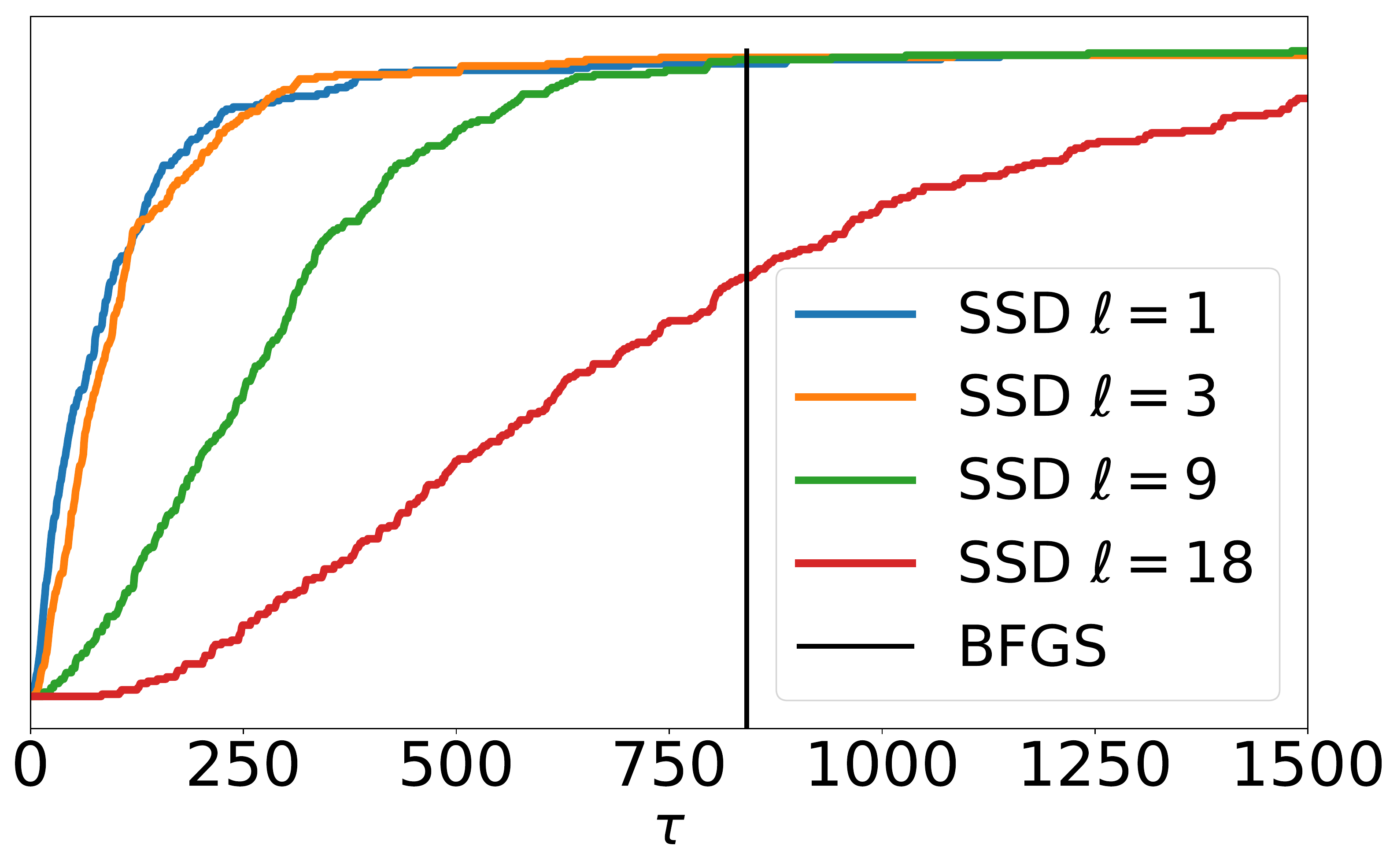}
        \vspace*{-1em}
      \caption{Left: 30-dimensional problem. Right: 60-dimensional problem. $M_{\ell}$ is the number of function evaluations required to attain a cut-off threshold for various values of $\ell$. For a fixed initialization BFGS is non-random, represented by the vertical line. Gradient descent, not pictured, has a vertical line at $\tau=2850$ and $\tau=22828$ for $p=30$ and $p=60$, respectively.}
      \label{fig: SSD-PerformanceProfile}
\end{figure}

  The data set contains 200 training points. We use the kernel $K(z_i,z_j;\boldsymbol{\theta}) = \theta_A \exp(-\abs{z_i-z_j}^2/2\theta_{l}^2)$. The parameters to be learned are the amplitude $\theta_A$, the lengthscale $\theta_l$,the noise variance $\sigma^2$, and the locations of $p-3$ inducing points in $\reals$ so that our optimization problem is $p$-dimensional. For testing purposes we can control the dimension of the optimization problem by increasing the number of inducing points. We acknowledge that for such a simple problem the use of a zeroth order method is not required as tools such as autograd \cite{maclaurin2015autograd} for Python are able to find the gradient of $\ell(\boldsymbol{\Theta};\y)$ with respect to each of the parameters by automatic differentiation; in general though, the gradient may not be available, even through automatic differentiation packages as not all functions have been incorporated into such packages. It is also possible to take individual derivatives by hand, but this is impractical in more complex, high-dimensional problems.  Once again, we use performance profiles \cite{dolan2002benchmarking} to determine the effect of varying the parameters for different problem sizes. We run each parameterization 300 times. Results for SSD for are shown in Figure \ref{fig: SSD-PerformanceProfile} for 30 and 60 inducing points. 
\begin{figure}[ht]
      \centering
      \includegraphics[width=.32\linewidth]{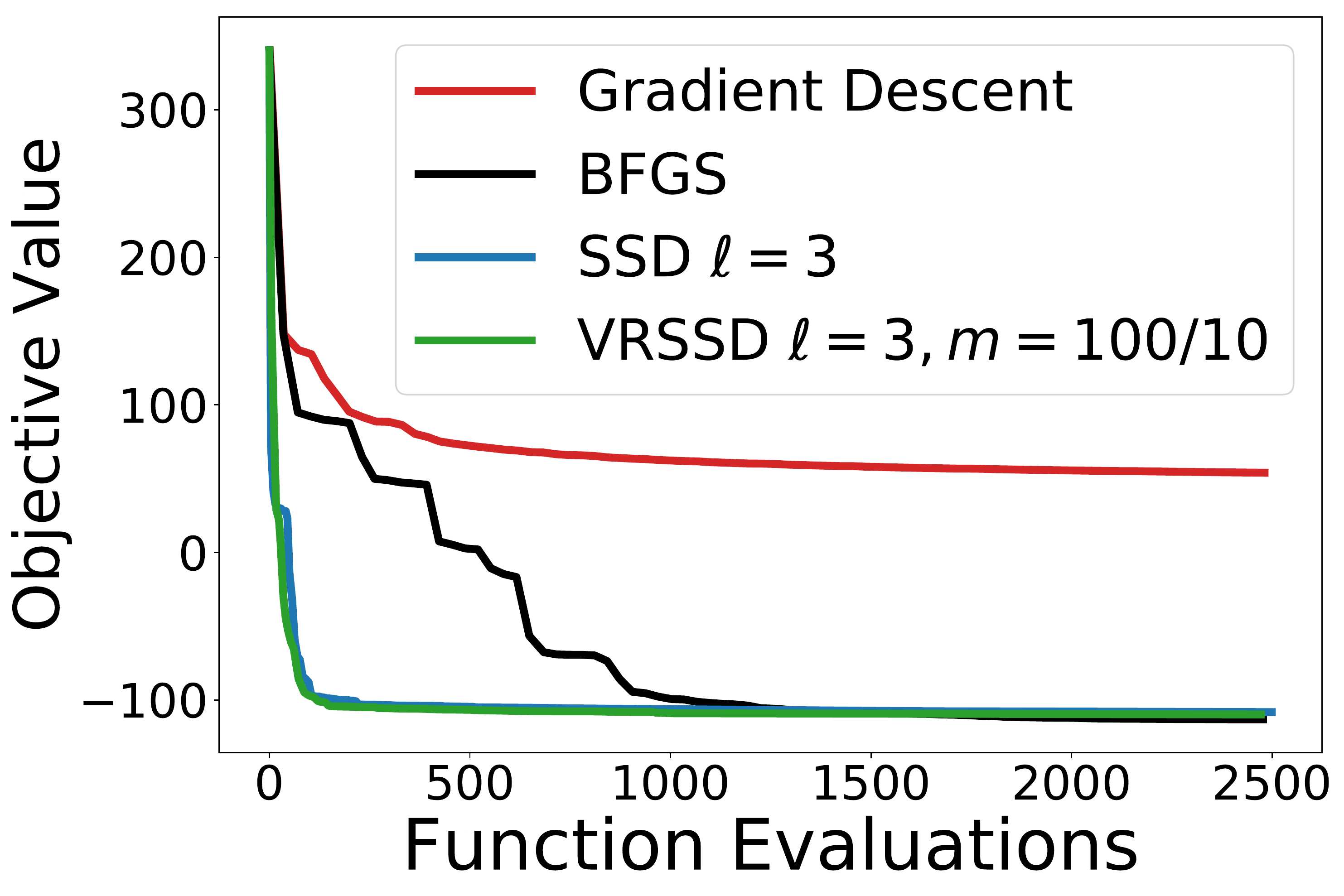}
      \includegraphics[width=.32\linewidth]{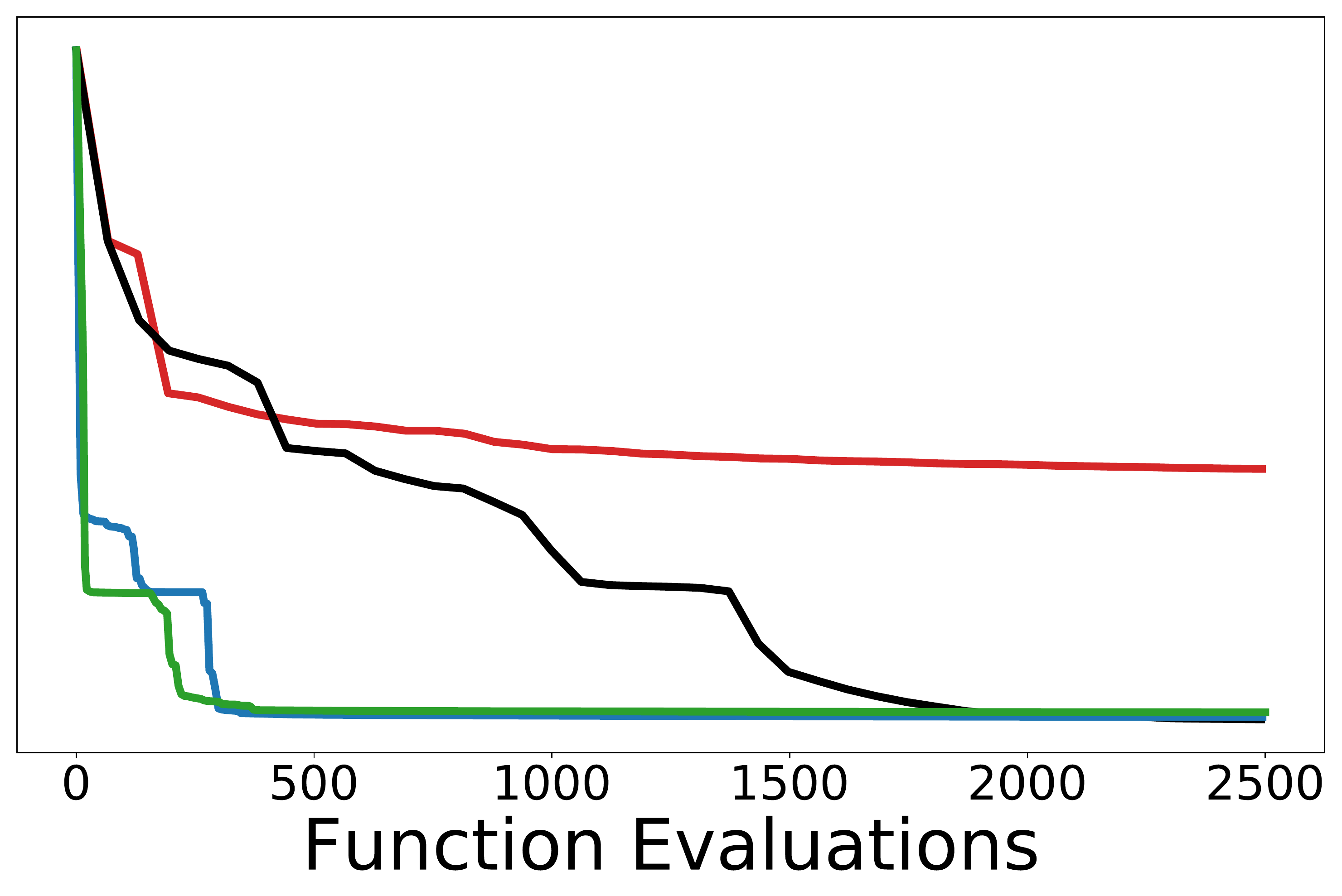}
      \includegraphics[width=.32\linewidth]{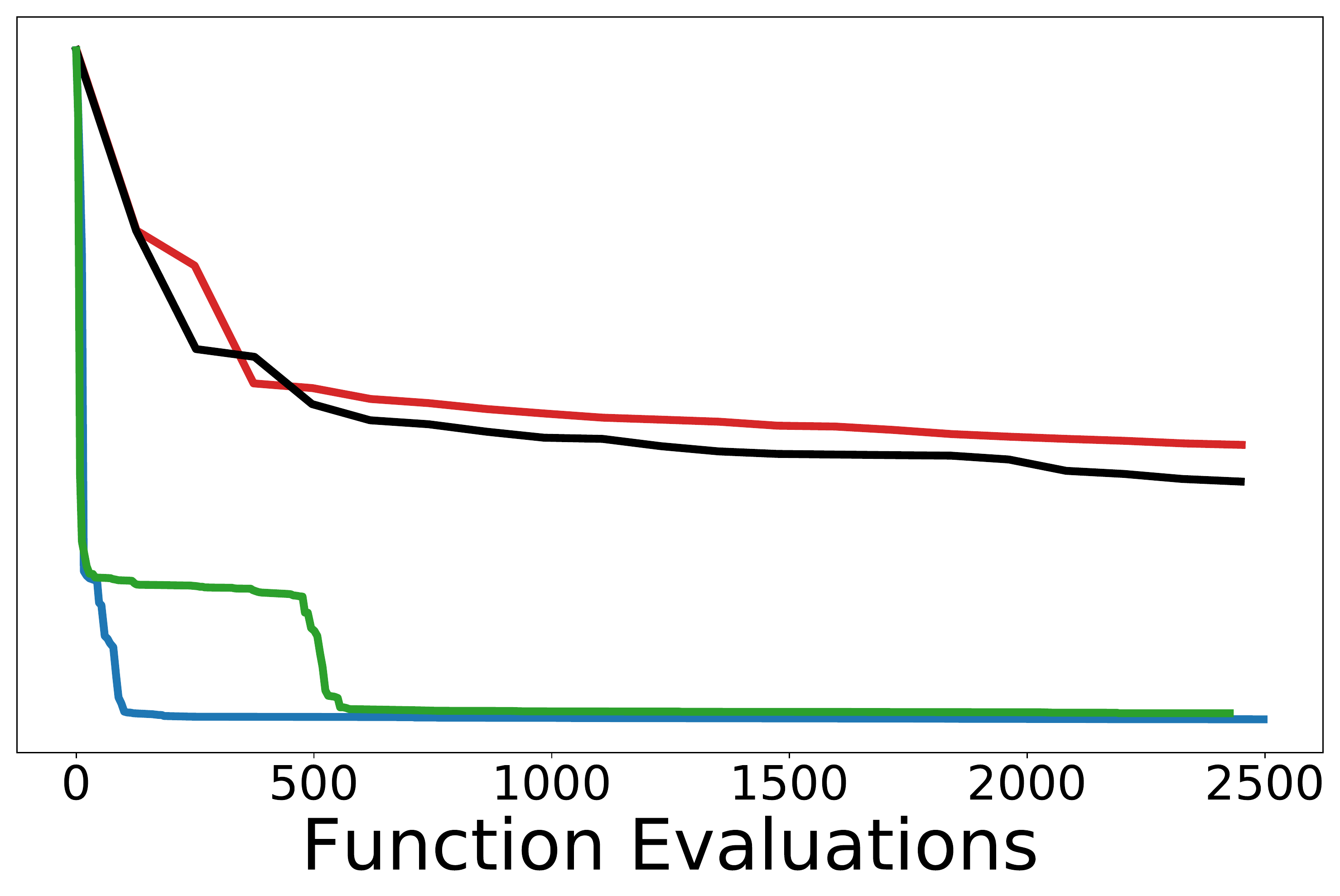}
      \vspace*{-1em}
      \caption{Left: 30-dimensional problem. Center: 60-dimensional problem. Right: 120-dimensional problem.}
      \label{fig: GP-performance}
  \end{figure}

 The cut-off threshold is $95\%$ of the distance between the objective function at the parameter initialization and at the optima, as found by BFGS. Clearly, $\ell=18$ is not a good option in this case. Similarly, $\ell=9$ can be ruled as it underperforms $\ell=1$ and $\ell=3$ approximately 90\% (resp. 99\%) of the time in the 30- (resp. 60-) dimensional problem. The case $\ell=1$ has the best single performance: in the fastest trial it is roughly 100 (resp. 800) times faster than BFGS for the 30- (resp. 60-) dimensional problem, but the variance of the performance for $\ell = 1$ is high, and about 1\% of the time it performs at least 10 times slower than BFGS (not pictured). On the other hand, $\ell=3$ beats BFGS by a similar factor and seems to be insulated from the high variance observed for $\ell=1$. Note also that in 60 dimensions $\ell=3$ is approximately three times faster than BFGS in 90\% of the trials, and about 100 times faster in 40\% of trials. A few trials of $\ell=1$ and $\ell=3$ found their way to a local minima, resulting in the methods not achieving the target threshold.

Figure \ref{fig: GP-performance} shows the performance for single trials on several problems of different dimension. For high-dimensional problems the performance of BFGS degrades because at each iteration the Hessian update is only rank two. This acts doubly against BFGS as high-dimensional problems have more expensive function evaluations, and more iterations are required before the Hessian is reasonably approximated. In contrast, the performance of SSD does not appear to degrade noticeably as the dimension increases. The performance of VRSSD is impacted by the occasional expensive full gradient evaluations but not to the same extent as BFGS. The value for the memory parameter in VRSSD is labeled $100/10$ to signify that 100 iterations of SSD are taken before any variance reduction is used, after which the memory parameter becomes 10. This is because SSD descends very quickly until it is near the solution, after which it is beneficial to make use of full gradient information.
\subsection{Shape optimization}\label{subsect: experiments-plate}

 \begin{figure}
      \centering
      \includegraphics[trim=90mm 79mm 90mm 55mm, clip, width=.45\textwidth]{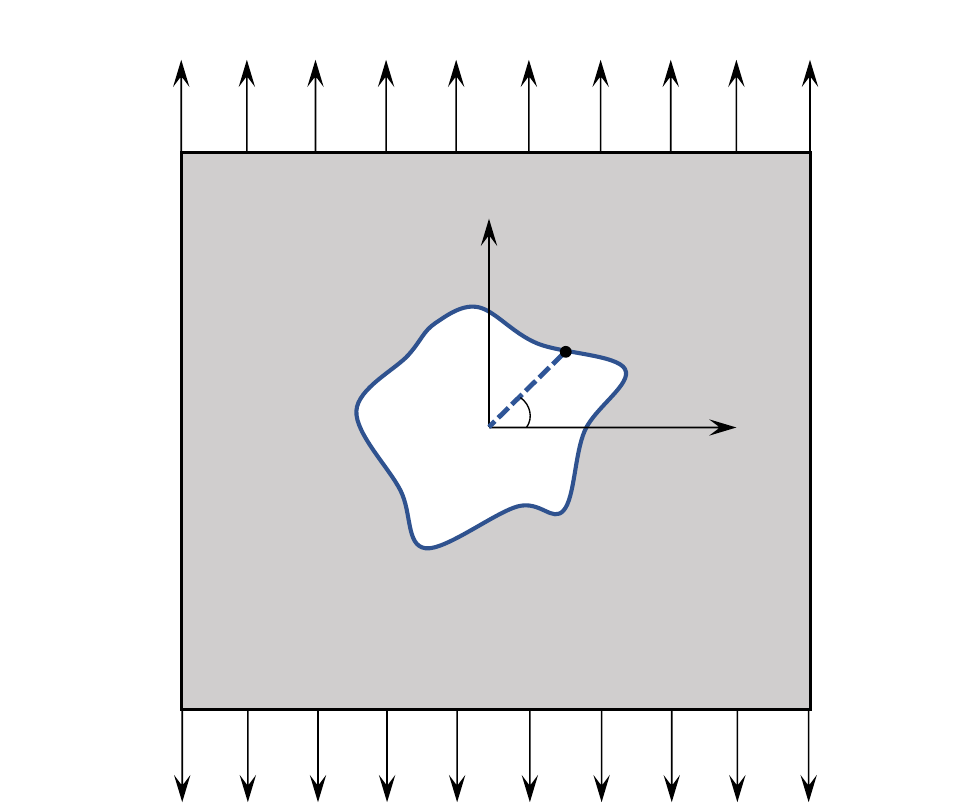}
      \includegraphics[trim=0mm 0mm 0mm 0mm, clip, width=.34\textwidth]{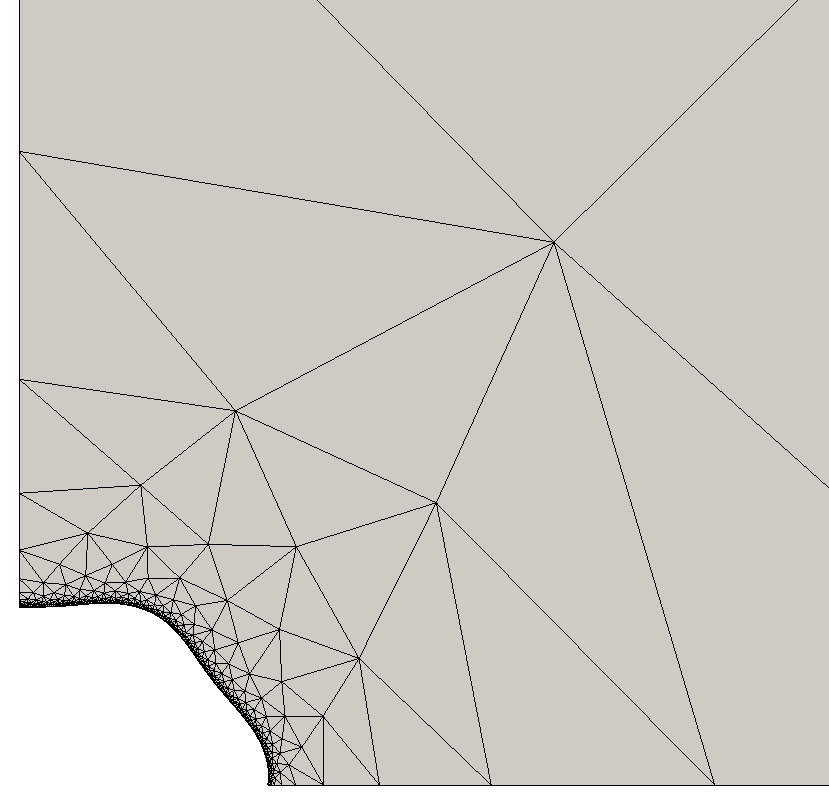}
      \put(-202,66){\footnotesize $\theta$}
      \put(-210,74){\footnotesize $r$}
      \put(-216,120){\footnotesize $\sigma_0$}
      \put(-216,0){\footnotesize $\sigma_0$}
      \put(-166,63){\footnotesize $x$}
      \put(-215,104){\footnotesize $y$}
      \vspace*{-.8em}
      \caption{Left: Schematic of the linear elasticity problem used in the shape optimization example of Section \ref{subsect: experiments-plate}. Right: Conforming finite element mesh used to solve for maximum stress $\sigma_y$ along the $y$ direction. Only a quarter of the plate corresponding to $\theta\in[0,\pi/2]$ is modeled.}
      \label{fig: HoleGeometry}
  \end{figure}
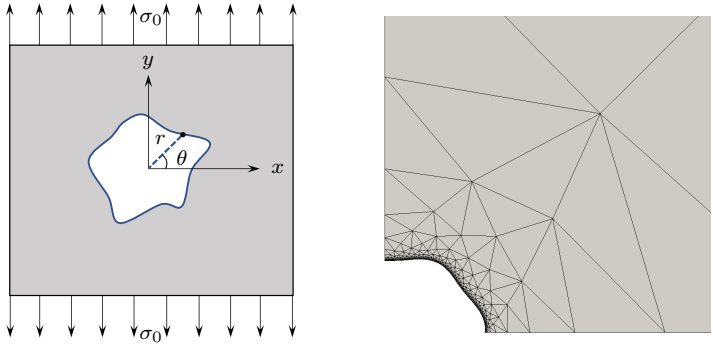

We consider a shape optimization problem involving a linear, elastic structure. Consider a square plate of size $250\times 250$ with a hole subject to uniform boundary traction $\sigma_0$=1, as illustrated in Fig.~\ref{fig: HoleGeometry}. 
We adopt a discretize-then-optimize approach to solving the PDE-constrained optimization problem. The discretization and optimization steps do not generally commute and an optimize-then-discretize approach  may be preferable for some types of problems \cite[\S 2.9]{gunzburger2003perspectives}, but we do not pursue this question here.

Our goal is to identify a shape of the hole that minimizes the maximum stress $\sigma_y$ along the $y$ direction over a quarter of the plate corresponding to $\theta\in[0,\pi/2]$. To this end, we parameterize the radius of the hole for a given $\theta$ (see Figure \ref{fig: HoleGeometry}) via
 \begin{equation}\label{eqn: parametric-radius}
     r(\theta) = 1 + \delta \medmath{\sum_{i=1}^d} i^{-1/2}\left( \xi_i \sin(i \theta) + \nu_i \cos(i \theta)\right),
 \end{equation}
 where $\delta \in (0,0.5/\sum_{i=1}^d i^{-1/2})$ is a user-defined parameter controlling the potential deviation from an n-gon of radius 1. For large values of $\delta$ it is possible to have a negative radius. The parameters that dictate the shape are $\bm{\xi} \in \reals^{d}$ and $\boldsymbol{\nu} \in \reals^{d}$. Subscripts indicate the index of the vector. We set $\delta = 0.3/\sum_{i=1}^d i^{-1/2}$ so that the minimum possible radius of any particular control point is 0.2 at the initialization. We initialize the entries of $\bm{\xi}$ and $\bm{\nu}$ uniformly at random between 0 and 1. For each instance of $\bm{\xi}$ and $\bm{\nu}$ -- equivalently $r(\theta)$ -- we generate a conforming triangular finite element mesh of the plate that we subsequently use within the FEniCS package \cite{logg2012automated} to solve for the maximum stress $\sigma_y$. A mesh refinement study is performed to ensure the spatial discretizatin errors are negligible. As we only model a quarter of the plate, we apply symmetry boundary conditions so that $y$ and $x$ displacements along $\theta=0$ and $\theta=\pi/2$ are zero. The Young's modulus and Poisson's ratio of the plate material are set to $E=1000$ and $\nu =0.3$, respectively.  
 
\begin{figure}[ht]
      \centering
      \includegraphics[width=.42\linewidth]{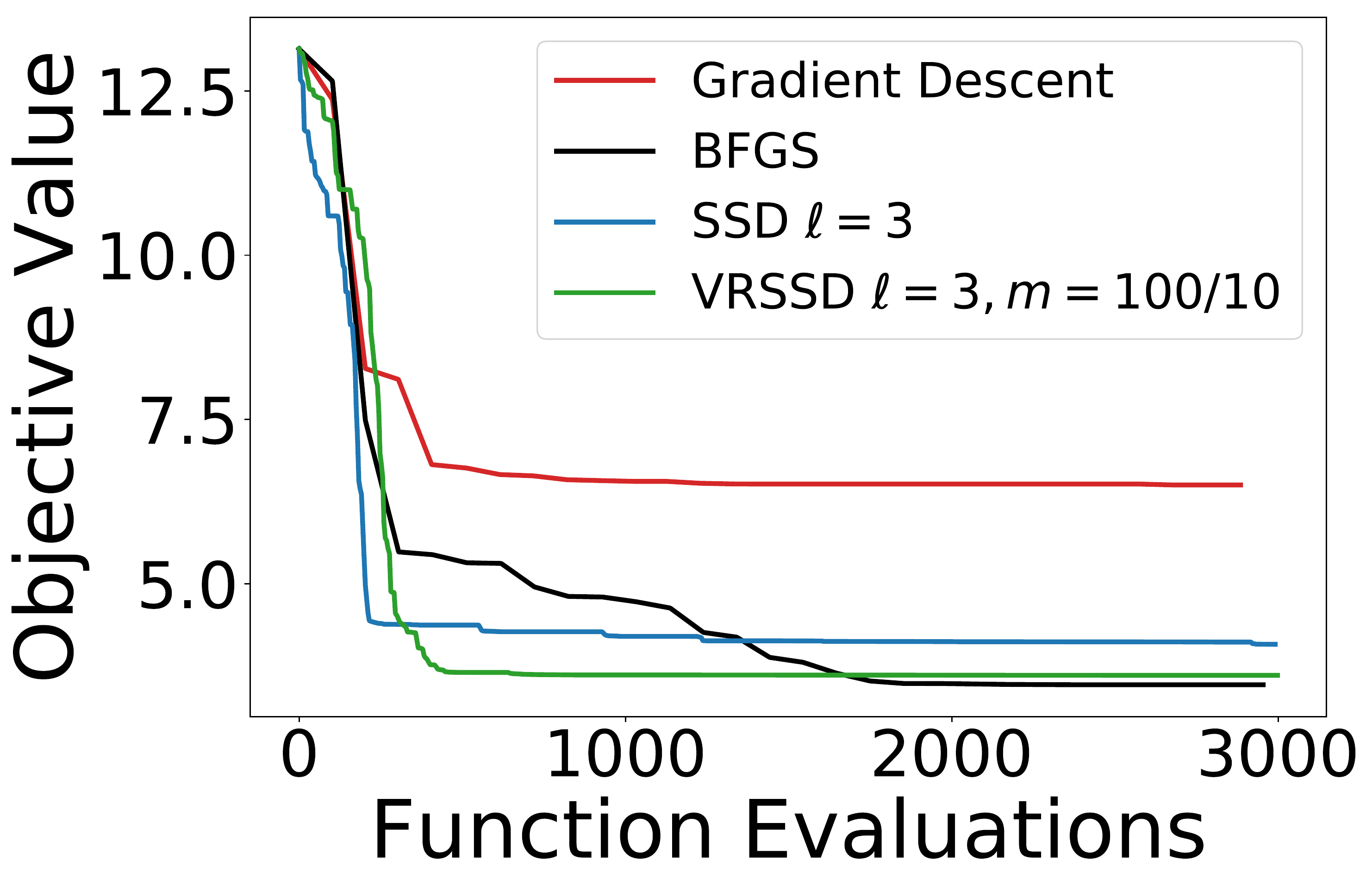}
      \includegraphics[width=.40\linewidth]{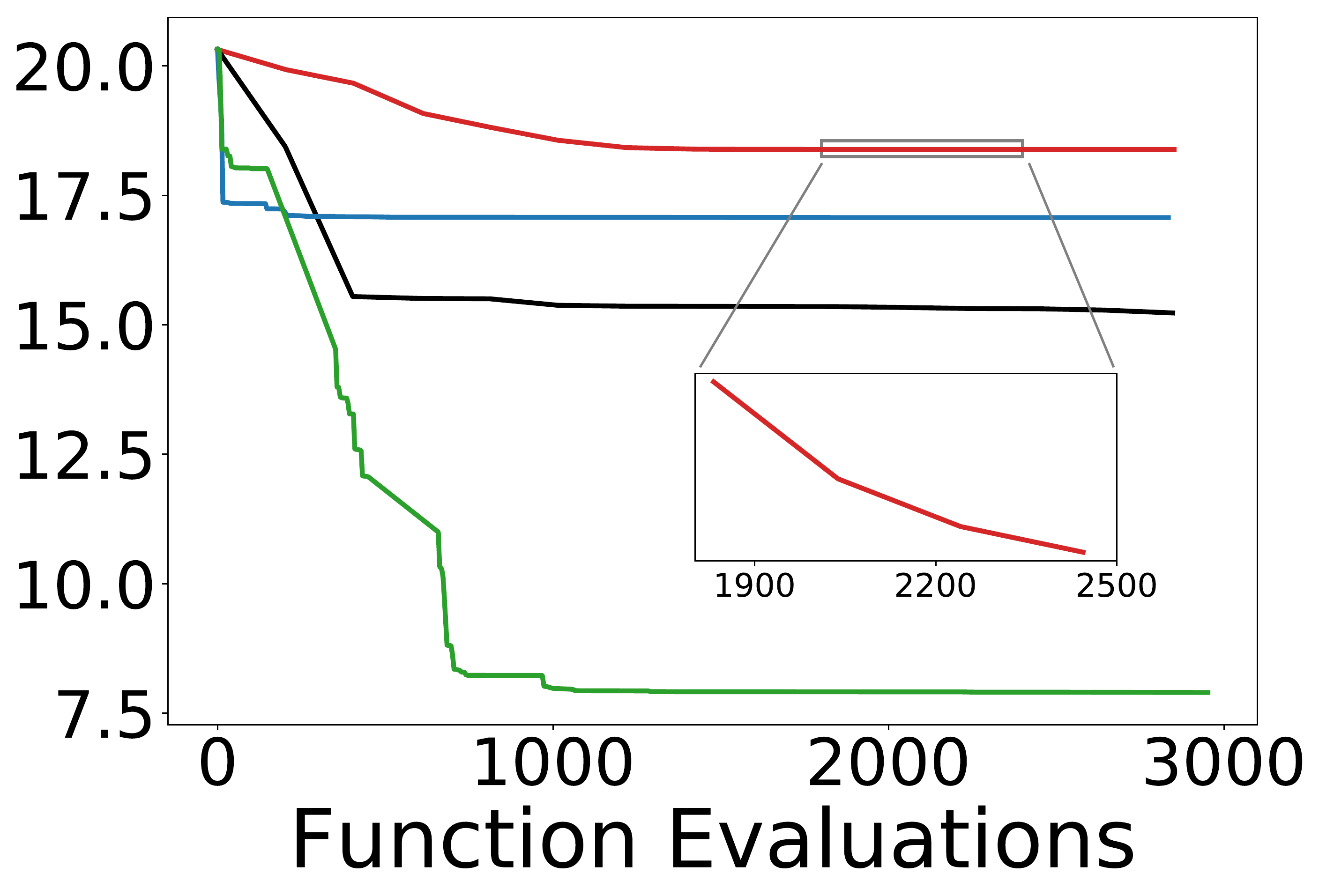}\vspace*{-1em}
      \caption{Left: 100 dimensional problem. Right: 200 dimensional problem; inset: zoom of gradient descent performance with the y-axis removed due to the minuscule absolute changes in the values.}
      \label{fig: SSD-Plate}
  \end{figure}

 The parametric radius defined by \eqref{eqn: parametric-radius} allows to scale the complexity of the problem arbitrarily by increasing the dimension $d$. In effect, if $d$ is large then the problem becomes ill-conditioned since $\xi_d$ and $\nu_d$ each make at most $\delta d^{-1/2}$ additive contribution to the radius. Such ill-posedness suggests that gradient descent ought to perform poorly as it does not account for the curvature of the objective function.  We expect BFGS to outperform gradient descent once it has iterated sufficiently to provide a good approximation of the Hessian. Based upon the intrinsic dimensionality results presented in Section \ref{subsect:synthetic-data} we anticipate SSD and VRSSD to also outperform gradient descent even though they do not explicitly account for the curvature either. Note that each function evaluation requires a PDE-solve meaning that gradient descent and BFGS require $d+1$ PDE-solves per iteration. Though a conforming finite element mesh is used to reduce the computational burden, the cost of so many PDE-solves makes this problem intractable in high-dimensions unless the resolution of the mesh is very low. On the other hand, SSD and VRSSD require far fewer PDE-solves per iteration provided $\ell \ll d$ (Figure \ref{fig: SSD-PerformanceProfile} suggests that $\ell<10$ may be reasonable). 
As mentioned above, the objective value is the maximum stress in the y-direction, $\sigma_y$, over the plate. In Figure \ref{fig: SSD-Plate} we minimize the objective for a hole with shape governed by \eqref{eqn: parametric-radius} for problems with $d=100$ and $d=200$ using gradient descent and BFGS, as well as SSD with $\ell = 3$. For VRSSD we use the same hybrid method as described in Section \ref{subsect: experiments-GP}, with identical rationale.
 As expected, BFGS outperforms gradient descent, but it is slow to provide a reasonable estimate of the Hessian and hence, relatively slow to converge. On low- or medium-dimensional problems this is not an issue, but since BFGS makes a rank-2 update to the Hessian at every iteration, and $d$-dimensional problems of the form explained in this section with radius specified by \eqref{eqn: parametric-radius} have a rank-$d$ Hessian. Thus, high-dimensional problems may require many iterations before the Hessian is well-approximated. 
 
 In the 100-dimensional problem in the left plot of Figure \ref{fig: SSD-Plate} we notice convergence of SSD and VRSSD after few function evaluations; so few, in fact, that the difference in performance between SSD and VRSSD likely has more to do with random chance than variance reduction. In multiple restarts (not shown), similar performance was observed wherein the methods converge before any variance reduction occurs. In all cases, convergence occurs by the time gradient descent of BFGS are able to take more than a couple of steps.
 
 In the 200-dimensional problem in the right plot of Figure \ref{fig: SSD-Plate} the progress of gradient descent, SSD, and BFGS stall prematurely, appearing to settle into local optima at different locations despite being initialized identically. 
 The figure is somewhat deceptive as each method is indeed still descending along shallow valleys of the objective function as shown in the inset. 
 VRSSD does not appear to be as susceptible to getting caught in such valleys for reasons we do not yet fully understand. 
 We hypothesize that there are two phenomena at play. First, as the $\ell$-dimensional subspace along which SSD and VRSSD descend changes with each iteration, they explore parameter space more erratically than deterministic methods, making them less likely to get funnelled into long, shallow basins; this is intuitively similar to the recent line of research which suggests that noisy perturbation of iterative algorithms helps avoid saddle points \cite{ge2015escaping}. 
 Second, VRSSD uses past information from a full-gradient evaluation that push it in directions that may not be optimal from a single-step perspective but that may help it traverse these valleys in a way that is heuristically similar to the heavy ball method of Polyak \cite{polyak1964some}; this phenomena has recently been explored in the context of non-convex minimization, e.g., \cite{jin2017accelerated}. Since SSD still appears to get caught in these low-gradient valleys, it is apparently the combination of these two reasons, the second of them, or some as-yet-undescribed property of VRSSD that allows it to avoid them. 
 It could also be the case that this particular problem is amenable to VRSSD and further investigation on different problems is required to understand why VRSSD appears to perform so well. 

\section{Conclusions}\label{sect: Conclusions}
We have presented a generalization of coordinate descent and gradient descent algorithms, and have shown asymptotic convergence at a linear rate under various convexity assumptions. We have demonstrated empirical improvements compared to the \emph{status quo} for several practical problems that do not necessarily adhere to the assumptions required by the theory. We have provided analysis that puts the ideas of the SVRG algorithm \cite{SVRG} in a stochastic subspace framework, making it a feasible alternative to deterministic gradient-based methods for problems that cannot be represented with a finite-sum structure. This analysis opens the doors to adapting other stochastic iterative methods designed for an ERM framework into a more general form.   

In the future we plan to adapt other first-order stochastic optimization techniques such as stochastic quasi-Newton \cite{Mokhtari2015GlobalCO} to the stochastic subspace descent framework. We plan to perform analysis to determine what properties can be expected as a result of using different random matrices (e.g., Haar-distributed, Gaussian, Rademacher) for the subspace projection. A clever adaptive scheme could be beneficial for determining the decent directions, an idea that has been discussed at length in the coordinate descent literature \cite{richtarik2014iteration, nesterov2012efficiency}. Parallelizing our methods to calculate the $\ell$ directional derivatives at each iteration simultaneously is straight-forward, but we would like to explore the feasibility of asynchronous parallelization as has been discussed in the coordinate descent case (see, e.g., \cite{peng2016arock}). A computationally straightforward extension may allow sketching methods (see e.g. \cite{pilanci2015randomized}) to improve our results with minimal programmatic overhead, but analysis must be conducted to confirm the theoretical properties of such modifications. The analysis we have performed thus far does not explain some phenomena we have observed in our experiments; in particular, it appears as though the intrinsic dimension rather than the ambient dimension dictates the convergence rate of our methods, and that VRSSD is resistant to getting trapped in the shallow basins common to ill-conditioned non-convex problems. 

Finally, recent work on a universal ``catalyst'' scheme \cite{lin2015universal} also applies to our method, allowing for Nesterov-style acceleration without requiring additional knowledge of the Lipschitz constants along any particular direction.

\appendix
\section{Proofs of main results}
\subsection*{Proof of Theorem \ref{thm:convergence}} Define the filtration (that is, sequence of $\sigma$-algebras) $\mathcal{F}_k = \sigma(\mathbf{P}_1, \ldots, \mathbf{P}_{k-1}), ~ k>1$, and $\mathcal{F}_1 = \{ \emptyset , \Omega\}$.
	Because $f$ is continuously-differentiable with a $\lambda$-Lipschitz derivative, it follows that  
	\begin{equation}\label{Taylor}
	f(\x_{k+1}) \leq  f(\x_k) + \nabla f(\x_k)^\top (\x_{k+1}-\x_k) + \frac{\lambda}{2}\norm{\x_{k+1}-\x_k}^2.
	\end{equation}
Let $f_e(\x) = f(\x)-f_*$ be the error for a particular $\x$.   Then, \eqref{eq: iterations} and \eqref{Taylor} yield:
	\begin{align}\label{delineq}
	\nonumber	f_e(\x_{k+1}) - f_e(\x_k)  &\leq \alpha_{\lambda} \langle \nabla f(\x_k),~ \mathbf{P}_k\mathbf{P}_k^\top \nabla f(\x_k) \rangle \quad \text{with} \quad \alpha_{\lambda} = \alpha + d\alpha^2\lambda/2\ell,
	\end{align}
	where we have used the fact that $\mathbf{P}_k\mathbf{P}_k^\top \mathbf{P}_k\mathbf{P}_k^\top = (d/\ell)\mathbf{P}_k\mathbf{P}_k^\top $. We can choose  $0< \alpha < 2\ell/d\lambda$ so that $\alpha_{\lambda} < 0$. With this choice the right-hand side is  non-positive and the errors are non-increasing. Since the error is bounded below by zero the sequence converges almost surely. Furthermore, since the sequence is bounded above by $f_e(\x_0)$, Lebesgue's dominated convergence implies convergence of the sequence in $L^1$. To find the actual limit we take conditional expectations of both sides to get
	\begin{equation*}
	\Expectation [f_e(\x_{k+1}) \mid \mathcal{F}_k ] \leq \alpha_{\lambda} \Expectation [ \langle \nabla f(\x_k), \mathbf{P}_k\mathbf{P}_k^\top  \nabla f(\x_k) \rangle  \mid \mathcal{F}_k ] +f_e(\x_k),
	\end{equation*}
    which leads to
	\begin{equation} \label{intermediate}
	\Expectation \left[f_e(\x_{k+1}) \mid \mathcal{F}_k\right] \leq  \alpha_{\lambda} \norm{\nabla f(\x_k)}^2 + f_e(\x_{k}),
	\end{equation}
and since $\alpha_{\lambda} < 0$ the PL-inequality yields
	\begin{equation*}
	\Expectation \left[f_e(\x_{k+1}) \mid \mathcal{F}_k\right] \leq 2\gamma  \alpha_{\lambda} f_e(\x_k) + f_e(\x_{k})= \left(1+2\gamma\alpha_{\lambda}\right)f_e(\x_k).
	\end{equation*}
	With recursive application we obtain $\Expectation \left[f(\x_{k+1})-f_* \mid \mathcal{F}_k \right]\leq  \left(1-2\gamma\alpha_{\lambda}\right)^{k+1} \left(f(\x_0) - f_*\right)$.
	Taking an expectation of both sides yields
	\begin{equation*}
	    	\Expectation f(\x_{k+1})-f_* \leq  \left(1+2\gamma\alpha_{\lambda}\right)^{k+1} \left(f(\x_0) - f_*\right).
	\end{equation*}
	Any step-size $0 < \alpha < 2\ell/d\lambda$ forces the right-hand side to converge to zero. Thus, since $	f_e(\x_k) \convas \X$ for some $\X \in L^1$ and $f_e(\x_k)\overset{L^1}{\to} 0$, we have both $f(\x_k) \convas f_*$ and $f(\x_k) \overset{L^1}{\to} f_*$. 
\qed
\subsection*{Proof of Corollary \ref{corr:strong-convexity}  part (i)}
By strong-convexity, the PL-inequality, and Theorem \ref{thm:convergence} we obtain $f(\x_k) \convas f(\x_*)$ and $f(\x_k) - f(\x_*) \geq \frac{\gamma}{2} \norm{\x_*-\x_k}$. Since the left-hand side converges a.s.\ to zero and $\gamma > 0$, we have $\x_k \convas \x_*$.
\subsection*{Proof of Corollary \ref{corr:strong-convexity}  part (ii)}
 Rearranging the terms in equation \eqref{intermediate} we have $\alpha_{\lambda} ^{-1}\Expectation \left[f_e(\x_k) - f_e(\x_{k+1})\mid \mathcal{F}_k \right] \geq \norm{\nabla f(\x_k)}^2$. Combining this with Lipschitz continuity yields
	$
	2\gamma f_e(\x_k) \leq \norm{\nabla f(\x_k)}^2 \leq \alpha_{\lambda}^{-1} \Expectation\left[f_e(\x_{k}) - f_e(\x_{k+1}) \mid \mathcal{F}_k\right]. $
	That is,
	\begin{equation}\label{nearly}
	\Expectation \left[f_e(\x_{k+1}) \mid \mathcal{F}_k\right] \leq \left(1+2\gamma\alpha_{\lambda}\right)f_e(\x_k).
	\end{equation}
	With repeated recursion	we get $\Expectation \left[f_e(\x_{k+1}) \mid \mathcal{F}_k\right] \leq \left(1+2\gamma\alpha_{\lambda}\right)^{k+1}f_e(\x_0)$, the right-hand side of which achieves a minimum when $\alpha=\ell/d\lambda$, resulting in $\Expectation \left[f_e(\x_{k+1}) \mid \mathcal{F}_k\right] \leq \left(1-\ell\gamma / d \lambda \right)^{k+1}f_e(\x_0)$. An expectation of both sides proves the corollary.
\qed

	\subsection*{Proof of Theorem \ref{thm: convergence-convex}}
	We follow the proof of Theorem \ref{thm:convergence} until \eqref{intermediate}, then we rearrange terms to obtain,
	\begin{align}\label{eq:simplified-Lipschitz}
	\Expectation \left[f(\x_{k+1}) \mid \mathcal{F}_k \right] \leq f(\x_k) + \alpha_{\lambda} \norm{\nabla f(\x_k)}^2,
	\end{align}
	and then by convexity, $\norm{\nabla f(\x_k)} \geq f_e(\x_k)/R$. Plugging this into equation \eqref{eq:simplified-Lipschitz} and letting $\alpha = \ell/d\lambda$ results in 
	\begin{equation}\label{eq: required}
 \Expectation [f_e(\x_{k+1}) \mid \mathcal{F}_k]- f_e(\x_k)  \leq  -\alpha  f_e(\x_k)^2/2 R^2,
	\end{equation}
	and one more expectation yields
	\begin{equation*}
	    \Expectation[f_e(\x_{k+1}) - f_e(\x_k)] \leq -\alpha \Expectation f_e(\x_k)^2/2R^2
	    \leq -\alpha \left(\Expectation f_e(\x_k)\right)^2/2R^2.
	\end{equation*}
    We note that $\alpha/2 R^2 \geq 1/\Expectation f_e(\x_k) - 1/\Expectation f_e(\x_{k+1})$,
	because $\Expectation f_e(\x_{k+1}) < \Expectation f_e(\x_k)$. Then \eqref{eq:simplified-Lipschitz} leads to
	\begin{equation}\label{eq: before recursion}
	\frac{1}{\Expectation f_e(\x_{k+1})} \geq 	\frac{1}{\Expectation f_e(\x_{k})} + \frac{\alpha}{2 R^2}.
	\end{equation}
	Applying \eqref{eq: before recursion} recursively, and replacing $\alpha$ with $\ell/d\lambda$ we obtain
	$\Expectation  f_e(\x_{k+1}) \leq 2d\lambda R^2/k \ell$. 
\qed

\subsection*{Proof of Lemma \ref{thm:SVRG-AS}}
By Markov's inequality $\mathbb{P}(\norm{f(\x_k) - f(\x_*)} > \epsilon) \leq C\beta^k$.	Sum both sides to get $\sum_{k=1}^{\infty}	\mathbb{P}(\norm{f(\x_k) - f(\x_*)} > \epsilon) \leq C/(1- \beta)$.
Therefore $f(\x_n) \convas f(\x_*)$. For the convergence of the iterates, note that if $f$ is strongly convex it has a unique minimizer.
	\qed

\subsection*{Proof of Theorem \ref{thm:VarianceReducedRandomGradient} part (i)}
We define the following filtration that encodes information at each step $s,k$:
\begin{equation*}
\mathcal{F}_{s,k} = \sigma(\textbf{P}_{1,1}, \ldots, \textbf{P}_{1,m}, \ldots, \textbf{P}_{s-1,1},\ldots,  \textbf{P}_{s-1,m}, \textbf{P}_{s,1},\ldots,  \textbf{P}_{s,k-1}  , J_1, \ldots J_{s-1}),
\end{equation*}
which makes $\x_{s,k}$ and $\tilde{\x}_{s-1}$ $\mathcal{F}_{s,k}$-measurable. Note that 
\begin{equation}\label{doubleexpect}
\Expectation f(\tilde{\x}_s) = \Expectation~ \Expectation (f(\tilde{\x}_s)\mid \mathcal{F}_{s,m}) = \Expectation (1/m)\medmath{\sum_{k=0}^{m-1}} f(\x_{s,k}).
\end{equation}
The recursion for epoch $s$ can be written as 
	\begin{equation}\label{eqn: t4-recursion}
	\x_{s,k} = \x_{s,k-1} - \alpha~ \v,
	\end{equation}
	where $\v  =  \P\P^\top \nabla f(\x_{s,k-1}) -  \left(\P\P^\top -\I\right)\nabla f(\tilde{\x}_{s-1})$.
   Equation \eqref{eqn: t4-recursion} then leads to
	\begin{equation}\label{eqn: t4-normsquared}
	\norm{\x_{s,k}-\x_*}^2 = \norm{\x_{s,k-1}-\x_*}^2 - 2\alpha (\x_{s,k-1}-\x_*)^\top \v + \alpha^2 \norm{\v}^2,
	\end{equation}
which in turn leads to
	\begin{equation}\label{eqn: t4-expected normsquared}
\CondExp{\norm{\x_{s,k}-\x_*}^2} =\norm{\x_{s,k-1}-\x_*}^2 - 2\alpha (\x_{s,k-1}-\x_*)^\top \CondExp{\v} + \alpha^2 \CondExp{\norm{\v}^2}.
	\end{equation}
By definition, $\CondExp{\v} = \nabla f(\x_{s,k-1})$. For the last term we have:
\begin{align*}
     \CondExp{\norm{\v}^2} = \rho \norm{\nabla f(\x_{s,k-1})}^2 + (\rho-1) \norm{f(\tilde{\x}_{s-1})}^2 - 2(\rho-1) \nabla f(\x_{s,k-1})^\top f(\tilde{\x}_{s-1}).
\end{align*}
Plugging this into \eqref{eqn: t4-expected normsquared}, taking another expectation, summing from $k=1$ to $k=m$, and using convexity yields
\begin{align}\label{eqn: t4-second-intermediate-step-after-sum}
  \nonumber   \Expectation \left[\norm{\x_{s,m} - \x_*}^2 \right] \leq&   \Expectation \left[\norm{\x_{s,0} - \x_*}^2 \right] +2\alpha^2 \lambda m(\rho-1)f_e(\tilde{\x}_{s-1})  -2\alpha(1- \alpha \lambda\rho)  \medmath{\sum_{k=1}^m} \Expectation f_e(\x_{s,k-1})
 \\ &-  2\alpha^2(\rho-1)\medmath{\sum_{k=1}^m} \Expectation \left[  \nabla f(\x_{s,k-1})^\top f(\tilde{\x}_{s-1})\right].
\end{align}
By $\gamma$-strong convexity, $\norm{\x_{s,0}-\x_*}^2 \leq 2/\gamma \left(f(\x_{s,0})-f(\x_*)\right)$, and then \eqref{eqn: t4-second-intermediate-step-after-sum} leads to
\begin{align}\label{eqn: t4-penultimate step}
  \nonumber  0 \leq{}&  \frac{2}{\gamma}(1 + \alpha^2 \gamma \lambda m(\rho-1)) \Expectation f_e(\tilde{\x}_{s-1}) -2\alpha(1- \alpha \lambda\rho)  \medmath{\sum_{k=1}^m} \Expectation f_e(\x_{s,k-1}).
\end{align}
By \eqref{doubleexpect}, $\Expectation \sum_{k=1}^{m} f(\x_{s,k-1})=m\Expectation f(\tilde{\x}_s)$. Thus, we have
\begin{equation*}
 \Expectation f_e(\tilde{\x}_s)\leq \beta f_e(\tilde{\x}_{s-1}), \quad \beta =  \frac{1}{\alpha \gamma m (1 - \alpha \lambda \rho)} + \frac{\alpha \lambda (\rho-1)}{1- \alpha \lambda \rho}.
\end{equation*}
	 By recursion we get the desired result. The proof of a.s. convergence follows as in the proof of Theorem \ref{thm:SVRG-AS}.
\qed

\subsection*{Proof of Theorem \ref{thm:VarianceReducedRandomGradient} part (ii)}

The proof of part (ii) of Theorem \ref{thm:VarianceReducedRandomGradient} follows closely that of part (i), but $\eta_{s,k-1}$ must be accounted for. This time, the recursion for epoch $s$ is written as $\x_{s,k} = \x_{s,k-1} - \alpha \v$,
	with 
	\begin{equation*}
	\v  =  \P\P^\top \nabla f(\x_{s,k-1}) -  \eta_{s,k-1} \left(\P\P^\top -\I\right)\nabla f(\tilde{\x}_{s-1}).
	\end{equation*}
    By subtracting $\x_*$ from both sides of the recursion we get
	\begin{equation}\label{eqn: normsquared}
	\norm{\x_{s,k}-\x_*}^2 = \norm{\x_{s,k-1}-\x_*}^2 - 2\alpha (\x_{s,k-1}-\x_*)^\top \v + \alpha^2 \norm{\v}^2,
	\end{equation}
which leads to
	\begin{equation}\label{eqn: expected normsquared}
\CondExp{\norm{\x_{s,k}-\x_*}^2} =\norm{\x_{s,k-1}-\x_*}^2 - 2\alpha (\x_{s,k-1}-\x_*)^\top \CondExp{\v} + \alpha^2 \CondExp{\norm{\v}^2}.
	\end{equation}
Again, $\CondExp{\v} = \nabla f(\x_{s,k-1})$, but this time
\begin{align}\label{eqn: v-squared}
    \nonumber \CondExp{\norm{\v}^2} =& \rho \norm{\nabla f(\x_{s,k-1})}^2 + (\rho-1) \eta_{s,k-1}^2 \norm{f(\tilde{\x}_{s-1})}^2\\ &- 2(\rho-1) \eta_{s,k-1}\nabla f(\x_{s,k-1})^\top f(\tilde{\x}_{s-1})
\end{align}
Plug equation \eqref{eqn: v-squared} into equation \eqref{eqn: expected normsquared} to get
	\begin{align*}
\nonumber \CondExp{\norm{\x_{s,k}-\x_*}^2} ={}&\norm{\x_{s,k-1}-\x_*}^2 - 2\alpha (\x_{s,k-1}-\x_*)^\top \nabla f(\x_{s,k-1}) \\
&+ \alpha^2 
\left(\rho \norm{\nabla f(\x_{s,k-1})}^2 + (\rho-1) \eta_{s,k-1}^2 \norm{f(\tilde{\x}_{s-1})}^2 \right. \\& \left.- 2(\rho-1) \eta_{s,k-1}\nabla f(\x_{s,k-1})^\top f(\tilde{\x}_{s-1})
\right).
	\end{align*}
The value that minimizes the quadratic function of $\eta_{s,k-1}$ in the parentheses is
\begin{equation}\label{eqn: eta-opt}
    \eta_{s,k-1} = \nabla f(\x_{s,k-1})^\top \nabla f(\tilde{\x}_{s-1})/\norm{f(\tilde{\x}_{s-1})}^2,
\end{equation}
which by convexity leads to 
\begin{align}\label{eqn: second-intermediate-step-after-sum}
  \nonumber   \Expectation \left[\norm{\x_{s,m} - \x_*}^2 \right] \leq{}&   \Expectation \left[\norm{\x_{s,0} - \x_*}^2 \right] -2\alpha (1 - 2\alpha\lambda\rho)  \medmath{\sum_{k=1}^m} \Expectation f_e(\x_{s,k-1}) \\
    &- \alpha^2(\rho-1) \medmath{\sum_{k=1}^m} \Expectation\left[\left( \nabla f(\x_{s,k-1})^\top \nabla f(\tilde{\x}_{s-1})/\norm{f(\tilde{\x}_{s-1})} \right)^2\right].
\end{align}
By $\gamma$-strong convexity, $\norm{\x_{s,0}-\x_*}^2 \leq 2/\gamma \left(f(\x_{s,0})-f(\x_*)\right)$. We obtain
\begin{align}\label{eqn: penultimate step}
  \nonumber  0 \leq{}&  \frac{2}{\gamma} \Expectation f_e(\x_{s,0})  -2\alpha (1 - 2\alpha\lambda\rho)  \medmath{\sum_{k=1}^m} \Expectation f_e(\x_{s,k-1}).
\end{align}
Further simplification results from equation \eqref{doubleexpect} which suggests that \\ $\Expectation \sum_{k=1}^{m} f(\x_{s,k-1})=m\Expectation f(\tilde{\x}_s)$. We move this term to the left-hand side so that
\begin{equation*}
 \Expectation f_e(\tilde{\x}_s) - f(\x_*) \leq \beta \Expectation f_e(\tilde{\x}_{s-1}), \quad 	\beta = \frac{1}{\alpha \gamma m (1 - \alpha \lambda \rho)}.
\end{equation*}
	By recursion we get the desired result. The proof of a.s. convergence follows in the same manner as in the proof of Theorem \ref{thm:SVRG-AS}.
	
\section*{Acknowledgments}
We thank Lior Horesh for his insightful suggestions regarding derivative free optimization, and Gregory Fasshauer for fruitful discussions on sparse Gaussian processes. We also thank Subhayan De for providing the FEniCS code used in Section \ref{subsect: experiments-plate}.

\bibliographystyle{siamplain}
\bibliography{RandomizedBib,thesisBecker}

\begin{thebibliography}{10}

\bibitem{Abacioglu2001}
{\sc Y.~Abacioglu, D.~Oliver, and A.~Reynolds}, {\em Efficient reservoir
  history matching using subspace vectors}, Computat. Geosci., 5 (2001),
  pp.~151--172.

\bibitem{allen2016even}
{\sc Z.~Allen-Zhu, Z.~Qu, P.~Richt{\'a}rik, and Y.~Yuan}, {\em Even faster
  accelerated coordinate descent using non-uniform sampling}, in ICML, 2016,
  pp.~1110--1119.

\bibitem{Bell2007LessonsFT}
{\sc R.~M. Bell and Y.~Koren}, {\em {Lessons from the Netflix prize
  challenge}}, SIGKDD Explorations, 9 (2007), pp.~75--79.

\bibitem{berahas2018derivative}
{\sc A.~S. Berahas, R.~H. Byrd, and J.~Nocedal}, {\em Derivative-free
  optimization of noisy functions via quasi-{N}ewton methods}, arXiv preprint
  arXiv:1803.10173,  (2018).

\bibitem{TakacNIPS16}
{\sc A.~S. Berahas, J.~Nocedal, and M.~Takac}, {\em A multi-batch {L-BFGS}
  method for machine learning}, in NIPS, vol.~29, 2016.

\bibitem{bertsekas1999nonlinear}
{\sc D.~P. Bertsekas}, {\em {Nonlinear Programming}}, Athena scientific
  Belmont, 1999.

\bibitem{Hit-And-Run}
{\sc D.~Bertsimas and S.~Vempala}, {\em Solving convex programs by random
  walks}, J. ACM, 51 (2004), pp.~540--556.

\bibitem{biegler2003large}
{\sc L.~T. Biegler, O.~Ghattas, M.~Heinkenschloss, and B.~van
  Bloemen~Waanders}, {\em Large-scale pde-constrained optimization: an
  introduction}, in Large-Scale PDE-Constrained Optimization, Springer, 2003,
  pp.~3--13.

\bibitem{biros2005parallel}
{\sc G.~Biros and O.~Ghattas}, {\em Parallel {Lagrange--Newton--Krylov--Schur}
  methods for {PDE}-constrained optimization. part i: The {Krylov--Schur}
  solver}, SIAM J. Sci. Comput., 27 (2005), pp.~687--713.

\bibitem{biros2005parallel2}
{\sc G.~Biros and O.~Ghattas}, {\em Parallel {Lagrange--Newton--Krylov--Schur}
  methods for {PDE}-constrained optimization. part ii: The {Lagrange--Newton}
  solver and its application to optimal control of steady viscous flows}, SIAM
  J. Sci. Comput., 27 (2005), pp.~714--739.

\bibitem{WRCR:WRCR23173}
{\sc E.~K. Bjarkason, O.~J. Maclaren, J.~P. O'Sullivan, and M.~J. O'Sullivan},
  {\em Randomized truncated {SVD} {L}evenberg-{M}arquardt approach to
  geothermal natural state and history matching}, Water Resour. Res., 54
  (2018), pp.~2376--2404.

\bibitem{NocedalBottou}
{\sc L.~Bottou, F.~E. Curtis, and J.~Nocedal}, {\em Optimization methods for
  large-scale machine learning}, SIAM Review, 60 (2018), pp.~223--311.

\bibitem{BoydVandenbergheBook}
{\sc S.~Boyd and L.~Vandenberghe}, {\em Convex Optimization}, Cambridge
  University Press, 2004.

\bibitem{bui2013computational}
{\sc T.~Bui-Thanh, O.~Ghattas, J.~Martin, and G.~Stadler}, {\em A computational
  framework for infinite-dimensional {B}ayesian inverse problems part i: The
  linearized case, with application to global seismic inversion}, SIAM J. Sci.
  Comput., 35 (2013), pp.~A2494--A2523.

\bibitem{ByrdPaper}
{\sc R.~H. Byrd, S.~L. Hansen, J.~Nocedal, and Y.~Singer}, {\em A stochastic
  quasi-{N}ewton method for large-scale optimization}, SIAM J. Optim., 26
  (2016), pp.~1008--1031.

\bibitem{ceruti2014danco}
{\sc C.~Ceruti, S.~Bassis, A.~Rozza, G.~Lombardi, E.~Casiraghi, and
  P.~Campadelli}, {\em Danco: An intrinsic dimensionality estimator exploiting
  angle and norm concentration}, Pattern recognition, 47 (2014),
  pp.~2569--2581.

\bibitem{ChenWildSTARS2015}
{\sc R.~Chen and S.~Wild}, {\em Randomized derivative-free optimization of
  noisy convex functions}, arXiv preprint arXiv:1507.03332,  (2015).

\bibitem{choromanski2018structured}
{\sc K.~Choromanski, M.~Rowland, V.~Sindhwani, R.~E. Turner, and A.~Weller},
  {\em Structured evolution with compact architectures for scalable policy
  optimization}, in ICML, 2018.

\bibitem{connDFO}
{\sc A.~R. Conn, K.~Scheinberg, and L.~N. Vicente}, {\em Introduction to
  Derivative-Free Optimization}, vol.~8, SIAM, 2009.

\bibitem{LikelihoodInformedDimensionReduction_Alessio}
{\sc T.~Cui, J.~Martin, Y.~M. Marzouk, A.~Solonen, and A.~Spantini}, {\em
  Likelihood-informed dimension reduction for nonlinear inverse problems},
  Inverse Problems, 30 (2014), p.~114015.

\bibitem{RandomCuttingPlane}
{\sc F.~Dabbene, P.~S. Shcherbakov, and B.~T. Polyak}, {\em A randomized
  cutting plane method with probabilistic geometric convergence}, SIAM J.
  Optim., 20 (2010), pp.~3185--3207.

\bibitem{SAGA}
{\sc A.~Defazio, F.~Bach, and S.~Lacoste-Julien}, {\em {SAGA}: A fast
  incremental gradient method with support for non-strongly convex composite
  objectives}, NIPS, 27 (2014), pp.~1646--1654.

\bibitem{FINITO}
{\sc A.~Defazio, J.~Domke, et~al.}, {\em Finito: A faster, permutable
  incremental gradient method for big data problems}, in ICML, 2014,
  pp.~1125--1133.

\bibitem{dolan2002benchmarking}
{\sc E.~D. Dolan and J.~J. Mor{\'e}}, {\em Benchmarking optimization software
  with performance profiles}, Math. Program., 91 (2002), pp.~201--213.

\bibitem{duchi2015optimal}
{\sc J.~C. Duchi, M.~I. Jordan, M.~J. Wainwright, and A.~Wibisono}, {\em
  Optimal rates for zero-order convex optimization: The power of two function
  evaluations}, IEEE T. Inform. Theory, 61 (2015), pp.~2788--2806.

\bibitem{dvurechensky2018accelerated}
{\sc P.~Dvurechensky, A.~Gasnikov, and E.~Gorbunov}, {\em An accelerated
  directional derivative method for smooth stochastic convex optimization},
  arXiv preprint arXiv:1804.02394,  (2018).

\bibitem{dvurechensky2017randomized}
{\sc P.~Dvurechensky, A.~Gasnikov, and A.~Tiurin}, {\em Randomized similar
  triangles method: A unifying framework for accelerated randomized
  optimization methods (coordinate descent, directional search, derivative-free
  method)}, arXiv preprint arXiv:1707.08486,  (2017).

\bibitem{ErmolievWets1988}
{\sc Y.~Ermoliev and R.-B. Wets}, {\em Numerical Techniques for Stochastic
  Optimization}, Springer-Verlag, 1988.

\bibitem{fefferman2016testing}
{\sc C.~Fefferman, S.~Mitter, and H.~Narayanan}, {\em Testing the manifold
  hypothesis}, J. Am. Math. Soc., 29 (2016), pp.~983--1049.

\bibitem{ferguson2017course}
{\sc T.~S. Ferguson}, {\em A Course in Large Sample Theory}, Routledge, 2017.

\bibitem{Flath2011}
{\sc H.~Flath, L.~Wilcox, V.~Ak\c{c}elik, J.~Hill, B.~{Van Bloemen Waanders},
  and O.~Ghattas}, {\em Fast algorithms for {B}ayesian uncertainty
  quantification in large-scale linear inverse problems based on low-rank
  partial hessian approximations}, SIAM J. Sci. Comput., 33 (2011),
  pp.~407--432.

\bibitem{GavianoRandomSearch1975}
{\sc M.~Gaviano}, {\em Some general results on convergence of random search
  algorithms in minimization problems}, in Towards Global Optimisation, 1975,
  pp.~149--157.

\bibitem{ge2015escaping}
{\sc R.~Ge, F.~Huang, C.~Jin, and Y.~Yuan}, {\em Escaping from saddle points:
  online stochastic gradient for tensor decomposition}, in Conference on
  Learning Theory, 2015, pp.~797--842.

\bibitem{GhadimiLan13}
{\sc S.~Ghadimi and G.~Lan}, {\em Stochastic first- and zeroth-order methods
  for nonconvex stochastic programming}, SIAM J. Optim.,  (2013),
  pp.~2341--2368.

\bibitem{GowerBFGS}
{\sc R.~M. Gower, D.~Goldfarb, and P.~Richt\'{a}rik}, {\em Stochastic block
  {BFGS}: Squeezing more curvature out of data}, in ICML, vol.~48, 2016,
  pp.~1869--1878.

\bibitem{griewank2008evaluating}
{\sc A.~Griewank and A.~Walther}, {\em Evaluating derivatives: principles and
  techniques of algorithmic differentiation}, vol.~105, SIAM, 2~ed., 2008.

\bibitem{gunzburger2003perspectives}
{\sc M.~D. Gunzburger}, {\em Perspectives in flow control and optimization},
  vol.~5, SIAM, 2003.

\bibitem{haber2012effective}
{\sc E.~Haber, M.~Chung, and F.~Herrmann}, {\em An effective method for
  parameter estimation with pde constraints with multiple right-hand sides},
  SIAM J. Optim., 22 (2012), pp.~739--757.

\bibitem{haber2012numerical}
{\sc E.~Haber, Z.~Magnant, C.~Lucero, and L.~Tenorio}, {\em Numerical methods
  for a-optimal designs with a sparsity constraint for ill-posed inverse
  problems}, Comput. Optim. Appl., 52 (2012), pp.~293--314.

\bibitem{cg851}
{\sc W.~W. Hager and H.~Zhang}, {\em Algorithm 851: {CG\_DESCENT}, a conjugate
  gradient method with guaranteed descent}, ACM Trans. Math. Software, 32
  (2006), pp.~113--137.

\bibitem{Karimi2016LinearCO}
{\sc {Hamed Karimi and Julie Nutini and Mark W. Schmidt}}, {\em Linear
  convergence of gradient and proximal-gradient methods under the
  {P}olyak-{L}ojasiewicz condition}, in ECML/PKDD, 2016.

\bibitem{CMA-ES2001}
{\sc N.~Hansen and A.~Ostermeier}, {\em Completely derandomized self-adaptation
  in evolution strategies}, Evol. Comput., 9 (2001), pp.~159--195.

\bibitem{SEGA}
{\sc F.~Hanzely, K.~Mishchenko, and P.~Richtarik}, {\em {SEGA}: Variance
  reduction via gradient sketching}, in NIPS, 2018, pp.~2086--2097.

\bibitem{hanzely2018accelerated}
{\sc F.~Hanzely and P.~Richt{\'a}rik}, {\em Accelerated coordinate descent with
  arbitrary sampling and best rates for minibatches}, arXiv preprint
  arXiv:1809.09354,  (2018).

\bibitem{horesh2010optimal}
{\sc L.~Horesh, E.~Haber, and L.~Tenorio}, {\em Optimal experimental design for
  the large-scale nonlinear ill-posed problem of impedance imaging},
  Large-Scale Inverse Problems and Quantification of Uncertainty,  (2010),
  pp.~273--290.

\bibitem{isaac2015scalable}
{\sc T.~Isaac, N.~Petra, G.~Stadler, and O.~Ghattas}, {\em Scalable and
  efficient algorithms for the propagation of uncertainty from data through
  inference to prediction for large-scale problems, with application to flow of
  the antarctic ice sheet}, J. Comput. Phys., 296 (2015), pp.~348--368.

\bibitem{jin2017accelerated}
{\sc C.~Jin, P.~Netrapalli, and M.~I. Jordan}, {\em Accelerated gradient
  descent escapes saddle points faster than gradient descent}, arXiv preprint
  arXiv:1711.10456,  (2017).

\bibitem{SVRG}
{\sc R.~Johnson and T.~Zhang}, {\em Accelerating stochastic gradient descent
  using predictive variance reduction}, NIPS, 26 (2013), pp.~315--323.

\bibitem{SimulatedAnnealing}
{\sc S.~Kirkpatrick, C.~D. Gelatt, and M.~P. Vecchi}, {\em Optimization by
  simulated annealing}, Science, 220 (1983), pp.~671--680.

\bibitem{LewisDirectionalSearch}
{\sc D.~Leventhal and A.~Lewis}, {\em Randomized {H}essian estimation and
  directional search}, Optimization, 60 (2011), pp.~329--345.

\bibitem{levina2005maximum}
{\sc E.~Levina and P.~J. Bickel}, {\em Maximum likelihood estimation of
  intrinsic dimension}, in NIPS, vol.~18, 2005, pp.~777--784.

\bibitem{lin2015universal}
{\sc H.~Lin, J.~Mairal, and Z.~Harchaoui}, {\em A universal catalyst for
  first-order optimization}, in NIPS, vol.~28, 2015, pp.~3384--3392.

\bibitem{liu2018zeroth}
{\sc S.~Liu, B.~Kailkhura, P.-Y. Chen, P.~Ting, S.~Chang, and L.~Amini}, {\em
  Zeroth-order stochastic variance reduction for nonconvex optimization}, in
  NIPS, vol.~31, 2018, pp.~3731--3741.

\bibitem{logg2012automated}
{\sc A.~Logg, K.-A. Mardal, and G.~Wells}, {\em Automated solution of
  differential equations by the finite element method: The FEniCS book},
  vol.~84, Springer Science \& Business Media, 2012.

\bibitem{maclaurin2015autograd}
{\sc D.~Maclaurin, D.~Duvenaud, and R.~P. Adams}, {\em Autograd: Effortless
  gradients in numpy}, in ICML 2015 AutoML Workshop, 2015.

\bibitem{MISO}
{\sc J.~Mairal}, {\em Optimization with first-order surrogate functions}, in
  ICML, 2013, pp.~783--791.

\bibitem{mezzadri2006generate}
{\sc F.~Mezzadri}, {\em How to generate random matrices from the classical
  compact groups}, in Notices of the American Mathematical Society, vol.~54,
  2006.

\bibitem{Mokhtari}
{\sc A.~Mokhtari and A.~Ribeiro}, {\em {RES}: Regularized stochastic {BFGS}
  algorithm}, {IEEE T. Signal Proces.}, 62 (2014), pp.~6089--6104.

\bibitem{Mokhtari2015GlobalCO}
{\sc A.~Mokhtari and A.~Ribeiro}, {\em Global convergence of online limited
  memory {BFGS}}, JMLR, 16 (2015), pp.~3151--3181.

\bibitem{StochasticLBFGS_Moritz}
{\sc P.~Moritz, R.~Nishihara, and M.~I. Jordan}, {\em A linearly-convergent
  stochastic {L-BFGS} algorithm}, in AISTATS, 2016.

\bibitem{Nesterov1983}
{\sc Y.~Nesterov}, {\em {A Method of Solving a Convex Programming Problem with
  Convergence Rate $\mathcal{O}(1/k^2)$.}}, {Soviet Mathematics Doklady}, 27
  (1983), pp.~372--376.

\bibitem{nesterov2011random}
{\sc Y.~Nesterov}, {\em Random gradient-free minimization of convex functions},
  tech. report, Universit{\'e} catholique de Louvain, Center for Operations
  Research and Econometrics (CORE), 2011.

\bibitem{nesterov2012efficiency}
{\sc Y.~Nesterov}, {\em Efficiency of coordinate descent methods on huge-scale
  optimization problems}, SIAM J. Optim., 22 (2012), pp.~341--362.

\bibitem{nesterov2013introductory}
{\sc Y.~Nesterov}, {\em Introductory Lectures on Convex Optimization: A Basic
  Course}, vol.~87, Springer Science \& Business Media, 2013.

\bibitem{nesterov2017random}
{\sc Y.~Nesterov and V.~Spokoiny}, {\em Random gradient-free minimization of
  convex functions}, Foundations of Computational Mathematics, 17 (2017),
  pp.~527--566.
\newblock First appeard as CORE discussion paper 2011.

\bibitem{NASA_NielsenDiskin}
{\sc E.~J. Nielsen and B.~Diskin}, {\em Discrete adjoint-based design for
  unsteady turbulent flows on dynamic overset unstructured grids}, AIAA
  Journal, 51 (2013), pp.~1355--1373.

\bibitem{NocedalWright}
{\sc J.~Nocedal and S.~Wright}, {\em {Numerical Optimization}},
  Springer-Verlag, 2~ed., 1999.

\bibitem{peng2016arock}
{\sc Z.~Peng, Y.~Xu, M.~Yan, and W.~Yin}, {\em Arock: an algorithmic framework
  for asynchronous parallel coordinate updates}, SIAM J. Sci. Comput., 38
  (2016), pp.~A2851--A2879.

\bibitem{petra2014computational}
{\sc N.~Petra, J.~Martin, G.~Stadler, and O.~Ghattas}, {\em A computational
  framework for infinite-dimensional {B}ayesian inverse problems, part ii:
  Stochastic {N}ewton mcmc with application to ice sheet flow inverse
  problems}, SIAM J. Sci. Comput., 36 (2014), pp.~A1525--A1555.

\bibitem{pilanci2015randomized}
{\sc M.~Pilanci and M.~J. Wainwright}, {\em Randomized sketches of convex
  programs with sharp guarantees}, IEEE T. Inform. Theory, 61 (2015),
  pp.~5096--5115.

\bibitem{polyak1964some}
{\sc B.~T. Polyak}, {\em Some methods of speeding up the convergence of
  iteration methods}, USSR Comp Math. Math. Phys.+, 4 (1964), pp.~1--17.

\bibitem{powell1973search}
{\sc M.~J. Powell}, {\em On search directions for minimization algorithms},
  Math. Program., 4 (1973), pp.~193--201.

\bibitem{richtarik2014iteration}
{\sc P.~Richt{\'a}rik and M.~Tak{\'a}c}, {\em Iteration complexity of
  randomized block-coordinate descent methods for minimizing a composite
  function}, Math. Program., 144 (2014), pp.~1--38.

\bibitem{SAG}
{\sc N.~L. Roux, M.~Schmidt, and F.~R. Bach}, {\em A stochastic gradient method
  with an exponential convergence rate for finite training sets}, in NIPS,
  vol.~25, 2012, pp.~2663--2671.

\bibitem{salimans2017evolution}
{\sc T.~Salimans, J.~Ho, X.~Chen, S.~Sidor, and I.~Sutskever}, {\em Evolution
  strategies as a scalable alternative to reinforcement learning}, arXiv
  preprint arXiv:1703.03864,  (2017).

\bibitem{Schraudolph2007ASQ}
{\sc N.~N. Schraudolph, J.~Yu, and S.~G{\"u}nter}, {\em A stochastic
  quasi-{N}ewton method for online convex optimization}, in AISTATS, 2007.

\bibitem{SDCA}
{\sc S.~Shalev-Shwartz and T.~Zhang}, {\em Stochastic dual coordinate ascent
  methods for regularized loss minimization}, JMLR, 14 (2013), pp.~567--599.

\bibitem{snelson2006sparse}
{\sc E.~Snelson and Z.~Ghahramani}, {\em Sparse {G}aussian processes using
  pseudo-inputs}, in NIPS, 2006, pp.~1257--1264.

\bibitem{SolisWets81}
{\sc F.~Solis and R.~J.-B. Wets}, {\em Minimization by random search
  techniques}, Math. Oper. Res., 6 (1981), pp.~19--30.

\bibitem{spantini2017goal}
{\sc A.~Spantini, T.~Cui, K.~Willcox, L.~Tenorio, and Y.~Marzouk}, {\em
  Goal-oriented optimal approximations of {B}ayesian linear inverse problems},
  SIAM J. Sci. Comput., 39 (2017), pp.~S167--S196.

\bibitem{spantini2015optimal}
{\sc A.~Spantini, A.~Solonen, T.~Cui, J.~Martin, L.~Tenorio, and Y.~Marzouk},
  {\em Optimal low-rank approximations of {B}ayesian linear inverse problems},
  SIAM J. Sci. Comput., 37 (2015), pp.~A2451--A2487.

\bibitem{stich2013optimization}
{\sc S.~U. Stich, C.~Muller, and B.~Gartner}, {\em Optimization of convex
  functions with random pursuit}, SIAM J. Optim., 23 (2013), pp.~1284--1309.

\bibitem{Tan2016BarzilaiBorweinSS}
{\sc C.~Tan, N.~S. Aybat, Y.-H. Dai, and Y.~Qian}, {\em {B}arzilai-{B}orwein
  step size for stochastic gradient descent}, in NIPS, vol.~29, 2016.

\bibitem{LuisBook}
{\sc L.~Tenorio}, {\em An Introduction to Data Analysis and Uncertainty
  Quantification for Inverse Problems}, Mathematics in Industry, SIAM, 2017.

\bibitem{titsias2009variational}
{\sc M.~Titsias}, {\em Variational learning of inducing variables in sparse
  {G}aussian processes}, in AISTATS, 2009, pp.~567--574.

\bibitem{tropp2017practical}
{\sc J.~A. Tropp, A.~Yurtsever, M.~Udell, and V.~Cevher}, {\em Practical
  sketching algorithms for low-rank matrix approximation}, SIAM J. Matrix Anal.
  Appl., 38 (2017), pp.~1454--1485.

\bibitem{udell2019big}
{\sc M.~Udell and A.~Townsend}, {\em Why are big data matrices approximately
  low rank?}, SIAM J. Math. Data Sci., 1 (2019), pp.~144--160.

\bibitem{NIPS2002_2203}
{\sc P.~Vincent and Y.~Bengio}, {\em Manifold {P}arzen windows}, in NIPS,
  S.~Becker, S.~Thrun, and K.~Obermayer, eds., vol.~16, MIT Press, 2003,
  pp.~849--856.

\bibitem{wang2009minimal}
{\sc Q.~Wang, P.~Moin, and G.~Iaccarino}, {\em Minimal repetition dynamic
  checkpointing algorithm for unsteady adjoint calculation}, SIAM J. Sci.
  Comput., 31 (2009), pp.~2549--2567.

\bibitem{MaGoldfarbStochasticQuasiNewton}
{\sc X.~Wang, S.Ma, D.~Goldfarb, and W.~Liu}, {\em Stochastic quasi-{N}ewton
  methods for nonconvex stochastic optimization}, SIAM J. Optim., 27 (2017),
  pp.~927--956.

\bibitem{warga1963minimizing}
{\sc J.~Warga}, {\em Minimizing certain convex functions}, Journal of the
  Society for Industrial and Applied Mathematics, 11 (1963), pp.~588--593.

\bibitem{Rasmussen2009GaussianPF}
{\sc C.~K. Williams and C.~E. Rasmussen}, {\em Gaussian Processes for Machine
  Learning}, vol.~2, MIT Press Cambridge, MA, 2006.

\bibitem{williams2001using}
{\sc C.~K. Williams and M.~Seeger}, {\em Using the {N}ystr{\"o}m method to
  speed up kernel machines}, in NIPS, vol.~14, 2001, pp.~682--688.

\bibitem{wright2015coordinate}
{\sc S.~J. Wright}, {\em Coordinate descent algorithms}, Math. Program., 151
  (2015), pp.~3--34.

\end{thebibliography}

\newpage
\section*{Stochastic Subspace Descent Supplementary Material}
\subsection*{S1. Known results}
We summarize results from optimization and probability used in the proofs of convergence.
\subsubsection*{Optimization}
 \begin{lemma}[A property of Lipschitz-continuous first derivatives]\label{Supp-Lipschitz}
Let $f: \reals^d \to \reals$ be continuously-differentiable with $\lambda$-Lipschitz 	first derivative. Then for any $\x, \y \in \reals^d$,
\begin{equation*}
\abs{f(\y) - f(\x) - \langle \nabla f(\x), \y-\x\rangle} \leq \frac{\lambda}{2} \norm{\y-\x}^2.
\end{equation*}
\end{lemma}
The proof can be found in \cite[pg. 25-26]{nesterov2013introductory}.
 \begin{lemma}[Strong-convexity implies PL-inequality]\label{PLLemma}
 	Let $f :\reals^d \to \reals$ be a $\gamma$-strongly convex, continuously-differentiable function. Denote the unique minimizer of $f$ by  $\x_*$. Then:
 	\begin{equation}\label{Supp-PLineq}
 	f(\x) - f(\x_*) \leq \frac{1}{2\gamma} \norm{\nabla f(\x)}^2.
 	\end{equation}
 \end{lemma}
 Inequality \eqref{Supp-PLineq} is the PL-inequality. The proof is simple and can be found in the Appendix of \cite{Karimi2016LinearCO}.
 
\begin{lemma}[A characterization of convex functions]\label{Supp-lemma: convex bound}
A function $f:\reals^d \to \reals$ is convex with $\lambda$-Lipschitz derivative if and only if for all $\x, \y \in \reals^d$
\begin{equation*}
    \langle \nabla f(\x) - \nabla f(\y), \x-\y \rangle  \geq \frac{1}{\lambda} \norm{\nabla f(\x) - \nabla f(\y)}^2 \geq 0.
\end{equation*}
\end{lemma} 
The proof can be found in \cite[pg. 56]{nesterov2013introductory}. If the function is strongly convex then the bound is sharpened according to the following theorem, the proof of which can be found in \cite[pg. 62]{nesterov2013introductory}.
\begin{theorem}[Strongly-convex functions]\label{Supp-thm: strongly-convex bound}
If a function $f:\reals^d \to \reals$ is $\gamma$-strongly convex with $\lambda$-Lipschitz derivative then for all $\x, \y \in \reals^d$
\begin{equation*}
    \langle \nabla f(\x) - \nabla f(\y), \x-\y \rangle  \geq \frac{\gamma \lambda}{\gamma + \lambda} \norm{\x-\y}^2 + \frac{1}{\gamma + \lambda} \norm{\nabla f(\x) - \nabla f(\y)}^2.
\end{equation*}
\end{theorem}
\subsubsection*{Probability}
 A sequence of random vectors $(\x_n)$ is said to converge a.s. to a random vector $\x$ if and only if  $\mathbb{P}(\x_n \to \x) = 1$. It is said to converge in probability if and only if for all $\epsilon > 0$, $\mathbf{P}(\norm{\x_n - \x} < \epsilon ) \to 0$ as $n \to \infty$. It converges in $L^r$ for some integer $r >0$ if and only if $\mathbb{E}\norm{\x_n - \x}^r \rightarrow 0$. We will write $\x_n \convas \x$, $\x_n \convp \x$, and $\x_n \overset{L^r}{\rightarrow}\x$, respectively. 

It is readily apparent that $\x_n \convas \x$ if and only if for all $\epsilon >0$
 	\begin{equation*}
 	\lim_{n\to\infty}\mathbb{P}\left(\bigcap_{k\geq n}\{\norm{\x_k-\x}<\epsilon\}\right) = 1
 	\end{equation*}
 	as $n \to \infty$ \cite{ferguson2017course}. It follows from this characterization  that $\x_n \convp \x$ if $\x_n \convas \x$. By subadditivity,
 	\begin{equation}\label{Supp-eq: fast convp}
 	    \sum_n \mathbb{P} (\norm{\x_n - \x} > \epsilon) < \infty \quad \forall ~\epsilon > 0 \implies \x_n \convas \x.
 	\end{equation}
   Finally, by Markov's inequality, if $\x_n \overset{L^r}{\to}\x$, then $\x_n \convp \x$.
 
\subsection*{S2. Proofs of main results}
We provide expanded proofs of the main results.
\subsubsection*{Proof of Theorem \ref{thm:convergence}} Define the filtration $\mathcal{F}_k = \sigma(\mathbf{P}_1, \ldots, \mathbf{P}_{k-1}), ~ k>1$, and $\mathcal{F}_1 = \{ \emptyset , \Omega\}$.
	Because $f$ is continuously-differentiable with a $\lambda$-Lipschitz derivative, Lemma \ref{Supp-Lipschitz} yields  
	\begin{equation}\label{Supp-Taylor}
	f(\x_{k+1}) \leq  f(\x_k) + \nabla f(\x_k)^{\T}(\x_{k+1}-\x_k) + \frac{\lambda}{2}\norm{\x_{k+1}-\x_k}^2
	\end{equation}
Let $f_e(\x) = f(\x)-f_*$ be the error for a particular $\x$.   Then, equations \eqref{eq: iterations} and \eqref{Supp-Taylor} yield:
	\begin{align}\label{Supp-delineq}
	\nonumber	f_e(\x_{k+1}) - f_e(\x_k)  &\leq -\alpha \langle \nabla f(\x_k),~ \mathbf{P}_k\mathbf{P}_k^{\T}\nabla f(\x_k) \rangle + \frac{\alpha^2\lambda}{2}  \langle  \mathbf{P}_k\mathbf{P}_k^{\T}\nabla f(\x_k),~ \mathbf{P}_k\mathbf{P}_k^{\T}\nabla f(\x_k) \rangle \\
	\nonumber	&= -\alpha \langle \nabla f(\x_k),~ \mathbf{P}_k\mathbf{P}_k^{\T}\nabla f(\x_k) \rangle + \frac{d\alpha^2\lambda}{2\ell}  \langle \nabla f(\x_k),~ \mathbf{P}_k\mathbf{P}_k^{\T}\nabla f(\x_k) \rangle\\
	\nonumber	&= \left(-\alpha +\frac{d\alpha^2\lambda}{2\ell}\right)  \langle \nabla f(\x_k),~ \mathbf{P}_k\mathbf{P}_k^{\T}\nabla f(\x_k) \rangle \\
		\nonumber	&= \alpha_{\lambda} \langle \nabla f(\x_k),~ \mathbf{P}_k\mathbf{P}_k^{\T}\nabla f(\x_k) \rangle \quad \text{with} \quad \alpha_{\lambda} = \left(-\alpha + \frac{d\alpha^2\lambda}{2\ell} \right),
	\end{align}
	where we have used the fact that $\mathbf{P}_k\mathbf{P}_k^{\T}\mathbf{P}_k\mathbf{P}_k^{\T}= (d/\ell)\mathbf{P}_k\mathbf{P}_k^{\T}$. We can choose $\alpha$ such that $\alpha_{\lambda} <0$; that is, $0< \alpha < 2\ell/d\lambda$. With this choice the right-hand side is  non-positive and the errors are non-increasing. Since the error is bounded below by zero the sequence converges almost surely but the actual limit remains to be determined. Furthermore, since the sequence is bounded above by $f_e(\x_0)$, Dominated convergence implies convergence of the sequence in $L^1$. To find the limit we take conditional expectations of both sides to get
	\begin{equation*}
	\Expectation \left[f_e(\x_{k+1}) \mid \mathcal{F}_k \right] \leq \alpha_{\lambda} \Expectation \left[ \langle \nabla f(\x_k), \mathbf{P}_k\mathbf{P}_k^{\T} \nabla f(\x_k) \rangle  \bigg| \mathcal{F}_k \right] +f_e(\x_k).
	\end{equation*}
	Since $\x_k$ is $\mathcal{F}_k$-measurable, the gradients can be pulled out of the conditional expectation, the expectation can be passed through due to linearity. Noting that $\Expectation \left[\mathbf{P}_k\mathbf{P}_k^{\T} \mid \mathcal{F}_k \right]=\I_d$,  this leaves us with
	\begin{equation} \label{Supp-intermediate}
	\Expectation \left[f_e(\x_{k+1}) \mid \mathcal{F}_k\right] \leq  \alpha_{\lambda} \norm{\nabla f(\x_k)}^2 + f_e(\x_{k}) .
	\end{equation}
  The first term on the right-hand side in equation \eqref{Supp-intermediate} is negative so the PL-inequality leads to:	
	\begin{equation*}
	\Expectation \left[f_e(\x_{k+1}) \mid \mathcal{F}_k\right] \leq 2\gamma  \alpha_{\lambda} f_e(\x_k) + f_e(\x_{k})= \left(1+2\gamma\alpha_{\lambda}\right)f_e(\x_k).
	\end{equation*}
	Recursive application yields
	\begin{equation*}
	\Expectation \left[f(\x_{k+1})-f_* \mid \mathcal{F}_k \right]\leq  \left(1-2\gamma\alpha_{\lambda}\right)^{k+1} \left(f(\x_0) - f_*\right).
	\end{equation*}
	Taking an expectation of both sides yields
	\begin{equation*}
	    	\Expectation f(\x_{k+1})-f_* \leq  \left(1+2\gamma\alpha_{\lambda}\right)^{k+1} \left(f(\x_0) - f_*\right).
	\end{equation*}
	Any step-size $0 < \alpha < 2\ell/d\lambda$ forces the right-hand side to converge to zero. Thus, since $	f_e(\x_k) \convas \X$ for some $\X \in L^1$ and $f_e(\x_k)\overset{L^1}{\to} 0$, we have both $f(\x_k) \convas f_*$ and $f(\x_k) \overset{L^1}{\to} f_*$. 
\qed
\subsubsection*{Proof of Corollary \ref{corr:strong-convexity}  part (i)}
By Lemma \ref{PLLemma}, the strong-convexity assumed by this corollary implies the assumptions of Theorem \ref{thm:convergence}, so that $f(\x_k) \convas f(\x_*)$. Now by strong-convexity, 
	\begin{equation*}
	f(\x_k) - f(\x_*) \geq \frac{\gamma}{2} \norm{\x_*-\x_k}.
	\end{equation*}
	Since the left-hand side converges almost surely to zero and $\gamma > 0$, we have $\x_k \convas \x_*$.
\subsubsection*{Proof of Corollary \ref{corr:strong-convexity}  part (ii)}

	We begin by providing an upper bound for $\norm{\nabla f(\x_k)}^2$. Rearranging the terms in equation \eqref{Supp-intermediate} we have
	\begin{equation} \label{Supp-UpperBound}
	\alpha_{\lambda} ^{-1}\Expectation \left[f_e(\x_k) - f_e(\x_{k+1})\mid \mathcal{F}_k \right] \geq \norm{\nabla f(\x_k)}^2.
	\end{equation}
Combing \eqref{Supp-UpperBound} with \eqref{Supp-PLineq} we get 
	\begin{equation*}
	2\gamma f_e(\x_k) \leq \norm{\nabla f(\x_k)}^2 \leq \alpha_{\lambda}^{-1} \Expectation\left[f_e(\x_{k}) - f_e(\x_{k+1}) \mid \mathcal{F}_k\right].
	\end{equation*}
	That is,
	\begin{equation}\label{Supp-nearly}
	\Expectation \left[f_e(\x_{k+1}) \mid \mathcal{F}_k\right] \leq \left(1+2\gamma\alpha_{\lambda}\right)f_e(\x_k).
	\end{equation}
	Repeated recursion then yields	$\Expectation \left[f_e(\x_{k+1}) \mid \mathcal{F}_k\right] \leq \left(1+2\gamma\alpha_{\lambda}\right)^{k+1}f_e(\x_0).$ Recalling that $\alpha_{\lambda} = -\alpha + d\alpha^2\lambda /(2\ell)$, the minimum of the right-hand side is achieved when $\alpha=\ell/(d\lambda)$, resulting in
	\begin{equation*}
		\Expectation \left[f_e(\x_{k+1}) \mid \mathcal{F}_k\right] \leq \left(1-\ell\gamma / (d \lambda) \right)^{k+1}f_e(\x_0).
	\end{equation*}
	Taking an expectation of both sides proves the corollary.
\qed

	\subsubsection*{Proof of Theorem \ref{thm: convergence-convex}}
	We follow the proof of Theorem \ref{thm:convergence} until equation \eqref{Supp-intermediate}, then we rearrange terms to obtain,
	\begin{align}\label{Supp-eq:simplified-Lipschitz}
	\Expectation \left[f(\x_{k+1}) \mid \mathcal{F}_k \right] \leq f(\x_k) + \alpha_{\lambda} \norm{\nabla f(\x_k)}^2.
	\end{align}
	By convexity, for any $\x_* \in \mathcal{D}$,
	\begin{equation*}
	f(\x_k)- f_*  \leq \nabla f(\x_k)^{\T} \left(\x_k - \x_*\right) \leq \norm{\nabla f(\x_k)} \norm{\x_k - \x_*}\leq R \norm{\nabla f(\x_k)}.
	\end{equation*}
 That is,
	\begin{equation}\label{Supp-eq: grad lower bound}
	\norm{\nabla f(\x_k)} \geq f_e(\x_k)/R.   
	\end{equation}
	Plugging \eqref{Supp-eq: grad lower bound} into equation \eqref{Supp-eq:simplified-Lipschitz} and letting $\alpha = \ell/(d\lambda)$ results in 
	\begin{equation}\label{Supp-eq: required}
 \Expectation [f_e(\x_{k+1}) \mid \mathcal{F}_k]- f_e(\x_k)  \leq  -\frac{\alpha  f_e(\x_k)^2}{2 R^2},
	\end{equation}
	and one more expectation yields
	\begin{equation*}
	    \Expectation[f_e(\x_{k+1}) - f_e(\x_k)] \leq -\frac{\alpha \Expectation f_e(\x_k)^2}{2R^2} 
	    \leq -\frac{\alpha \left(\Expectation f_e(\x_k)\right)^2}{2R^2}.
	\end{equation*}
We note that
	\begin{align*}
 -\frac{\alpha}{2 R^2} &\geq \frac{\Expectation f_e(\x_{k+1})- \Expectation f_e(\x_k) }{\left(\Expectation f_e(\x_k)\right)^2} \\ 
 &\geq \frac{\Expectation  f_e(\x_{k+1}) - \Expectation f_e(\x_k)}{\Expectation f_e(\x_k)\Expectation f_e(\x_{k+1})}\\
&=  \frac{1}{\Expectation f_e(\x_k)} - \frac{1}{\Expectation f_e(\x_{k+1})},
	\end{align*}
	where the second inequality is a result of $\Expectation f_e(\x_{k+1}) < \Expectation f_e(\x_k)$, a consequence of \eqref{Supp-eq:simplified-Lipschitz}. A rearrangement leads to
	\begin{equation}\label{Supp-eq: before recursion}
	\frac{1}{\Expectation f_e(\x_{k+1})} \geq 	\frac{1}{\Expectation f_e(\x_{k})} + \frac{\alpha}{2 R^2}.
	\end{equation}
	Applying \eqref{Supp-eq: before recursion} recursively,
	\begin{equation*}
	\frac{1}{\Expectation f_e(\x_{k+1})} \geq 	\frac{1}{\Expectation f_e(\x_{0})} + \frac{\alpha k}{2 R^2} \geq \frac{\alpha k}{2 R^2}.
	\end{equation*}
	Replace $\alpha$ with $\ell/(d\lambda)$ and cross-multiply to get the desired result
	\begin{equation*}
	\Expectation  f_e(\x_{k+1}) \leq \frac{2d\lambda R^2}{k \ell }. 
	\end{equation*}
\qed

\subsubsection*{Proof of Lemma \ref{thm:SVRG-AS}}
	We begin with the convergence of the function evaluations. $f(\x_k) \overset{L^1}{\longrightarrow} f(\x^*)$ implies $f(\x_k) \convp f(\x_*)$ at a rate that is at least equivalent. In particular,
	\begin{equation*}
	\mathbb{P}(\norm{f(\x_k) - f(\x_*)} > \epsilon) \leq \beta^k \left[f(\tilde{\x}_0) - f(\x_*)\right], \qquad \beta<1.
	\end{equation*}
	Sum both sides to get
	\begin{equation*}
	\sum_{k=1}^{\infty}	\mathbb{P}(\norm{f(\x_k) - f(\x_*)} > \epsilon) \leq \sum_{k=1}^{\infty} \beta^k \left[f(\x_0) - f(\x_*)\right] \leq \frac{\left[f(\x_0) - f(\x_*)\right]}{1- \beta} < \infty.
	\end{equation*}
	Where the second inequality is the sum of the geometric series since $\beta< 1$. Therefore the result follows from \eqref{Supp-eq: fast convp}. For the convergence of the iterates, note that if $f$ is strongly convex it has a unique minimizer.
	\qed

\subsubsection*{Proof of Theorem \ref{thm:VarianceReducedRandomGradient} part (i)}
We define the following filtration that encodes information at each step $s,k$:
\begin{equation*}
\mathcal{F}_{s,k} = \sigma(\textbf{P}_{1,1}, \ldots, \textbf{P}_{1,m}, \ldots, \textbf{P}_{s-1,1},\ldots,  \textbf{P}_{s-1,m}, \textbf{P}_{s,1},\ldots,  \textbf{P}_{s,k-1}  , J_1, \ldots J_{s-1}),
\end{equation*}
which makes $\x_{s,k}$ and $\tilde{\x}_{s-1}$ $\mathcal{F}_{s,k}$-measurable. Note that 
\begin{equation}\label{Supp-doubleexpect}
\Expectation f(\tilde{\x}_s) = \Expectation~ \Expectation (f(\tilde{\x}_s)\mid \mathcal{F}_{s,m}) = \Expectation \frac{1}{m}\sum_{k=0}^{m-1} f(\x_{s,k}).
\end{equation}
The recursion for epoch $s$ can be written as 
	\begin{equation}\label{Supp-eqn: t4-recursion}
	\x_{s,k} = \x_{s,k-1} - \alpha~ \v,
	\end{equation}
	where 
	\begin{equation*}
	\v  =  \P\P^{\T}\nabla f(\x_{s,k-1}) -  \left(\P\P^{\T}-\I\right)\nabla f(\tilde{\x}_{s-1}).
	\end{equation*}
   Equation \eqref{Supp-eqn: t4-recursion} then leads to
	\begin{equation}\label{Supp-eqn: t4-normsquared}
	\norm{\x_{s,k}-\x_*}^2 = \norm{\x_{s,k-1}-\x_*}^2 - 2\alpha (\x_{s,k-1}-\x_*)^{\T}\v + \alpha^2 \norm{\v}^2,
	\end{equation}
which in turn leads to
	\begin{equation}\label{Supp-eqn: t4-expected normsquared}
\CondExp{\norm{\x_{s,k}-\x_*}^2} =\norm{\x_{s,k-1}-\x_*}^2 - 2\alpha (\x_{s,k-1}-\x_*)^{\T}\CondExp{\v} + \alpha^2 \CondExp{\norm{\v}^2}.
	\end{equation}
By definition, $\CondExp{\v} = \nabla f(\x_{s,k-1})$. For the last term we have:
\begin{align}\label{Supp-eqn: t4-v-squared}
    \nonumber \CondExp{\norm{\v}^2} ={}& \CondExp{\rho \nabla f(\x_{s,k-1})^{\T} \P\P^{\T} \nabla f(\x_{s,k-1})} \\
 \nonumber   &+\CondExp{ (\rho-2) \nabla f(\tilde{\x}_{s-1})^{\T} \left(\P\P^{\T} + \frac{1}{\rho-2}\I \right) \nabla f(\tilde{\x}_{s-1})} \\
\nonumber    &- \CondExp{2 (\rho-1) \nabla f(\x_{s,k-1})^{\T} \nabla f(\tilde{\x}_{s-1})} \\
    ={}& \rho \norm{\nabla f(\x_{s,k-1})}^2 + (\rho-1) \norm{f(\tilde{\x}_{s-1})}^2 - 2(\rho-1) \nabla f(\x_{s,k-1})^{\T}f(\tilde{\x}_{s-1}).
\end{align}
We plug equation \eqref{Supp-eqn: t4-v-squared} into equation \eqref{Supp-eqn: t4-expected normsquared} to get
\begin{align}\label{Supp-eqn: t4-pluggedin}
\nonumber \CondExp{\norm{\x_{s,k}-\x_*}^2} =&\norm{\x_{s,k-1}-\x_*}^2 - 2\alpha (\x_{s,k-1}-\x_*)^{\T}\nabla f(\x_{s,k-1}) \\
\nonumber &+ \alpha^2 \left( \rho \norm{\nabla f(\x_{s,k-1})}^2 + (\rho-1) \norm{f(\tilde{\x}_{s-1})}^2 \right. \\& \left. - 2(\rho-1)  \nabla f(\x_{s,k-1})^\top f(\tilde{\x}_{s-1}) \vphantom{\norm{\nabla f(\x_{s,k-1})}^2}\right).
\end{align}
Taking an expectation of both sides of the equality in equation \eqref{Supp-eqn: t4-pluggedin} and summing from $k=1$ to $k=m$ yields
\begin{align}\label{Supp-eqn: t4-intermediate-step-after-sum}
  \nonumber  \sum_{k=1}^m \Expectation \left[\norm{\x_{s,k} - \x_*}^2 \right] ={}&  \sum_{k=1}^m \Expectation \left[\norm{\x_{s,k-1} - \x_*}^2 \right] -2\alpha  \sum_{k=1}^m \Expectation \left[(\x_{s,k-1}-\x_*)^{\T}\nabla f(\x_{s,k-1}) \right] \\
\nonumber    &+  \alpha^2 \rho \sum_{k=1}^m \Expectation\left[\norm{\nabla f(\x_{s,k-1})}^2\right] + \alpha^2 m(\rho-1) \norm{f(\tilde{\x}_{s-1})}^2\\
    &-  2\alpha^2(\rho-1)\sum_{k=1}^m \Expectation \left[  \nabla f(\x_{s,k-1})^{\T}f(\tilde{\x}_{s-1})
\right].
\end{align}
There are several inequalities and facts that substantially simplify equation \eqref{Supp-eqn: t4-intermediate-step-after-sum}. By convexity: $(\x_{s,k-1} - \x_*)^{\T}\nabla f(\x_{s,k-1}) \geq f(\x_{s,k-1})-f(\x_*)$, and $2\lambda (f(\x_{s,k-1} - f(\x_*)) \geq \norm{\nabla f(\x_{s,k-1})}^2$. Combining these facts yields
\begin{align}\label{Supp-eqn: t4-second-intermediate-step-after-sum}
  \nonumber   \Expectation \left[\norm{\x_{s,m} - \x_*}^2 \right] \leq{}&   \Expectation \left[\norm{\x_{s,0} - \x_*}^2 \right] -2\alpha  \sum_{k=1}^m \Expectation \left[f(\x_{s,k-1})-f(\x_*) \right] \\
 \nonumber   &+  2\alpha^2 \lambda\rho \sum_{k=1}^m   \Expectation \left[f(\x_{s,k-1})-f(\x_*) \right] 
 +2\alpha^2 \lambda m(\rho-1)( f(\tilde{\x}_{s-1})-f(\x_*)) \\
    &-  2\alpha^2(\rho-1)\sum_{k=1}^m \Expectation \left[  \nabla f(\x_{s,k-1})^{\T}f(\tilde{\x}_{s-1})\right].
\end{align}
By $\gamma$-strong convexity, $\norm{\x_{s,0}-\x_*}^2 \leq 2/\gamma \left(f(\x_{s,0})-f(\x_*)\right)$ and $\x_{s,0} = \tilde{\x}_{s-1}$. We remove the term on the left-hand side because it is strictly positive, and the final term on the right-hand side because it is strictly negative. Equation \eqref{Supp-eqn: t4-second-intermediate-step-after-sum} then becomes
\begin{align}\label{Supp-eqn: t4-penultimate step}
  \nonumber  0 \leq{}&  \frac{2}{\gamma} \Expectation \left[f(\tilde{x}_{s-1})-f(\x_*) \right]  -2\alpha  \sum_{k=1}^m \Expectation \left[f(\x_{s,k-1})-f(\x_*) \right] \\
 \nonumber   &+  2\alpha^2 \lambda\rho \sum_{k=1}^m   \Expectation \left[f(\x_{s,k-1})-f(\x_*) \right] + 2\alpha^2 \lambda m(\rho-1)( f(\tilde{\x}_{s-1})-f(\x_*)).\\
\end{align}
Furthermore, by \eqref{Supp-doubleexpect}, $\Expectation \sum_{k=1}^{m} f(\x_{s,k-1})=m\Expectation f(\tilde{\x}_s)$. Therefore,
\begin{equation*}
 \Expectation\left[f(\tilde{\x}_s) - f(\x_*) \right] \leq \left( \frac{1}{\alpha \gamma m (1 - \alpha \lambda \rho)} + \frac{\alpha \lambda (\rho-1)}{1- \alpha \lambda \rho}\right)\Expectation \left[f(\tilde{\x}_{s-1})-f(\x_*) \right].
\end{equation*}
It follows that 
	\begin{equation*}
	\beta =  \frac{1}{\alpha \gamma m (1 - \alpha \lambda \rho)} + \frac{\alpha \lambda (\rho-1)}{1- \alpha \lambda \rho},
	\end{equation*}
	and by recursion we get the desired result. The proof of a.s. convergence follows as in the proof of Theorem \ref{thm:SVRG-AS}.
\qed

\subsubsection*{Proof of Theorem \ref{thm:VarianceReducedRandomGradient} part (ii)}

The proof of part (ii) of Theorem \ref{thm:VarianceReducedRandomGradient} follows closely that of part (i), but $\eta_{s,k-1}$ must be accounted for. This time, the recursion for epoch $s$ is written as 
	\begin{equation}\label{Supp-recursion}
	\x_{s,k} = \x_{s,k-1} - \alpha \v,
	\end{equation}
	with 
	\begin{equation*}
	\v  =  \P\P^{\T}\nabla f(\x_{s,k-1}) -  \eta_{s,k-1} \left(\P\P^{\T}-\I\right)\nabla f(\tilde{\x}_{s-1}).
	\end{equation*}
    By subtracting $\x_*$ from both sides of equation \eqref{Supp-recursion} and taking a norm squared we get
	\begin{equation}\label{Supp-eqn: normsquared}
	\norm{\x_{s,k}-\x_*}^2 = \norm{\x_{s,k-1}-\x_*}^2 - 2\alpha (\x_{s,k-1}-\x_*)^{\T}\v + \alpha^2 \norm{\v}^2,
	\end{equation}
which leads to
	\begin{equation}\label{Supp-eqn: expected normsquared}
\CondExp{\norm{\x_{s,k}-\x_*}^2} =\norm{\x_{s,k-1}-\x_*}^2 - 2\alpha (\x_{s,k-1}-\x_*)^{\T}\CondExp{\v} + \alpha^2 \CondExp{\norm{\v}^2}.
	\end{equation}
By definition, $\CondExp{\v} = \nabla f(\x_{s,k-1})$, so it remains to find $\CondExp{\norm{\v}^2}$. Assume first that $\eta_{s,k-1}$ is $\mathcal{F}_{s,k}$ measurable. Its optimal form based on this assumption is given by equation \eqref{Supp-eqn: eta-opt} as described next. Notice that 
\begin{align}\label{Supp-eqn: v-squared}
    \nonumber \CondExp{\norm{\v}^2} ={}& \CondExp{\rho \nabla f(\x_{s,k-1})^{\T} \P\P^{\T} \nabla f(\x_{s,k-1})} \\
 \nonumber   &+\CondExp{ (\rho-2) \eta_{s,k-1}^2 \nabla f(\tilde{\x}_{s-1})^{\T} \left(\P\P^{\T} + \frac{1}{\rho-2}\I \right) \nabla f(\tilde{\x}_{s-1})} \\
\nonumber    &- \CondExp{2 (\rho-1) \eta_{s,k-1} \nabla f(\x_{s,k-1})^{\T} \nabla f(\tilde{\x}_{s-1})} \\
    ={}& \rho \norm{\nabla f(\x_{s,k-1})}^2 + (\rho-1) \eta_{s,k-1}^2 \norm{f(\tilde{\x}_{s-1})}^2 - 2(\rho-1) \eta_{s,k-1}\nabla f(\x_{s,k-1})^{\T}f(\tilde{\x}_{s-1})
\end{align}
Plug equation \eqref{Supp-eqn: v-squared} into equation \eqref{Supp-eqn: expected normsquared} to get
	\begin{align}\label{Supp-eqn: pluggedin}
\nonumber \CondExp{\norm{\x_{s,k}-\x_*}^2} ={}&\norm{\x_{s,k-1}-\x_*}^2 - 2\alpha (\x_{s,k-1}-\x_*)^{\T}\nabla f(\x_{s,k-1}) \\
&+ \alpha^2 
\left(\rho \norm{\nabla f(\x_{s,k-1})}^2 + (\rho-1) \eta_{s,k-1}^2 \norm{f(\tilde{\x}_{s-1})}^2 - 2(\rho-1) \eta_{s,k-1}\nabla f(\x_{s,k-1})^{\T}f(\tilde{\x}_{s-1})
\right).
	\end{align}
Equation \eqref{Supp-eqn: pluggedin} is quadratic with respect to $\eta_{s,k-1}$ so by taking a derivative with respect to $\eta_{s,k-1}$ it is possible to find the value which minimizes the expected conditional  norm squared error:
\begin{equation}\label{Supp-eqn: intermediate-eta-opt}
    0 = 2\alpha^2(\rho-1)\eta_{s,k-1} \norm{f(\tilde{\x}_{s-1})}^2 - 2\alpha^2(\rho-1)\nabla f(\x_{s,k-1})^{\T}f(\tilde{\x}_{s-1}).
\end{equation}
It is simple to see from equation \eqref{Supp-eqn: intermediate-eta-opt} that 
\begin{equation}\label{Supp-eqn: eta-opt}
    \eta_{s,k-1} = \frac{\nabla f(\x_{s,k-1})^{\T}\nabla f(\tilde{\x}_{s-1})}{\norm{f(\tilde{\x}_{s-1})}^2}
\end{equation}
is the optimal $\eta_{s,k-1}$. Plug equation \eqref{Supp-eqn: eta-opt} back into equation \eqref{Supp-eqn: pluggedin} and notice there are several term cancellations
\begin{align}\label{Supp-eqn: intermediate-step-before-sum}
    \nonumber \CondExp{\norm{\x_{s,k}-\x_*}^2} ={}&\norm{\x_{s,k-1}-\x_*}^2 - 2\alpha (\x_{s,k-1}-\x_*)^{\T}\nabla f(\x_{s,k-1}) \\
    &+ \alpha^2 \rho \norm{\nabla f(\x_{s,k-1}}^2 - \alpha^2(\rho-1) \left( \frac{\nabla f(\x_{s,k-1})^{\T}\nabla f(\tilde{\x}_{s-1})}{\norm{f(\tilde{\x}_{s-1})}} \right)^2.
\end{align}
We take an expectation of both sides of the equality in equation \eqref{Supp-eqn: intermediate-step-before-sum}, and since the interest is in iterations of $s$ rather than $k$ we sum from $k=1$ to $k=m$
\begin{align}\label{Supp-eqn: intermediate-step-after-sum}
  \nonumber  \sum_{k=1}^m \Expectation \left[\norm{\x_{s,k} - \x_*}^2 \right] ={}&  \sum_{k=1}^m \Expectation \left[\norm{\x_{s,k-1} - \x_*}^2 \right] -2\alpha  \sum_{k=1}^m \Expectation \left[(\x_{s,k-1}-\x_*)^{\T}\nabla f(\x_{s,k-1}) \right] \\
    &+  \alpha^2 \rho \sum_{k=1}^m \Expectation\left[\norm{\nabla f(\x_{s,k-1})}^2\right] - \alpha^2(\rho-1) \sum_{k=1}^m \Expectation\left[\left( \frac{\nabla f(\x_{s,k-1})^{\T}\nabla f(\tilde{\x}_{s-1})}{\norm{f(\tilde{\x}_{s-1})}} \right)^2\right].
\end{align}
Consider that the series on the left-hand side minus the first series on the right forms a telescoping series. Recall also that by convexity, $(\x_{s,k-1} - \x_*)^{\T}\nabla f(\x_{s,k-1}) \geq f(\x_{s,k-1})-f(\x_*)$. Also by convexity, $2\lambda (f(\x_{s,k-1} - f(\x_*)) \geq \norm{\nabla f(\x_{s,k-1})}^2$. Combining these facts yields
\begin{align}\label{Supp-eqn: second-intermediate-step-after-sum}
  \nonumber   \Expectation \left[\norm{\x_{s,m} - \x_*}^2 \right] \leq{}&   \Expectation \left[\norm{\x_{s,0} - \x_*}^2 \right] -2\alpha  \sum_{k=1}^m \Expectation \left[f(\x_{s,k-1})-f(\x_*) \right] \\
 \nonumber   &+  2\alpha^2 \lambda\rho \sum_{k=1}^m   \Expectation \left[f(\x_{s,k-1})-f(\x_*) \right] \\
    &- \alpha^2(\rho-1) \sum_{k=1}^m \Expectation\left[\left( \frac{\nabla f(\x_{s,k-1})^{\T}\nabla f(\tilde{\x}_{s-1})}{\norm{f(\tilde{\x}_{s-1})}} \right)^2\right].
\end{align}
By $\gamma$-strong convexity, $\norm{\x_{s,0}-\x_*}^2 \leq 2/\gamma \left(f(\x_{s,0})-f(\x_*)\right)$. We remove the term on the left-hand side because it is strictly positive, and the final term on the right-hand side because it is strictly negative. Equation \eqref{Supp-eqn: second-intermediate-step-after-sum} becomes
\begin{align}\label{Supp-eqn: penultimate step}
  \nonumber  0 \leq{}&  \frac{2}{\gamma} \Expectation \left[f(\x_{s,0})-f(\x_*) \right]  -2\alpha  \sum_{k=1}^m \Expectation \left[f(\x_{s,k-1})-f(\x_*) \right] \\
 \nonumber   &+  2\alpha^2 \lambda\rho \sum_{k=1}^m   \Expectation \left[f(\x_{s,k-1})-f(\x_*) \right]. 
\end{align}
Further simplification results from equation \eqref{Supp-doubleexpect} which suggests that \\ $\Expectation \sum_{k=1}^{m} f(\x_{s,k-1})=m\Expectation f(\tilde{\x}_s)$. We move this term to the left-hand side so that
\begin{equation*}
 \Expectation\left[f(\tilde{\x}_s) - f(\x_*) \right] \leq \frac{1}{\alpha \gamma m (1 - \alpha \lambda \rho)}\Expectation \left[f(\tilde{\x}_{s-1})-f(\x_*) \right].
\end{equation*}
It follows that 
	\begin{equation*}
	\beta = \frac{1}{\alpha \gamma m (1 - \alpha \lambda \rho)},
	\end{equation*}
	and by recursion we get the desired result. The proof of almost sure convergence follows in the same manner as in the proof of Theorem \ref{thm:SVRG-AS}.

\subsection*{S3. Miscellaneous Results}

\subsubsection*{Haar matrix construction}
The pseudocode provided in Algorithm \ref{Haar} allows to draw from a scaled, Haar distributed matrix.

\begin{algorithm}[H]
	\caption{Generate a scaled, Haar distributed matrix (based on \cite{mezzadri2006generate})}\label{Haar}
	\begin{algorithmic}
		\Inputs{$\ell,d$ \Comment{Dimensions of desired matrix, $d > \ell$}}
		\Outputs{$\textbf{P}\in \reals^{d\times \ell}$ such that: $ \textbf{P}^{\T}\textbf{P} =\frac{d}{\ell} I_{\ell}$, 
			$\Expectation~\textbf{P}\textbf{P}^{\T} = I_d$,
			 $\textbf{P}$ has orthogonal columns}
		\State \textbf{Initialize} $X \in \reals^{d \times \ell}$
		\State \textbf{Set} $X_{i,j} \sim \mathcal{N}(0,1)$
		\State \textbf{Calculate} thin QR decomposition of $\X=\textbf{QR}$
		\State \textbf{Let} $\boldsymbol{\Lambda} = \begin{pmatrix} \mathrm{R}_{1,1}/\abs{\mathrm{R}_{1,1}} && \\ & \ddots & \\ && \mathrm{R}_{\ell,\ell}/\abs{\mathrm{R}_{\ell,\ell}}  \end{pmatrix}$
		\State $\textbf{P} =\sqrt{\frac{d}{\ell}}\mathbf{Q}\boldsymbol{\Lambda}$

	\end{algorithmic}
\end{algorithm}
\subsubsection*{Variance of estimators} \label{sec:P}
Consider drawing a random matrix $\bP \in \R^{d \times \ell}$, for $\ell=1$, in one of two ways: the first (i.e., random coordinate descent), by selecting it to be a scaled canonical basis vector $\bP=\sqrt{d}\textbf{e}_i$ where $i$ is uniform in $\{1,\ldots,d\}$, 
and the second (i.e., random subspace), by selecting $\bP$ as an appropriate scaling of the first column of a random matrix from the Haar distribution on orthogonal $d \times d$ matrices.
In both cases it is elementary to show $\Expectation (\bP\bP^\T) = I_{d}$, as well as $\bP^\T\bP = d$, so both satisfy the prerequisites for our theorem. Which type of random matrix is better?  For simplicity in the present discussion we replace the Haar distribution in the second scheme with a Gaussian distribution; that is, one in which the entries of $\mathbf{P}$ are independent and identically distributed $\mathcal{N}(0,1)$. The scaled Haar distribution, for $\ell=1$, is the same as a standard normal Gaussian vector that has been properly normalized. 

Fix $\x = \textbf{e}_1$. It is straightforward to show that $\mathbb{V}\mathrm{ar}(\|\bP^T \textbf{e}_1\|^2) \to 2$ as $d \to \infty$ when $\mathbf{P}$ is generated by Algorithm \ref{Haar}. In the ``coordinate descent'' case $\mathbb{V}\mathrm{ar}(\|\bP^T \textbf{e}_1\|^2)=d-1$, whereas  $\mathbb{V}\mathrm{ar}(\|\bP^T \textbf{e}_1\|^2)=2$ in the ``random subspace'' case.
That is, the random subspace approach has much lower variance in high-dimensional settings, hence it has much better deviation bounds.
Note that in the random subspace case, the variance is independent of the particular $\x$ chosen, due to the orthogonal invariance of the Gaussian distribution (which also holds for the Haar distribution), so there is no ``worst-case'' vector for this scheme, whereas our choice $\x = \textbf{e}_1$ is an adversarial choice for the ``random coordinate'' scheme.

\subsubsection*{The control variate parameter}
In order to understand the role of $\eta$ more clearly, we consider three possible scenarios:  $\eta = 0$ corresponds to standard SSD as outlined in Theorem \ref{thm:convergence}, $\eta = 1$ corresponds to part ($i$) of Theorem \ref{thm:VarianceReducedRandomGradient}, and finally the optimal $\eta$ is given by \eqref{eqn: eta-opt} and used in part ($ii$) of Theorem \ref{thm:VarianceReducedRandomGradient}. With the understanding that the true gradient is the optimal search direction with respect to one-step improvement, we would like our search direction to vary around it as little as possible. To this end, we ascertain the conditional mean-square error defined as the expectation of
\begin{equation*}
    \norm{\nabla f(\x_{s,k-1}) - \left(\P\P^{\T}\nabla f(\x_{s,k-1}) -\eta_{s,k-1} \left(\P\P^{\T}-\I\right)\nabla f(\tilde{\x}_{s-1})\right)}^2,
\end{equation*}
given information at the current iteration (i.e., conditional on $\sigma-$algebras defined in the appendix). We denote this conditional error as CMSE. 
\begin{itemize}
    \item Case $\eta = 0$ (Theorem \ref{thm:convergence})
   \begin{equation*}  \text{CMSE}= (\rho-1) \norm{\nabla f(\x_{s,k-1})}^2
   \end{equation*}
    \item Case $\eta = 1$ (Theorem \ref{thm:VarianceReducedRandomGradient} part (i))
       \begin{equation*}  \text{CMSE}=  (\rho-1)\left(\norm{\nabla f(\x_{s,k-1})}^2 + \norm{\nabla f(\tilde{\x}_{s-1})}^2 - 2\ \nabla f(\x_{s,k-1})^{\T}\nabla f(\tilde{\x}_{s-1}) \right).
   \end{equation*}
    \item Case $\eta = \nabla f(\tilde{\x}_{s-1})^{\T} \nabla f(\x_{s,k-1}) / \norm{\nabla f(\tilde{\x}_{s-1})}^2$ (Theorem \ref{thm:VarianceReducedRandomGradient} part (ii))
    \begin{align}\label{eqn: eta-analysis}
        \text{CMSE}&=(\rho-1)\left(\norm{\nabla f(\x_{s,k-1})}^2- \left(\frac{\nabla f(\tilde{\x}_{s-1})^{\T} \nabla f(\x_{s,k-1})}{\norm{\nabla f(\tilde{\x}_{s-1})}}\right)^2 \right).
    \end{align}
\end{itemize}
These identities follow from the proofs in the appendix. Equation \eqref{eqn: eta-analysis} shows that the optimal $\eta$ never increases the variance compared to SSD and can in fact drop it to zero, though an increase in the variance is possible if $\eta=1$ (e.g., in the event $\nabla f(\x_{s,k-1})^{\T}\nabla f(\tilde{\x}_{s-1})$ is negative, which may occur if $m$ is too large and the iterates traverse substantial portions of parameter space between full gradient calculations). These results are in agreement with Theorem \ref{thm:VarianceReducedRandomGradient} which suggests that the per-iteration rate of convergence using an optimal $\eta$ is strictly less than the rate when $\eta=1$.

\end{document}